\documentclass[12pt]{article}
\usepackage{amssymb, amsmath}
\begin{document}
\renewcommand{\theequation}{\arabic{section}.\arabic{equation}}
\newtheorem{theorem}{Theorem}[section]
\newtheorem{definition}{Definition}[section]
\newtheorem{lemma}{Lemma}[section]
\newtheorem{pro}{Proposition}[section]
\newtheorem{cor}{Corollary}[section]
\newcommand{\n}{\nonumber}
\newcommand{\tv}{\tilde{v}}
\renewcommand{\l}{\lambda}
\newcommand{\lo}{\lambda_1}
\newcommand{\lt}{\lambda_2}
\newcommand{\ltt}{\lambda_3}
\newcommand{\s}{\sigma}
\newcommand{\tw}{\tilde{\omega}}
\renewcommand{\t}{\theta}
\renewcommand{\th}{\theta}
\newcommand{\w}{\omega}
\renewcommand{\b}{\dot{B}^0_{\infty ,1}}
\newcommand{\e}{\varepsilon}
\renewcommand{\a}{\alpha}
\renewcommand{\l}{\lambda}
\newcommand{\vare}{\varepsilon}
\renewcommand{\o}{\omega}
\renewcommand{\O}{\Omega}
\newcommand{\Om}{\Omega}
\newcommand{\nao}{\nabla^\bot}
\newcommand{\thee}{\theta^\vare}
\newcommand{\bn}{\|_{B^{\alpha, \infty}_{\frac{9}{5}}}}
\newcommand{\pn}{\|_{L^{\frac{9}{5}}}}
\newcommand{\pnn}{\|_{L^{\frac{9}{2}}}}
\newcommand{\ve}{v^{\varepsilon}}
\renewcommand{\oe}{\omega^{\varepsilon}}
\renewcommand{\a}{\alpha}
\newcommand{\intr}{\int_{\Bbb R^n}}
\newcommand{\intrr}{\int_{\Bbb T^3}}
\newcommand{\intt}{\int_{\Bbb R^3}}
\newcommand{\intom}{\int_{\Omega}}
\newcommand{\bb}{\begin{equation}}
\newcommand{\ee}{\end{equation}}
\newcommand{\bq}{\begin{eqnarray}}
\newcommand{\eq}{\end{eqnarray}}
\newcommand{\bqn}{\begin{eqnarray*}}
\newcommand{\eqn}{\end{eqnarray*}}
\title{Incompressible Euler Equations: the blow-up problem and  related results}
\author{Dongho Chae\thanks{This research was supported partially
by  KRF Grant(MOEHRD, Basic Research Promotion Fund). Part of the
work was done, while the author was visiting RIMS, Kyoto University.
He would like to  thank to Professor Hisashi Okamoto
for his hospitality during the visit.} \\
Department of Mathematics\\
              Sungkyunkwan University\\
              Suwon 440-746, Korea\\
  e-mail: {\it chae@skku.edu}}
 \date{}
  \maketitle
\begin{abstract}
The question of spontaneous apparition of singularity in the 3D
incompressible Euler equations is one of the most important and
challenging open problems in mathematical fluid mechanics.  In this
survey article we review some of recent approaches to  the problem.
We first review Kato's classical local well-posedness result in the
Sobolev space and derive the celebrated Beale-Kato-Majda criterion
for finite time blow-up. Then, we discuss
 recent refinements of the criterion as well as geometric type
of theorems on the sufficiency condition for the regularity of
solutions. After that we review
 results excluding some of the scenarios leading to finite time
 singularities. We also survey studies of various simplified
 model problems.
 A dichotomy type of result between the finite time blow-up and the
 global in time regular dynamics is presented, and a spectral
 dynamics approach to study local in time behaviors of the
 enstrophy is also reviewed. Finally, progresses on the problem of
 optimal regularity for solutions to have  conserved
 quantities are presented.
\end{abstract}
\tableofcontents
\section{Introduction}
 \setcounter{equation}{0}

The motion of homogeneous incompressible ideal fluid in a domain
$\Omega \subset \Bbb R^n $ is described by the following system of
Euler equations.
 \[
\mathrm{ (E)}
 \left\{ \aligned
 &\frac{\partial v}{\partial t} +(v\cdot \nabla )v =-\nabla p,
 \quad (x,t)\in {\Omega}\times (0, \infty) \\
 &\textrm{div }\, v =0 , \quad (x,t)\in {\Omega}\times (0,
 \infty)\\
  &v(x,0)=v_0 (x), \quad x\in \Omega
  \endaligned
  \right.
  \]
 where $v=(v^1, v^2, \cdots , v^n )$, $v^j
=v^j (x, t)$, $j=1, 2, \cdots , n,$ is the velocity of the fluid
flows, $p=p(x,t)$ is the scalar pressure, and $v_0 (x)$ is a given
initial velocity field satisfying div $v_0=0$.  Here we use the
standard notion of vector calculus, denoting
$$\nabla p =\left(\frac{\partial p}{\partial x_1 },
\frac{\partial p}{\partial x_2 }, \cdots , \frac{\partial
p}{\partial x_n } \right),\quad  (v\cdot \nabla )v^j=\sum_{k=1}^n
v^k \frac{\partial v^j}{\partial x_k},\quad \textrm{div }\, v
=\sum_{k=1}^n \frac{\partial v^k}{\partial x_k }.
$$
 The first equation of (E)  follows from the balance of momentum for
each portion of fluid, while the second equation can be derived from
the conservation of mass of fluid during its motion, combined with
the homogeneity(constant density) assumption on the fluid. The
system (E) is first derived by L. Euler in 1755(\cite{eul}). Unless
otherwise stated,  we are concerned on the Cauchy problem of the
system (E) on $\Omega= \Bbb R^n$, but many of the results presented
here are obviously valid also for
 $\Omega =\Bbb R^n/\Bbb Z^n$(periodic domain), and even for the bounded domain
 with the smooth boundary with the boundary condition $v\cdot \nu
 =0$, where $\nu $ is the outward unit normal vector.
 We also suppose $n=2$ or $3$ throughout this paper.
{\em In this article our aim to survey recent results on the
 mathematical aspects the 3D Euler equations closely related to
the problem of spontaneous apparition of singularity starting from a
  classical solutions having finite energy.}
If we add the dissipation term $\mu \Delta v=\mu \sum_{j=1}^n
 \frac{\partial^2 v}{\partial x_j^2} $, where $\mu>0$ is the
 viscosity coefficient, to the right hand side of the first equation of (E), then we have
  the Navier-Stokes equations, the regularity/singularity question of which is
   one of the seven millennium problems in mathematics.
 In this article we do not  treat the
 Navier-Stokes equations. For details of mathematical
 studies on the  Navier-Stokes equations
 see e.g. \cite{tem3, con10,lio, gald, lad, maj3, lem}.
We also omit other important topics such as existence and uniqueness
   questions of the weak solutions of the 2D Euler equations, and
   the related vortex patch problems, vortex sheet problems, and so
   on. These are well treated in the other
   papers and monographs(\cite{maj3, cha19a, che2, lio, sch, shn1, yud1, yud2, vis2, tad})
   and the references therein.
   For the survey related the stability
   question please see for example \cite{fri1} and references therein. For the
   results on  the regularity of the Euler equations with uniformly rotating
   external force we refer \cite{bab}, while for the numerical
   studies on the blow-up problem of the Euler equations there are
   many articles including \cite{ker1, ker2, hou3, bra, fri2, caf1, gra1,gra2,
   gre, pel}. For various mathematical  and physical aspects of the Euler equations
   there are many excellent books, review articles including
   \cite{arn, bren, che2, con1,con3, fri1, gib2, maj2, maj3, mar, cha12a, yu}.
   Obviously, the  references are not complete mainly due to
    author's ignorance.

\subsection{Basic properties}

In the study of the Euler equations the notion of vorticity,
$\omega=$curl $v$, plays important roles.  We can reformulate the
Euler system in terms of the vorticity fields only as follows. We
first consider the 3D case. Let us
 first rewrite the first equation of (E) as
 \bb\label{bas1}
\frac{\partial v}{\partial t} -v\times \mathrm{curl}\, v=-\nabla
(p+\frac12 |v|^2 ).
 \ee
Then, taking curl of (\ref{bas1}), and using elementary vector
identities, we obtain the following vorticity formulation:
  \bb\label{bas2}
 \frac{\partial \o }{\partial t} +(v\cdot \nabla )\o =\o \cdot \nabla
 v,
 \ee
 \bb\label{bas3}
  \textrm{div }\, v =0 , \quad  \textrm{curl }\, v =\o ,
 \ee\label{bas4}
 \bb \o(x,0)=\o_0 (x).
 \ee
 The linear elliptic system (\ref{bas3}) for $v$
  can be solved explicitly in terms of $\o$, assuming $\o$ decays sufficiently fast
  near spatial infinity, to provides us with the Biot-Savart law,
 \bb\label{bas5}
 v(x,t)=\frac{1}{4\pi}\int_{\Bbb R^3}
 \frac{(x-y)\times \o (y,t) }{|x-y|^3}dy.
 \ee
 Substituting this $v$ into (\ref{bas2}), we obtain an
 integro-differential system for $\o$. The term in the right hand side of
 (\ref{bas2}) is called the vortex stretching term, and is
 regarded as the main source of difficulties in the mathematical
 theory of the 3D Euler equations.
Let us introduce the deformation matrix
$S(x,t)=(S_{ij}(x,t))_{i,j=1}^3$ defined  as the symmetric part of
the velocity gradient matrix,
$$S_{ij} =\frac12 \left(\frac{\partial v_j}{\partial x_{i}}
 +\frac{\partial v_i}{\partial x_{j}}\right).
 $$
 From the Biot-Savart law in (\ref{bas5}) we can explicitly compute
 \bb\label{bas5a}
 S(x,t)=\frac{3}{8\pi} p.v. \int_{\Bbb R^3} \frac{ [ (y\times \o (x+y,t) )
 \otimes y +y\otimes
 (y\times
 \o (x+y,t))]}{ |y|^5}dy
 \ee
 (see e.g. \cite{maj3} for the details on the computation).
 The kernel in the convolution integral of (\ref{bas5a}) defines a
 singular integral operator of the Calderon-Zygmund type(see e.g.
 \cite{ste1, ste2} for more details).
 Since the vortex stretching term can be written as $ (\o \cdot \nabla )v = S \o
 $, we see that the singular integral operator and related
 harmonic analysis results could have  important roles to
study the Euler equations.

In the two dimensional case
 we take the vorticity  as the scalar, $\o =\frac{\partial v^2}{\partial x_1}
 -\frac{\partial v^1}{\partial x_2}$, and the evolution equation
 of $\o$ becomes
 \bb\label{bas6}
 \frac{\partial \o }{\partial t} +(v\cdot \nabla )\o =0,
 \ee
where the velocity is represented in terms of the vorticity by the
2D Biot-Savart law,
 \bb\label{bas7}
  v(x,t)=\frac{1}{2\pi} \int_{\Bbb R^2} \frac{(-y_2+x_2,
 y_1 -x_1
 )}{|x-y|^2} \o (y,t )dy.
 \ee
 Observe that there is no  vortex stretching term in (\ref{bas6}),
 which makes the proof of global regularity in 2D Euler equations easily accessible.
In many studies of the Euler equations it is convenient to introduce
the notion of `particle trajectory mapping', $X (\cdot ,t)$ defined
by
 \bb\label{bas8}
 \frac{\partial X (a, t )}{\partial t} =v (X (a , t
 ),t), \quad X(a , 0)=a , \quad a \in \Omega.
 \ee
 The mapping $X(\cdot, t)$ transforms from the location of the initial fluid particle
 to the location at time $t$, and the parameter $a$ is called the Lagrangian
 particle marker.
 If we denote the Jacobian of the transformation, det$(\nabla
 _a X (a ,t ))=J(a ,t)$, then we can show easily(see e.g.
 \cite{maj3} for the proof)
 that
 $$\frac{\partial J}{\partial t} =(\mathrm{div}\, v) J,$$
which implies that  the velocity field $v$ satisfies the
incompressibility, div $v=0$ if and only if the mapping $X (\cdot
,t)$ is volume preserving.  At this moment we note that, although
the Euler equations are originally derived by applying the physical
principles of mass conservation and  the momentum balance, we could
also derive them  by applying the least action principle to the
action defined by
$$
\mathcal{I} (A)=\frac12 \int_{t_1}^{t_2} \int_{\Omega}
\left|\frac{\partial X (x,t)}{\partial t}\right|^2 dxdt.
$$
Here,  $X(\cdot, t): \Omega \to \Omega\subset \Bbb R^n$ is a
parameterized family of volume preserving diffeomorphism. This
variational approach to the Euler equations implies that we can view
solutions of the Euler equations as a geodesic curve in the
$L^2(\Omega)$ metric on the infinite dimensional manifold of volume
preserving diffeomorphisms(see  e.g. \cite{arn, bren, ebi} and
references therein for more details on the geometric approaches to
the Euler equations).

The 3D Euler equations have many conserved quantities. We list some
important ones below.
\begin{itemize}
\item[(i)] Energy,
 $$
 E(t)=\frac12 \int_{\Omega} |v(x,t)|^2 dx.
$$
\item[(ii)] Helicity,
 $$
 H(t)=\int_{\Omega} v (x,t)\cdot \o (x,t)
dx.
 $$
\item[(iii)] Circulation,
$$
  \Gamma _{\mathcal{C}(t)}= \oint_{\mathcal{C}(t)} v\cdot
dl,
$$
 where $\mathcal{C}(t)=\{ X (a , t) \,| \, a \in \mathcal{C}\}$ is
a curve moving along with the fluid.
\item[(iv)] Impulse,
 $$
 I(t)=\frac12\int_{\Omega} x\times \o \,dx.
 $$
 \item[(v)] Moment of Impulse,
 $$
  M(t)= \frac13 \int_{\Omega} x\times (x\times \o ) \,dx.
  $$
 \end{itemize}
 The proof of conservations of the above quantities for the classical solutions
 can be done without difficulty using elementary vector calculus(for details
 see e.g. \cite{maj3, mar}). The helicity, in particular,
 represents the degree of knotedness of the vortex lines in the
 fluid, where the vortex lines are the integral curves of the
 vorticity fields. In \cite{arn} there are detailed discussions on this aspects
  and other topological implications of the helicity conservation.
 For the 2D Euler equations there is no analogue
 of helicity, while the circulation conservation is replaced by
 the vorticity flux integral,
$$
 \int_{D(t)} \o (x,t) dx,
 $$
 where $D(t)=\{ X (a , t) \, | \, a \in
D\subset \Omega\}$  is a planar region moving along the fluid in
$\Omega$.
 The impulse and the moment of impulse integrals in the 2E Euler equations are replace by
  $$
  \frac12\int_{\Omega} (x_2, -x_1 )\o dx \quad \mbox{and}\quad
 -\frac13 \int_{\Omega} |x|^2 \o dx\quad \mbox{respectively}.
 $$
 In the 2D Euler equations we have extra conserved
 quantities; namely for any continuous function $f$ the integral
 $$
  \int_\Omega f(\o (x,t)) dx
$$
 is conserved.
There are also many known explicit solutions to the Euler equations,
for which we just refer \cite{lam, maj3}. In the remained part of
this subsection we introduce some notations to be used later for 3D
Euler equations. Given velocity $v(x,t)$, and pressure $p(x,t)$, we
set the $3\times 3$ matrices,
$$
V_{ij}=\frac{\partial v_j}{\partial x_i},\quad
S_{ij}=\frac{V_{ij}+V_{ji}}{2},\quad A_{ij}=\frac{V_{ij}-V_{ji}}{2},
\quad P_{ij}=\frac{\partial ^2 p}{\partial x_i \partial x_j},
$$
with  $i,j=1,2,3$. We have the decomposition $V=(V_{ij})=S+A$, where
the symmetric part $S=(S_{ij})$ represents the deformation tensor of
the fluid introduced above, while the antisymmetric part
$A=(A_{ij})$ is related to the vorticity $\o$
 by the formula,
 \bb\label{basfor}
 A_{ij}= \frac12 \sum_{k=1}^3 \e
_{ijk} \o_k,\qquad \o_i = \sum_{j,k=1}^3\e_{ijk}A_{jk},
 \ee
 where $\e_{ijk}$ is the skewsymmetric tensor with the normalization
 $\e_{123}=1$. Note that $P=(P_{ij})$ is the hessian of the
 pressure. We also frequently use the
 notation for the vorticity direction field,
 $$ \xi (x,t)=\frac{\o (x,t)}{|\o(x,t)|}, $$
 defined whenever $\o (x,t)\neq 0$.
Computing partial derivatives $\partial/
\partial x_k$ of the first equation of (E), we obtain the matrix equation
 \bb\label{basfor1}
\frac{D V}{Dt}=-V^2 -P, \quad \frac{D}{Dt}=
 \frac{\partial}{\partial t}+(v\cdot \nabla )v.
\ee Taking symmetric part of this, we obtain
$$ \frac{D
S}{Dt}=-S^2-A^2-P,
 $$
 from which, using  the formula (\ref{basfor}),
 we have
 \bb\label{basfor2}
\frac{D S_{ij}}{Dt} = - \sum_{k=1}^3 S_{ik}S_{kj}
 +\frac14 (|\o|^2 \delta_{ij} -\o_i\o_j )- P_{ij},
 \ee
 where $\delta_{ij}=1$ if $i=j$, and $\delta_{ij}=0$ if $i\neq j$.
  The antisymmetric part of (\ref{basfor1}), on the other hand,  is
 $$
 \frac{DA}{Dt}=-SA -AS,
 $$
 which, using the formula (\ref{basfor}) again, we obtain easily
 \bb\label{basfor2a}
 \frac{D \o}{Dt} = S \o,
 \ee
 which is the  vorticity evolution equation (\ref{bas2}). Taking dot
 product (\ref{basfor2a}) with $\o$,
 we immediately have
 \bb\label{basfor3}
  \frac{D |\o |}{Dt} =\a |\o|,
  \ee
  where we set
  $$\a (x,t)= \left\{ \aligned &\sum_{i,j=1}^3
   \xi_i(x,t) S_{ij}(x,t) \xi_j (x,t)
  & \mbox{if $\o (x,t)\neq 0$}\\
    &0 & \mbox{if $\o (x,t)=0$}.
    \endaligned
    \right. $$
\subsection{Preliminaries}

Here we introduce some notations and function spaces to be used in
the later sections.
 Given $p\in [1, \infty]$, the Lebesgue space $L^p (\Bbb R^n )$,
  $p\in [1, \infty]$, is the Banach space defined
 by the norm
$$\|f\|_{L^p}:=\left\{
\aligned &\left(\intr |f(x)|^p dx \right)^{\frac{1}{p}}, \quad
 p\in [1, \infty )\\
&\mathrm{ess.} \sup_{x\in \Bbb R^n} |f(x)|,  \quad p=\infty .
\endaligned \right.
$$
For $j=1, \cdots, n$ the Riesz transform $R_j$ of $f$ is given by
$$
R_j (f)(x)=\frac{\Gamma (\frac{n+1}{2} )}{\pi ^{\frac{n+1}{2}}} p.v.
\int_{\Bbb R^n} \frac{x_j -y_j}{|x-y|^{n+1}} f(y)dy
$$
whenever the right hand side makes sense. The Hardy space
$\mathcal{H}^1 (\Bbb R^n)\subset L^1 (\Bbb R^n)$ is defined by
$$
f\in \mathcal{H}^1 (\Bbb R^n) \quad \mbox{if and only if}\quad
\|f\|_{\mathcal{H}^1}:= \|f\|_{L^1} +\sum_{j=1}^n \|R_j f\|_{L^1}
<\infty.
$$
 The space $BMO(\Bbb R^n)$ denotes the space of functions of
bounded mean oscillations, defined by
$$
f\in BMO (\Bbb R^n )\quad \mbox{if and only if}\quad
\|f\|_{BMO}:=\sup_{Q\subset \Bbb R^n} \frac{1}{\mathrm{Vol} (Q)}
\int_Q |f-f_Q |dx <\infty,
$$
where $f_Q =\frac{1}{\mathrm{Vol} (Q)}\int_Q f dx$. For more details
on the Hardy space and BMO we refer \cite{ste1, ste2}.
 Let us set  the
multi-index $\alpha :=(\alpha_1 , \alpha_2 , \cdots , \alpha_n )\in
(\Bbb
 Z_+ \cup \{ 0\} )^n $ with $|\alpha |=\alpha_1 +\alpha_2 +\cdots
 +\alpha_n$. Then, $D^\alpha :=D^{\alpha_1}_1 D^{\alpha_2}_2 \cdots
 D^{\alpha_n}_n$, where $D_j =\partial/\partial x_j$, $j=1,2,\cdots,
 n$.
Given $k\in \Bbb Z$ and  $p\in [1, \infty)$ the Sobolev space,
 $W^{k,p} (\Bbb R^n )$ is the Banach space of functions
  consisting of functions $f\in L^p (\Bbb R^n )$
 such that
 $$\|f\|_{W^{k, p}}:=\left(\intr |D^\alpha f(x)|^p dx
 \right)^{\frac{1}{p}} <\infty,$$
  where the derivatives are in the sense of distributions. For
  $p=\infty$ we replace the $L^p (\Bbb R^n )$ norm by the
  $L^\infty (\Bbb R^n)$ norm. In particular, we denote $H^m (\Bbb R^n) =W^{m,2} (\Bbb R^n )$.
 In order to handle the functions having fractional derivatives of order $s\in \Bbb R$,
  we use the
 Bessel potential space $L^{s}_p (\Bbb R^n) $ defined by the Banach spaces norm,
 $$\|f\|_{L^{s,p}}:=\|(1-\Delta )^{\frac{s}{2}} f\|_{L^p},
 $$
  where $(1-\Delta )^{\frac{s}{2}} f=\mathcal{F}^{-1} [ (1+|\xi |^2
  )^{\frac{s}{2}} \mathcal{F} (f)(\xi )]$. Here  $\mathcal{F} (\cdot)$ and
  $\mathcal{F}^{-1}(\cdot)$ denoting the Fourier transform and its
  inverse, defined by
$$
 \mathcal{F} (f)(\xi)=\hat{f} (\xi)
 =\frac{1}{(2\pi )^{n/2}}\int_{\Bbb R^n} e^{-ix\cdot \xi }
 f(x)dx,
  $$
 and
$$
 \mathcal{F}^{-1} (f)(x)=\check{f} (x)
 =\frac{1}{(2\pi )^{n/2}}\int_{\Bbb R^n} e^{ix\cdot \xi }
 f(\xi)d\xi ,
  $$
  whenever the integrals make sense.
 Next we
introduce the Besov spaces. We follow \cite{tri}(see also \cite{tay,
lem, che2, run}). Let $\mathfrak{S}$ be the Schwartz class of
rapidly decreasing
 functions.
We consider  $\varphi \in \mathfrak{S}$ satisfying
 $\textrm{Supp}\, \hat{\varphi} \subset
 \{\xi \in {\mathbb R}^n \, |\,\frac12 \leq |\xi|\leq
 2\}$,
 and $\hat{\varphi} (\xi)>0 $ if $\frac12 <|\xi|<2$.
 Setting $\hat{\varphi_j } =\hat{\varphi } (2^{-j} \xi )$ (In other words,
 $\varphi_j (x)=2^{jn} \varphi (2^j x )$.), we can adjust the
 normalization constant in front of $\hat{\varphi}$ so that
  $$
 \sum_{j\in \mathbb{Z}}  \hat{\varphi}_j (\xi )=1\quad \forall \xi \in
 {\mathbb R^n}\setminus \{ 0\}.
 $$

 Let $s\in \mathbb R$, $p,q
 \in [0, \infty]$. Given $f\in \mathfrak{S}'$, we denote
 $\Delta_j f=\varphi_j* f$.
 Then the homogeneous Besov semi-norm $\|
f\|_{\dot{B}^{s}_{p,q}}$ is defined by \[ \|
f\|_{\dot{B}^s_{p,q}}=\left\{ \aligned &\left( \sum_{j\in \Bbb Z}
2^{jqs} \| \varphi_j \ast f \|_{L^p}^q
\right)^{\frac{1}{q}}&\mbox{ if } q\in [1, \infty) \\
&\sup_{j\in \Bbb Z} \left(2^{js} \| \varphi_j \ast f \|_{L^p}
\right)&\mbox{ if }q=\infty.
\endaligned \right.
\]
For $(s, p,q)\in [0, \infty) \times [1, \infty]\times [1, \infty]$
the homogeneous Besov space $\dot{B}^{s}_{p,q}$ is a quasi-normed
space with the quasi-norm given by $\| \cdot \|_{\dot{B}^s_{p,q}}$.
For $s>0$ we define the inhomogeneous Besov space norm $\|
f\|_{{B}^s_{p,q}}$ of $f\in \mathfrak{S}'$ as
$\|f\|_{{B}^s_{p,q}}=\| f\|_{L^p}+\| f\|_{\dot{B}^s_{p,q}}$.
Similarly, for $(s, p,q)\in [0, \infty) \times [1, \infty)\times [1,
\infty]$, the homogeneous Triebel-Lizorkin semi-norm
$\|f\|_{\dot{F}^s_{p,q}}$ is defined by
 $$
\|f\|_{\dot{F}^s_{p,q}} =\left\{ \aligned &\left\| \left(\sum_{j\in
\mathbb{Z}} 2^{jqs}|\varphi_j \ast
f(\cdot)|^q\right)^{\frac1q}\right\|_{L^p}
& \mbox{if $q\in[1, \infty)$}\\
&\left\|\sup_{j\in \mathbb{Z}}\left( 2^{js} |\varphi_j \ast f (\cdot
)|\right) \right\|_{L^p} & \mbox{if $q=\infty$}
\endaligned \right. .
$$
 The homogeneous Triebel-Lizorkin space $\dot{F}^s_{p,q}$ is a quasi-normed
 space with the quasi-norm given by $\|\cdot \|_{\dot{F}^s_{p,q}}$.
 For $s>0$, $(p,q)\in [1, \infty)\times [1, \infty)$
 we define the inhomogeneous Triebel-Lizorkin space
  norm by
 $$
 \|f\|_{{F}^s_{p,q}} =\|f\|_{L^p} +\|f\|_{\dot{F}^s_{p,q}}.
 $$
 The inhomogeneous Triebel-Lizorkin space is  a Banach space
 equipped with the norm, $\|\cdot\|_{{F}^s_{p,q}}$.
We observe that $B^s_{p,p}(\Bbb R^n)=F^s_{p,p}(\Bbb R^n) $. The
Triebel-Lizorkin space is a generalization of many classical
 function spaces. Indeed, the followings are well established(see e.g.
 \cite{tri})
 $$F^0_{p, 2}(\Bbb R^n ) =\dot{F}^0_{p, 2}(\Bbb R^n)=L^p (\Bbb R^n ), \quad (1<p<\infty).$$
  $$\dot{F}^0_{1,2}(\Bbb R^n)=\mathcal{H}^1 (\Bbb R^n )\quad \mbox{and}
  \quad \dot{F}^0_{\infty,2}=BMO(\Bbb R^n). $$
 $$F^s_{p,2}(\Bbb R^n )=L^{s,p} (\Bbb R^n ).$$
We also note  sequence of continuous embeddings for the spaces close
to $L^\infty (\Bbb R^n)$(\cite{tri, jaw}).
 \bb\label{emb}
\dot{B}^{\frac{n}{p}}_{p,1} (\Bbb R^n)\hookrightarrow
\dot{B}^{0}_{\infty,1}(\Bbb R^n)
 \hookrightarrow L^\infty (\Bbb R^n) \hookrightarrow  BMO (\Bbb R^n)\hookrightarrow
 \dot{B}^{0}_{\infty,\infty}(\Bbb R^n).
\ee
  Given $0<s <1$, $1\leq p \leq \infty, 1\leq q \leq \infty$, we
introduce another function spaces $\mathcal{\dot{F}}^s_{p,q} $
defined by the
 seminorm,
 $$ \|f\|_{\mathcal{\dot{F}}^s_{p,q} } =\left\{
 \aligned \left\|
 \left(\intr\frac{|f(x)-f(x-y)|^q}{|y|^{n+s q}} dy\right)^{
 \frac{1}{q}}\right\|_{L^p (\Bbb R^n , dx)}
 & \mbox{if $1\leq p \leq \infty , 1\leq q <\infty$}\\
 \left\| ess\sup_{|y|>0 } \frac{|f(x) -f(x-y)|}{|y|^s } \right\|_{
 L^p (\Bbb R^n , dx)}
 & \mbox{if $1\leq p \leq \infty, q=\infty$}
 \endaligned \right. . $$
 On the other hand, the space $\mathcal{\dot{B}}^s_{p,q} $ is
defined by the
 seminorm,
 $$ \|f\|_{\mathcal{\dot{B}}^s_{p,q} } =\left\{
 \aligned
\left( \intr \frac{\|f(\cdot)-f(\cdot-y)\|^q_{L^p}}{|y|^{n+s q}}
dy\right)^{
 \frac{1}{q}}
 & \quad\mbox{if $1\leq p \leq \infty , 1\leq q <\infty$}\\
 ess\sup_{|y|>0 } \frac{\|f(\cdot) -f(\cdot-y)\|_{L^p }}{|y|^s }
 & \quad\mbox{if $1\leq p \leq \infty, q=\infty$}
 \endaligned \right. . $$
 Observe that, in particular,
 $\mathcal{\dot{F}}^{s}_{\infty,\infty}=
 \mathcal{\dot{B}}^{s}_{\infty,\infty}= C^s$, which
 is the usual H\"{o}lder
 seminormed space for $s\in \Bbb R_+\Bbb Z$. We also note that if $q=\infty$,
 $\mathcal{\dot{B}}^{s}_{p,\infty}=
 \mathcal{\dot{N}}^s_p$, which is the Nikolskii space.

The inhomogeneous version of those spaces, $\mathcal{{F}}^s_{p,q}
 $ and  $\mathcal{{B}}^s_{p,q} $
 are defined by their
 norms,
$$ \|f\|_{\mathcal{{F}}^s_{p,q} }= \|f\|_{L^p }
+\|f\|_{\mathcal{\dot{F}}^s_{p,q}}, \quad
\|f\|_{\mathcal{{B}}^s_{p,q} }= \|f\|_{L^p }
+\|f\|_{\mathcal{\dot{B}}^s_{p,q}}, $$ respectively.
 We note that for $0<s <1$, $ 2\leq p< \infty $, $q=2$,
 $\mathcal{{F}}^s_{p,2}
 \cong
 L^{p}_s (\Bbb R^n )$,
introduced above(see pp. 163,\cite{ste1}).
 If $\frac{n}{\min\{ p,q\} } <s <
 1$, $n<p<\infty$ and $n<q\leq \infty$, then
 $\mathcal{{F}}^s_{p,q} $ coincides with the Triebel-Lizorkin
 space $F^s_{p,q} (\Bbb R^n ) $ defined above(see pp. 101,
 \cite{tri}). On the other hand, for wider range of parameters, $0<s<1$, $0<p\leq \infty$,
 $0<q\leq \infty$, $\mathcal{{B}}^s_{p,q} $ coincides with the Besov
 space $B^s_{p,q} (\Bbb R^n )$ defined above.

\section{Local well-posedness and  blow-up criteria}
 \setcounter{equation}{0}

\subsection{Kato's local existence and the BKM criterion}

We  review briefly the key elements in the classical local existence
proof of solutions in the Sobolev space $H^m (\Bbb R^n)$, $m>n/2+1$,
essentially obtained by Kato in \cite{kat2}(see also \cite{maj3}).
After that we derive the celebrated  Beale, Kato and Majda's
criterion on finite time blow-up of the local solution in $H^m (\Bbb
R^n)$, $m>n/2+1$ in \cite{bea}.  Taking derivatives $D^\alpha $ on
the first equation of (E)
 and then taking $L^2$ inner
 product it with $D^\alpha v$, and summing over the multi-indices
 $\alpha $ with $|\alpha
 |\leq m$, we obtain
 \bqn
\lefteqn{\frac12 \frac{d}{dt} \|v \|_{H^m} ^2
 =-\sum_{|\alpha|\leq m} (D^\alpha (v\cdot \nabla )v
 -(v\cdot \nabla )D^\alpha v , D^\alpha v)_{L^2}}\hspace{.1in}\\
 &&-\sum_{|\alpha|\leq m}((v\cdot \nabla )
 D^\alpha v, D^\alpha v )_{L^2}-\sum_{|\alpha|\leq m} (D^\alpha \nabla\, p
 , D^\alpha v )_{L^2}\\
 &&\quad=I+II+III.
 \eqn
 Integrating by part, we obtain
  $$
  III=\sum_{|\alpha|\leq m}(D^\alpha  p
 , D^\alpha \mathrm{div} \, v)_{L^2} =0.
$$
Integrating by part again,  and using the
  fact div $v=0$, we have
$$
  II=-\frac12 \sum_{|\alpha|\leq m} \int_{\Bbb R^n} (v\cdot \nabla
  )|D^\alpha v |^2 dx
  =\frac12 \sum_{|\alpha|\leq m}\int_{\Bbb R^n} \mbox{div}\, v |D^\alpha
  v |^2 dx =0
  .
$$
 We  now use the so called {\it commutator type of
 estimate}(\cite{kla}),
$$
\sum_{|\alpha|\leq m}\|D^{\alpha} (fg)-fD^\alpha g\|_{L^2}
  \leq C (\|\nabla f\|_{L^\infty}
 \|g\|_{H^{m-1}} + \|f\|_{H^m} \|g\|_{L^\infty} ),
 $$
  and obtain
$$
 I\leq
 \sum_{|\alpha|\leq m} \|D^\alpha (v\cdot \nabla )v
 -(v\cdot \nabla )D^\alpha  v\|_{L^2} \|v\|_{H^{m}}
 \leq C\|\nabla v \|_{L^\infty} \|v \|_{H^{m}}^2 .
$$
 Summarizing the above estimates, I,II,III,  we have
 \bb\label{sum}
  \frac{d}{dt} \|v \|_{H^{m}} ^2
 \leq C \|\nabla v\|_{L^\infty} \|v \|_{H^{m}} ^2.
 \ee
Further estimate, using the {\it Sobolev inequality}, $\|\nabla
v\|_{L^\infty} \leq C \|v\|_{H^m} $ for $m>n/2+1$, gives
 $$
  \frac{d}{dt} \|v \|_{H^{m}} ^2
 \leq C  \|v \|_{H^{m}} ^3.
$$
  Thanks to Gr{o}nwall's lemma we have the local in time
  uniform estimate
\bb \label{ap1}
 \|v(t)\|_{H^m} \leq \frac{\|v_0\|_{H^m}}{1-Ct \|v_0\|_{H^m}
 }
 \leq 2 \|v_0\|_{H^m}
\ee
 for all $t\in [0, T ]$, where $T=\frac{1}{2C\|v_0\|_{H^m}}$.
 Using this estimate we can also deduce the estimate
\bb\label{ap2} \sup_{0\leq t\leq T} \left\| \frac{\partial
v}{\partial t
 }\right\|_{H^{m-1}} \leq C (\|v_0 \|_{H^m} )
 \ee
 directly from (E).
The estimates (\ref{ap1}) and (\ref{ap2}) are  the two key a priori
estimates for the construction of the local
 solutions. For actual elaboration of the proof we approximate the Euler
 system by mollification, Galerkin projection, or
 iteration of successive linear systems, and construct a sequence of
 smooth approximate solutions to (E), say $\{ v_k (\cdot ,t)\}_{k\in \Bbb N}$ corresponding to
 the initial data $ \{v_{0,k}\}_{k\in \Bbb N}$ respectively with $v_k \to v_0 $
  in $H^m (\Bbb R^n)$.
The estimates for the approximate solution sequence provides us with
the uniform estimates of $\{ v_k\}$ in $L^\infty ([0, T]; H^m (\Bbb
R^n )) \cap Lip ([0, T]; H^{m-1}(\Bbb R^n) )$. Then, applying the
standard Aubin-Nitche compactness lemma, we can pass to the limit
$k\to \infty$ in the equations for the approximate solutions, and
can show that the limit $v=v_\infty$ is a solution of the (E) in
$L^\infty ([0, T]); H^m (\Bbb R^n ))$. By further argument we can
actually show that the limit $v$ belongs to $ C([0, T]; H^m (\Bbb
R^n ))\cap AC([0, T]; H^{m-1} (\Bbb R^n ))$, where $AC([0,T]; X )$
denotes the space of $X$ valued absolutely continuous functions on
$[0,T]$. The general scheme of such existence proof is standard, and
is described in detail in \cite{maj1} in the general type of
hyperbolic conservation laws. The approximation  of the Euler system
by mollification was done for  the construction of local solution of
the Euler(and the
Navier-Stokes) system in \cite{maj3}. \\

Regarding the question of finite time blow-up of the local classical
solution in $H^m (\Bbb R^n )$, $m>n/2 +1$, constructed above, the
celebrated Beale-Kato-Majda theorem(called the BKM criterion) states
that \bb\label{bkm}
 \lim\sup_{t\nearrow T_*} \|v(t)\|_{H^s }=\infty \quad\mbox{if and only
 if}\quad
 \int_0 ^{T_*} \|\omega (s)\|_{L^\infty} ds =\infty.
 \ee
 We outline the proof of this theorem below(for more details see \cite{bea, maj3}).
We first recall the Beale-Kato-Majda's version of the {\it
logarithmic
  Sobolev inequality},
\bb\label{lo}
  \|\nabla v \|_{L^\infty}
   \leq C\|\o \|_{L^\infty} (1+\log (1+\|v
  \|_{H^{m}} ))+C \|\o\|_{L^2}
\ee
  for $m>n/2 +1$.
 Now suppose $\int_0 ^{T_*} \|\o (t)\|_{L^\infty } dt := M(T_*)< \infty .$
 Taking $L^2$ inner product the first equation of (E) with $\o$, then after
 integration by part we obtain
 $$
 \frac12\frac{d}{dt} \|\o \|_{L^2} ^2 =((\o\cdot \nabla )v ,\o
 )_{L^2}\leq \|\o\|_{L^\infty} \|\nabla v \|_{L^2} \|\o \|_{L^2}
 = \|\o\|_{L^\infty}\|\o \|_{L^2}^2,
 $$
 where we used the identity $\|\nabla v\|_{L^2}= \|\o\|_{L^2}$.
 Applying the Gronwall lemma, we obtain
 \bb\label{ineq1} \|\o (t)\|_{L^2} \leq \|\o_0 \|_{L^2}
 \exp\left( \int_0 ^{T_*}
 \|\o (s)\|_{L^\infty} ds \right)= \|\o_0 \|_{L^2}\exp[{M(T_*)}].
 \ee
 for all $t\in [0, T_*]$.
 Substituting (\ref{ineq1}) into (\ref{lo}), and combining this with
 (\ref{sum}), we have
 \[
  \frac{d}{dt} \|v \|_{H^{m}} ^2
 \leq C \left[ 1+\|\o  \|_{L^\infty} [1+\log (1+\|v
  \|_{H^{m}} )\right] \|v \|_{H^{m}} ^2
  \]
 Applying the Gr{o}nwall lemma we deduce
  \bb\label{fini}
  \|v (t)\|_{H^{m}} \leq \|v_0 \|_{H^{m}}
  \exp\left[ C_1 \exp\left(C_2 \int_0 ^{T_*}
  \|\o (\tau)\|_{L^\infty} d\tau\right)\right]
  \ee
  for all $t\in [0,T_*]$ and for some constants $C_1$ and  $C_2$ depending on $M(T_*)$.
 The inequality (\ref{fini}) provides the  with the necessity part of
 (\ref{bkm}).
The sufficiency part is an easy consequence of the Sobolev
inequality,
 $$ \int_0 ^{T_*} \|\omega (s)\|_{L^\infty} ds \leq
T_* \sup_{0\leq t\leq T_*} \|\nabla v(t)\|_{L^\infty} \\
\leq  C T_* \sup_{0\leq t\leq T_*}\|v(t)\|_{H^m}
 $$
 for $m>n/2+1$.
There are many other results of local well-posedness in various
function spaces(see \cite{cha1, cha2, cha3, cha6, che1,
che2,kat1,kat3,kat4, lic, tem1, tem2, vis1, vis2, yud1}). For the
local existence proved in terms of a geometric formulation see
\cite{ebi}. For the BKM criterion for solutions in the H\"{o}lder
space see \cite{bah}. Immediately after the BKM result appeared,
Ponce derive similar criterion in terms of the deformation
tensor(\cite{pon}). Recently, Constantin proved local well-posedness
and a blow-up criterion in terms of the active vector
formulation(\cite{con5}).

\subsection{Refinements of the BKM criterion}

The first refinement of the BKM criterion was done by Kozono and
Taniuchi in \cite{koz1}, where they proved
\begin{theorem}
Let $s>n/p+1$. A solution $v$ of the Euler equations belonging to
$C([0,T_*);W^{s,p}(\Bbb R^n)) \cap C^1 ([0, T_* ); W^{s-2, p}(\Bbb
R^n ))$ blows up at $T_*$ in $W^{s,p}(\Bbb R^n))$, namely
$$
\lim\sup_{t\nearrow T_*}\|v(t)\|_{W^{s,p}} =\infty \quad \mbox{if
and only if}\quad \int_0 ^{T_*} \|\o \|_{BMO} =\infty.
$$
\end{theorem}
The proof is based on the following version of the logarithmic
Sobolev inequality for $f\in W^{s,p} (\Bbb R^n)$, $s>n/p$,
$1<p<\infty$,
$$
\|f\|_{L^\infty} \leq C \left(1+\|f\|_{BMO} (1+\log^+
\|f\|_{W^{s,p}} ) \right).
$$
(see \cite{koz1} for details of the proof).  We recall  now the
embedding relations (\ref{emb}). Further refinement of the above
theorem is the following(see\cite{cha1, cha6}).
\begin{theorem}
 \begin{description}
 \item[(i)] {\rm{ (super-critical case)}} Let  $s>n/p+1, p\in (1, \infty), q\in [1,
 \infty]$.
  Then, the local in time solution $v\in C([0, T_*); B^s_{p, q} (\Bbb R^n))$
blows up at $T_* $ in $B^s_{p, q} (\Bbb R^n)$, namely
 $$
 \lim\sup_{t\nearrow T_*} \|v(t)\|_{B^s_{p,q}} =\infty
 \quad\mbox{
 if and only if}
 \quad
\int_0 ^{T_*}  \|\omega (t)\|_{\dot{B}^0_{\infty ,\infty}}dt
 =\infty .
$$
\item[(ii)] {\rm{(critical case)}} Let  $p\in (1, \infty)$.
Then, the local in time solution $v\in C([0, T_*); B^{n/p+1}_{p,1}
(\Bbb R^n))$ blows up at $T_* $ in $ B^{n/p+1}_{p,1} (\Bbb R^n))$,
namely
$$
 \lim\sup_{t\nearrow T_*} \|v(t)\|_{ B^{n/p+1}_{p,1}} =\infty
\quad\mbox{
 if and only if
 }\quad
\int_0 ^{T_*}  \|\omega (t)\|_{\dot{B}^0_{\infty ,1}}dt
 =\infty.
$$
\end{description}
\end{theorem}

The proof of (i) is based on the following version of the
logarithmic Sobolev inequality for $f\in B^s_{p, q} (\Bbb R^n)$ with
$s> n/p$ with $p\in (1, \infty )$, $q\in [1, \infty]$.
$$
 \|f\|_{L^\infty}
 \leq C(1+\|f\|_{\dot{B}^0_{\infty, \infty} }(\log^+ \|f\|_{B^s_{p,q}}  +1))
 $$
In \cite{koz2} Kozono, Ogawa and Taniuchi  obtained similar results
to
(i) above independently.\\
\ \\
In all of the above criteria, including the BKM theorem, we need to
control all of the three components of the vorticity vector to
obtain regularity. The following theorem proved in \cite{cha7}
states that actually we only need to control two components of the
vorticity in the slightly stronger norm than the $L^\infty$
norm(recall again the embedding (\ref{emb})).
\begin{theorem} Let $m>5/2$.
Suppose $v\in C([0,T_*);H^m (\Bbb R^3))$ is the local classical
solution  of (E) for some $T_1 >0$, corresponding to the initial
data $v_0 \in H^m (\Bbb R^3)$, and $\o=$ curl $v$ is its vorticity.
We decompose
 $\o=\tilde{\o}+\o^3 e_3$, where $\tilde{\o}=\o^1 e_1 +\o^2 e_2$,
 and $\{e_1, e_2, e_3\}$ is the canonical basis of $\Bbb R^3$. Then,
$$
  \lim\sup_{t\nearrow T_*} \|v(t)\|_{H^m} =\infty \quad \mbox{if
and only if} \quad \int_0 ^{T_*}\|\tilde{\o}
(t)\|_{\dot{B}^0_{\infty, 1}} ^2 dt =\infty.
 $$
\end{theorem}
Note that   $\tilde{\o}$ could be the projected component of $\o$
onto any plane in $\Bbb R^3$. For the solution $v=(v^1, v^2, 0)$ of
the Euler equations on the $x_1-x_2$ plane, the vorticity is $\o
=\o^3 e_3$ with $\o_3=\partial_{x_1} v^2 -\partial_{x_2} v^1$, and
$\tilde{\o}\equiv 0$. Hence, as a trivial application of the above
theorem we reproduce the well-known global in time  regularity
for the 2D Euler equations.\\
\ \\
Next we present recent results on the blow up criterion in terms of
 hessian of the pressure.  As in the introduction  we use  $P=(P_{ij})$,
$S=(S_{ij})$ and $\xi$ to denote  the hessian of the pressure, the
deformation tensor and the vorticity direction field respectively,
introduced in section 1. We also introduce the notations
$$ \frac{S\xi}{|S\xi |} =\zeta, \quad \zeta \cdot P \xi = \mu.
$$
The following is proved in \cite{cha15}.
\begin{theorem}
 If the solution $v(x,t)$ of the 3D Euler
system with $v_0 \in H^m (\Bbb R^3 )$, $m>5/2$, blows up at $T_*$,
namely $\lim\sup_{t\nearrow T_*} \|v(t)\|_{H^m} =\infty$, then
necessarily,
 $$
 \int_0 ^{T_*} \exp\left(\int_0 ^\tau \|\mu (s)\|_{L^\infty}ds\right)
 d\tau =\infty.
 $$
\end{theorem}

Similar criterion in terms of the hessian of pressure, but with
different detailed geometric configuration  from the above theorem
is obtained by Gibbon, Holm, Kerr and Roulstone in \cite{gib3}.
Below we denote $\xi_p=\xi \times P \xi$.
\begin{theorem}
Let $m \geq 3$ and $\Bbb T^3 =\Bbb R^3/\Bbb Z^3$ be a periodic box.
Then, there exists a global solution of the Euler equations $v\in
C([0, \infty ); H^m (\Bbb T^3 ))\cap C^1 ([0, \infty); H^{m-1} (\Bbb
T^3 ))$ if
$$\int_0 ^T \|\xi_p (t)\|_{L^\infty} dt <\infty, \quad \forall t\in (0,
T)
$$
excepting the case where $\xi$ becomes collinear with the
eigenvalues of $P$ at $T$.
\end{theorem}

Next, we consider the axisymmetric solution of the Euler equations,
which means velocity field $v(r, x_3, t)$, solving the Euler
equations, and having the representation
$$
v(r, x_3, t)=v^r (r,x_3,t)e_r +v^\theta (r,x_3,t)e_\theta +v^3
(r,x_3,t)e_3
$$
 in the cylindrical coordinate system, where
 $$e_r = (\frac{x_1}{r}, \frac{x_2}{r}, 0), \quad
 e_\theta = (-\frac{x_2}{r}, \frac{x_1}{r}, 0),\quad e_3=(0,0,1),\quad
 r=\sqrt{x_1^2 +x_2^2}.
 $$
In this case also the question of finite time blow-up  of solution
is wide open(see e.g. \cite{gra1, gra2,caf1} for studies in such
case). The vorticity $\o =$ curl $\,v$ is computed as
\[ \w = \w^r e_r+ \w^{\t}e_{\t} + \w^3 e_3 ,  \]
where
\[ \w^r = -\partial_{x_3} v^{\t},
   \quad \w^{\t} = \partial_{x_3} v^r - \partial_r v^3,
 \quad \w ^3 = \frac{1}{r}\partial_r(rv^{\t}).   \]
We denote
\[ \tv = v^r e_r + v^3 e_3 ,\qquad \tilde{\w} = \w^r e_r + \w^3 e_3 .\]
Hence, $  \w =  \tilde{\w}+\vec{\o}_\t $, where $\vec{\o}_\t
=\omega^\theta e_\theta$. The Euler equations for the axisymmetric
solution are
$$
 \left\{ \aligned
 &\frac{\partial v^r}{\partial t} +(\tv\cdot \tilde{\nabla} )v^r
 =-\frac{\partial p}{\partial r} ,\\
 &
 \frac{\partial v^\theta}{\partial t} +(\tv\cdot \tilde{\nabla}
 )v^\t
 =-\frac{v^r v^\t}{ r} ,\\
 &
 \frac{\partial v^3}{\partial t} +(\tv\cdot \tilde{\nabla} )v^3
 =-\frac{\partial p}{\partial x_3} ,\\
&
 \textrm{div }\, \tv =0 ,\\
&
 v(r,x_3,0)=v_0 (r,x_3),
\endaligned \right.
$$
where $\tilde{\nabla} =e_r \frac{\partial}{\partial r}
  +e_3  \frac{\partial}{\partial x_3}.$
In the axisymmetric Euler equations the vorticity formulation
becomes
$$
 \left\{\aligned &\frac{\partial\w^r}{\partial t} +(\tv \cdot
\tilde{\nabla} )
=\w^r(\tw \cdot \tilde{\nabla}) v^r     \\
& \frac{\partial\w^3}{\partial t} +(\tv \cdot
\tilde{\nabla} )=\w^3(\tw \cdot \tilde{\nabla}) v^3 \\
&
 \left[\frac{\partial}{\partial t} + \tv \cdot
\tilde{\nabla} \right] \left(\frac{\w ^{\t}}{r}\right)
     = (\tw \cdot \tilde{\nabla}) \left(\frac{v^{\t}}{r}\right)
     \\
&\textrm{div }\, \tv = 0 , \quad  \textrm{curl }\, \tv
      =\vec{\o}^\t.
\endaligned \right.
$$
In the case of axisymmetry we only need to control just one
component of vorticity(the angular component) to get the regularity
of solution. The following theorem is proved in \cite{cha21}.

 \begin{theorem} Let $v\in C([0,T_*);H^m (\Bbb R^3))$, $m>5/2$, be the
local classical axisymmetric solution of (E), corresponding to an
axisymmetric initial data $v_0 \in H^m (\Bbb R^3 )$. Then, the
solution blows up in $H^m (\Bbb R^3 )$ at $T_*$ if and only if for
all $(\gamma, p)\in (0,1)\times [1, \infty]$ we have
 \bq\label{axiscri}
 \lefteqn{\int_0 ^{T_*}
\|\o_\t  (t)\|_{L^\infty} dt +\int_0 ^{T_*} \exp\Big[ \int_0 ^t
\left\{ \|\o_\t (s)\|_{L^\infty} (1+\log^+ (\|\o _\t
(s)\|_{C^{\gamma}}
\|\o_\t (s)\|_{L^p} )) \right. }\hspace{2.in}\n \\
&&\left. +\|\o_\t (s)\log^+ r \|_{L^\infty}\right\}ds\Big] dt
=\infty.
 \eq
\end{theorem}
We observe that although we need to control only $\w_\t$ to get the
regularity, the its norm, which is in $C^\gamma$, is higher than the
$L^\infty$ norm used in the BKM criterion. If we use  the `critical'
Besov space $\dot{B}^0_{\infty,1}(\Bbb R^3)$ we can derive slightly
sharper criterion than Theorem 2.6 as follows(see \cite{cha7} for
the proof).
\begin{theorem} Let $v\in C([0,T_*);H^m (\Bbb R^3))$ be the
local classical axisymmetric solution of (E), corresponding to an
axisymmetric initial data $v_0 \in H^m (\Bbb R^3 )$. Then,
\bb\label{axiscri1}
  \lim\sup_{t\nearrow T_*} \|v(t)\|_{H^m} =\infty  \quad
\mbox{if and only if} \quad \int_0 ^{T_*}\|\vec{\o}_\theta
(t)\|_{\dot{B}^0_{\infty, 1}}  dt =\infty.
 \ee
\end{theorem}
We observe that contrary to (\ref{axiscri}) we do not need to
control the high regularity norm, the $C^{\gamma}$ norm of vorticity
in (\ref{axiscri1}). We can also have the regularity of the
axisymmetric Euler equation by controlling only one component of the
velocity, the swirl velocity $v^\theta$ as in the follows proved in
\cite{cha22}.

\begin{theorem} Let $v\in C([0,T_*);H^m (\Bbb R^3))$, $m>5/2$, be the
local classical axisymmetric solution of (E), corresponding to an
axisymmetric initial data $v_0 \in H^m (\Bbb R^3 )$. Then, the
solution blows up in $H^m (\Bbb R^3 )$ at $T_*$ if and only if
 $$ \int_0 ^{T_*}
 \left( \|\tilde{\nabla}v^\t \|_{L^\infty}
 +\left\|\frac{\partial v^\t}{\partial r} \right\|_{L^\infty}
 \left\| \frac{1}{r}\frac{\partial v^\t}{\partial
 x_3}\right\|_{L^\infty}\right) dt =\infty.
 $$
 \end{theorem}

\subsection{Constantin-Fefferman-Majda's and other related results}

In order to study the regularity problem of the 3D Navier-Stokes
equations Constantin and Fefferman investigated the geometric
structure of the integral kernel in the vortex stretching term more
carefully, and discovered the phenomena of `depletion effect' hidden
in the integration(\cite{con8}, see also \cite{con2} for detailed
exposition related to this fact). Later similar geometric structure
 of the vortex stretching term was studied extensively also in the blow-up
 problem of the 3D Euler equations by Constantin, Fefferman and
 Majda(\cite{con9}).
 Here we first present their results in detail, and  results in
\cite{cha10}, where the BKM criterion and the
Constantin-Fefferman-Majda's criterion are interpolated in some
sense. Besides those results presented in this subsection we also
mention that there are other interesting geometric approaches to the
Euler equations such as the quaternion formulation by
Gibbon(\cite{gib1, gib2, gib3}).  We begin with a definition in
\cite{con9}. Given a set $W\in \Bbb R^3$ and $r>0$ we use the
notation $B_r(W)=\{ y\in B_r (x)\, ;\, x\in W\}$.
\begin{definition}
A set $W_0 \subset \Bbb R^3$ is called smoothly directed if there
exists $\rho>0$ and $r, 0<r\leq \rho/2$ such that the following
three conditions are satisfied.
\end{definition}
\begin{itemize}
\item[(i)] For every $a\in W^*_0
=\{ q\in W_0 \, ; \, |\o_0 (q)|\neq 0\},$ and all $t\in [0,T)$, the
vorticity direction field $\xi (\cdot ,t)$ has a Lipshitz
extension(denoted by the same letter) to the Euclidean ball of
radius $4\rho$ centered at $X(a,t)$ and
$$
M=\lim_{t\to T} \sup_{a\in W^*_0} \int_0 ^t \|\nabla \xi (\cdot
,t)\|_{L^\infty (B_{4\rho} (X(a,t)))} dt <\infty.
$$
\item[(ii)] The inequality
$$ \sup_{B_{3r} (W_t )}|\o (x,t)|\leq m
\sup_{B_r (W_t )} |\o (x,t)|$$ holds for all $t\in [0,T)$ with
$m\geq 0$ constant.
\item[(iii)]  The inequality
$$ \sup_{B_{4\rho} (W_t )}|v (x,t)| \leq U$$
holds for all $t\in [0, T)$.
\end{itemize}
The assumption (i) means that the direction of vorticity is well
behaved in a neighborhood of a bunch of trajectories. The assumption
(ii) states that this neighborhood is large enough to capture the
local intensification of $\o$.  Under these assumptions the
following theorem is proved in \cite{con9}.

\begin{theorem}
Assume $W_0$ is smoothly directed. Then there exists $\tau >0$ and
$\Gamma$ such that
$$ \sup_{B_r (W_t )} |\o (x,t)|\leq \Gamma
\sup_{B_\rho (W_{t_0})} |\o (x,t_0 )| $$ holds for any $0\leq t_0
<T$ and $0\leq t-t_0 \leq \tau$.
\end{theorem}

They also introduced the notion of regularly directed set, closely
related to the geometric structure of the kernel defining vortex
stretching term.
\begin{definition}
We sat that a set $W_0$ is regularly directed if
 there exists $\rho >0$ such that
 $$ \sup_{a W^*_0} \int_0 ^T K_{\rho} (X(a,t))dt <\infty$$
 where
 $$ K_{\rho} (x)=\int_{|y|\leq \rho} |D(\hat{y},
 \xi (x+y), \xi (x)) | |\o (x+y)|\frac{dy}{|y|^3}
 $$
 and
 $$
 D(\hat{y}, \xi (x+y) ,\xi (x))=(\hat{y}\cdot \xi (x)) \mathrm{Det}
 (\hat{y}, \xi (x+y), \xi (x)).
 $$
\end{definition}
Under the above assumption on the regularly directed sets the
following is proved also in \cite{con9}.

\begin{theorem}
Assume $W_0$ is regularly directed. Then there exists a constant
$\Gamma $ such that
$$ \sup_{a \in W_0} |\o (X(a,t),t)| \leq \Gamma \sup_{a \in W_0}
|\o_0 (a)|
$$
holds for all $t\in [0,T]$.
\end{theorem}

The original studies by Constantin and Fefferman in \cite{con8}
about the Navier-Stokes equations, which motivated the above
theorems, are concerned mainly about the  regularity of solutions in
terms of  the vorticity direction fields $\xi$. We recall, on the
other hand, that the BKM type of criterion controls the magnitude of
vorticity to obtain regularity. Incorporation of both the direction
and the magnitude of vorticity to obtain regularity for the 3D
Navier-Stokes equations was first initiated  by Beir\~{a}o da Veiga
and Berselli in \cite{bei2}, and developed further by Beir\~{a}o da
Veiga in \cite{bei1}, and finally refined in an `optimal' form in
\cite{cha18}(see also \cite{cha20} for a localized version). We now
present the Euler equation version of the result in \cite{cha18}.

Below we use the notion of particle trajectory $X(a ,t)$, which is
defined  by the classical solution $v(x,t)$  of (E).  Let us denote
 $$
 \Omega_0 =\{ x\in \Bbb R^3 \, |\, \o_0 (x) \neq 0\},\quad
 \Omega_t = X(\O_0 ,t).
 $$
 We note that the direction field of the vorticity,
 $\xi (x,t)=\o (x,t)/ |\o (x,t)| $, is well-defined if
 $x\in \O_t$ for $v_0 \in C^1 (\Bbb R^3 )$ with $\O_0 \neq
 \emptyset$.
 The following is the
 main theorem proved in \cite{cha10}.
\begin{theorem}
Let $v(x,t)$ be the local classical solution to  $(E)$ with initial
data $v_0 \in H^m (\Bbb R^3)$, $m>5/2$, and $\o(x,t)=$curl $v(x,t)$.
We assume $\O_0 \neq\emptyset$.
 Then, the solution can be continued up to $T +\delta$ as the
 classical solution, namely
$v(t)\in C([0,T+\delta];{H^m} (\Bbb R^3 ))$ for some $\delta >0$,
 if there exists
  $p,p', q,q', s, r_1,r_2,r_3$ satisfying the following
conditions,
 \bb\label{th1}
  \frac{1}{p} +\frac{1}{p'} =1, \qquad  \frac{1}{q} +\frac{1}{q'} =1,
  \ee
  and
  \bb\label{th2}
\frac{1}{r_1}+\frac{p'}{r_2}\left(1-\frac{sq'}{3} \right)+
 \frac{1}{r_3} \left\{ 1-p' \left( 1-\frac{sq'}{3} \right)\right\}
 =1
 \ee
 with
 \bb\label{th3}
 0<s<1, \qquad 1\leq \frac{3}{sq'} <p \leq \infty, \qquad  \quad 1\leq
 q\leq\infty,
 \ee
 and
 \bb\label{th4}
 r_1\in [1, \infty],\,  r_2 \in \left[p'\left(1-\frac{sq'}{3}
 \right),\infty\right],\,  r_3 \in \left[1-p'\left(1-\frac{sq'}{3}
 \right), \infty\right]
\ee
 such that for direction field $\xi (x,t)$, and the magnitude of
 vorticity $|\o (x,t)|$ the followings hold;
\bb\label{reg}
 \int_0 ^{T} \|\xi (t )\|^{r_1}_{ \mathcal{\dot{F}}^s_{\infty,q} (\O_t) } dt <\infty,
 \ee
 and
 \bb\label{int}
  \int_0 ^{T}\|\o(t) \|^{r_2}_{L^{pq'}(\O_t)} dt +
   \int_0 ^{T}\|\o(t) \|^{r_3}_{L^{q'}(\O_t)} dt <\infty.
 \ee
\end{theorem}
\ \\
In order to get insight implied by the above theorem let us consider
the special case of $p=\infty, q=1$. In this case the conditions
(\ref{reg})-(\ref{int}) are satisfied if
 \bb\label{remm1}
 \xi (x,t )\in L^{r_1}(0,T; C^{ s} (\Bbb R^3 )),
 \ee
 \bb\label{remm2}
  \o
(x,t) \in L^{r_2} (0, T; L^{\infty} (\Bbb R^3 ))\cap L^{r_3} (0, T;
L^{\infty} (\Bbb R^3 )).
 \ee
with
 \bb\label{remm3}
 \frac{1}{r_1}+\frac{1}{r_2}\left(1-\frac{s}{3} \right)+
 \frac{s}{3r_3}
 =1.
 \ee
 Let us {\em formally} pass $s\to 0$
  in (\ref{remm1}) and (\ref{remm3}),
 and choose $r_1=\infty$ and $ r_2
 =r_3=1$, then we find that the conditions (\ref{remm1})-(\ref{remm2}) reduce to the
BKM condition,
 since the condition $\xi (x,t )\in L^{\infty}(0,T; C^{0} (\Bbb R^3
 ))\cong L^{\infty}((0,T)\times \Bbb R^3)$
 is obviously satisfied due to the fact that  $|\xi (x,t)|\equiv 1$.\\
 The other case of interest is $q'=3/s$, where
 (\ref{reg})-(\ref{int}) are satisfied if
 \bb\label{remm4}
 \xi (x,t )\in L^{r_1}(0,T;
 \mathcal{\dot{F}}^s_{\infty,\frac{3}{3-s}}(\Bbb R^3)
 ),\qquad
  |\o
(x,t)| \in L^{r_2} (0, T; L^{\frac{3}{s}} (\Bbb R^3 )).
 \ee
 with $1/r_1+1/r_2=1.$ The condition (\ref{remm4}) shows explicitly the
 mutual compensation between the regularity  of the
 direction field and the integrability  of the vorticity
 magnitude in order to control regularity/singularity of  solutions of the
 Euler equations.\\
 \ \\
 Next we review the result of non-blow-up conditions due to Deng,
 Hou and Yu\cite{den1,den2}.
 We consider a time $t$  and a vortex line segment $L_t$ such that
 the maximum of vorticity over the whole domain is comparable to the
 maximum of vorticity on over $L_t$, namely
 $$ \Omega (t):=\sup_{x\in \Bbb R^3}|\o (x,t)|\sim \max_{x\in L_t} |\o (x,t)|.$$
 We denote $ L(t):=\mbox{arc length of $L_t$}$; $\xi$, $\mathbf{n}$
 and $\kappa$ are
  the unit tangential and the unit normal vectors to $L_t$ and
  the curvature of $L_t$ respectively. We also use the  notations,
 \bqn
  U_\xi (t)&:=&\max_{x,y\in L_t} |(v\cdot \xi)(x,t)-(v\cdot \xi
 )(y,t)|,\\
  U_n (t)&:=&\max_{x\in L_t} |(v\cdot \mathbf{n})(x,t)|,\\
  M(t)&:=&\max_{x\in L_t} |(\nabla \cdot \xi)(x,t)|,\\
  K(t)&:=&\max_{x\in L_t} \kappa (x,t).
\eqn
 We denote by $X(A, s,t)$ the image by the trajectory map at time
 $t>s$ of fluid particles at $A$ at time $s$. Then, the following is
 proved in \cite{den2}.
 \begin{theorem}
 Assume that there is a family of vortex line segment $L_t$ and $T_0
 \in [0, T^*)$, such that $X(L_{t_1} , t_1 ,t_2 ) \supseteq L_{t_2}$
 for all $T_0 <t_1 <t_2 <T^*$. Also assume that $\Omega (t)$ is
 monotonically increasing and $\max_{x\in L_t} |\o (x,t)| \geq c_0
 \Omega (t)$ for some $c_0$ when $t$ is sufficiently close to $T^*$.
 Furthermore, we assume there are constants $C_U ,C_0 ,c_L$ such
 that
 \begin{enumerate}
 \item $[U\xi (t)+U_n (t) K(t) L(t)] \leq C_U (T^* -t)^{-A}$ for
 some constant $A\in (0,1)$,
 \item $M(t)L(t), K(t)L(t) \leq C_0$,
 \item $L(t) \geq c_L (T^* -t)^B$ for some constant $B\in (0, 1)$.
 \end{enumerate}
 Then there will be no blow-up in the 3D incompressible Euler flow
 up to time $T^*$, as long as $B< 1-A$.
 \end{theorem}
 In the endpoint case of $B=1-A$ they deduced the following
 theorem(\cite{den1}).
 \begin{theorem}
 Under the same assumption as in Theorem 2.10, there will be no
 blow-up in the Euler system up to time $T^*$ in the case $B=1-A$,
 as long as the following condition is satisfied:
 $$ R^3 K< y_1 \left( R^{A-1} (1-A)^{1-A} /(2-A)^{2-A} \right),
 $$
 where $R=e^{C_0} /c_0, K:=\frac{C_U c_0}{c_L (1-A)}$, and $y_1 (m)$
 denotes the smallest positive $y$ such that $m=y/(1+y)^{2-A}$.
\end{theorem}
 We refer \cite{den1,den2} for discussions on the various connections of Theorem 2.10 and Theorem 2.11
  with numerical computations.
\section{Blow-up scenarios}
 \setcounter{equation}{0}

\subsection{Vortex sheet collapse}

We recall that a vortex line is an integral curve of the vorticity,
and a vortex tube is a tubular neighborhood in $\Bbb R^3$ foliated
by vortex lines. Numerical simulations(see e.g. \cite{cho}) show
that vortex tubes grow and thinner(stretching),  and folds before
singularity happens. We review here the result by Cordoba and
Fefferman\cite{cor4} excluding a type of vortex tube collapse.

Let $Q=I_1 \times I_2 \times I_3\subset \Bbb R^3$ be a closed
rectangular box,  and  let $T>0$ be given. A regular tube is a
relatively open set $\Omega_t \subset Q$ parameterized by time $t\in
[0, T)$, having the form $\Omega_t =\{ (x_1, x_2 , x_3 )\in Q \, :\,
\theta (x_1,x_2,x_3, t)<0\}$ with $\theta \in C^1 (Q\times [0, T))$,
and satisfying the following properties:
$$ |\nabla _{x_1,x_2} \theta |\neq 0\, \mbox{for}\,
(x_1,x_2,x_3,t)\in Q\times [0, T), \theta (x_1,x_2,x_3,t) =0;
$$
$$\Omega_t (x_3):=\{ (x_1,x_2 )\in I_1 \times I_2 \, :\,
(x_1,x_2,x_3 )\in \Omega _t \}\, \mbox{is non-empty},
$$
for all $x_3 \in I_3, t\in [0, T)$;
$$ \mbox{closure} (\Omega_t (x_3 )) \subset \mbox{interior}
(I_1\times I_2 )
$$
for all $x_3 \in I_3, t\in [0, T)$.

Let $u(x,t)=(u_1 (x,t),u_2 (x,t),u_3 (x,t))$ be a $C^1$ velocity
field defined on $Q\times [0,T)$. We say that the regular tube
$\Omega_t $ moves with the velocity field $u$, if we have
$$
\left(\frac{\partial}{\partial t} +u\cdot \nabla_x \right)\theta =0
\, \mbox{whenever} \, (x,t)\in Q\times [0,T), \theta (x,t)=0.
$$
By the Helmholtz theorem we know that a vortex tube arising from a
3D Euler solution moves with the fluid velocity. The following
theorem proved by Cordoba and Fefferman(\cite{cor4}) says for the 3D
Euler equations that a vortex tube cannot reach zero thickness in
finite time, unless it bends and twists so violently that no part of
it forms a regular tube.

\begin{theorem}
Let $\Omega_t \subset Q (t\in [0,T))$ be a regular tube that moves
with $C^1$, divergence free velocity field $u(x,t)$.
$$ \mbox{If}\quad \int_0 ^T \sup_{x\in Q} |u(x,t)|dt <\infty,
\quad \mbox{then}\quad \lim\inf_{t\to T_-} \mbox{Vol}(\Omega_t )>0.
$$
\end{theorem}

\subsection{Squirt singularity}

 The theorem of excluding the regular vortex
tube collapse was generalized by Cordoba, Fefferman and de la
Lave(\cite{cor7}), which we review here. We first recall their
definition of squirt singularities. Let $\Omega \subset \Bbb R^n$ be
an open set. We denote $X_t (a)=X(a,t)$, which is a particle
trajectory generated by a $C^1$ vector field $u:\Omega \times
[0,T)\to \Bbb R^n$ such that div $u=0$. We also set $X_{t,s} (a)$ as
the position at time $t$ of the trajectory which at time $t=s$ is
$a$. We have obvious relations,
$$ X_t (a)=X_{t,0} (a), \quad X_{t,s} = X_t\circ X_s^{-1}, \quad
X_{t,s}\circ X_{s, s_1} =X_{t,s_1}.
$$
For $\mathcal{S} \subset \Omega$, we denote by
$$
X_{t,s}^\Omega \mathcal{S} =\{ x\in \Omega \, |\, x=X_t (a), \, a
\in \mathcal{S},\, X_s (a)\in \Omega ,\, 0\leq s\leq t\}.
$$
In other words, $X_{t,s}^\Omega \mathcal{S}$ is the evolution of the
set $\mathcal{S}$, starting at time $a$, after we eliminate the
trajectories which step out of $\Omega$ at some time. By the
incompressibility condition on $u$, we have
 that Vol($X_{t,s} \mathcal{S})$ is independent of $t$, and the
 function $t\mapsto$Vol($X_{t,s} \mathcal{S})$ is nonincreasing.

\begin{definition}
Let $\Omega_-,\Omega_+$ be open and bounded sets.
$\overline{\Omega_-} \subset \Omega_+$. Therefore, dist $(\Omega_- ,
\Bbb R^n  -\Omega_+ )\geq r >0$. We say that $u$ experiences a
squirt singularity in $\Omega_-$, at time $T>0$, when for every
$0\leq s <T$, we can find a set $\mathcal{S}_s \subset \Omega_+$
such that
\begin{itemize}
\item[(i)] $\mathcal{S}_s \cap \Omega_-$ has positive measure, $0\leq
s<T$,
 \item[(ii)] $\lim_{t\to T}\mathrm{Vol}(X_{t,s}^{\Omega^+}\mathcal{S}_s
 )=0$.
 \end{itemize}
 \end{definition}
 The physical intuition behind the above definition is that there
 is a region of positive volume so that all the fluid occupying it
 gets ejected from a slightly bigger  region in finite time.
Besides the vortex tube collapse singularity introduced in the
previous subsection the potato chip singularity and the saddle
collapse singularity, which will be defined below,  are also special
examples of the squirt singularity, connected with real fluid
mechanics phenomena.

\begin{definition}[potato chip singularity]
We say that $u$ experiences a potato chip singularity when we can
find continuous functions
$$f_\pm :\Bbb R^{n-1}\times [0, T) \to \Bbb R
$$
such that
 \bqn
 f_+ (x_1 ,\cdots ,x_{n-1}, t)\geq f_-(x_1 ,\cdots ,x_{n-1}, t),
 \, t\in [0, T], x_1, \cdots , x_{n-1} \in B_{2r} (\Pi x^0 ),\\
f_+ (x_1 ,\cdots ,x_{n-1},0 )\geq f_-(x_1 ,\cdots ,x_{n-1},0),
\quad x_1, \cdots , x_{n-1} \in B_{r} (\Pi x^0 ),\\
\lim_{t\to T_-} [f_+ (x_1 ,\cdots ,x_{n-1}, t)- f_-(x_1 ,\cdots
,x_{n-1}, t)]=0 \forall x_1, \cdots , x_{n-1} \in B_{2r} (\Pi x^0 )
\eqn
 and such that the surfaces
 $$ \Sigma _{\pm, t} =\{ x_n = f_\pm(x_1 ,\cdots ,x_{n-1}, t)
 \}\subset \Omega$$
 are transformed into each other by the flow
 $$X (\Sigma _{\pm, 0}, t) \supset  \Sigma _{\pm, t}. $$
 In the above $\Pi $ is projection on the first $n-1$ coordinates.
 \end{definition}
Previously to \cite{cor7} potato chip singularities were considered
in the 2D and 3D flows by C\'{o}rdoba and
 Fefferman(\cite{cor6},\cite{cor4a} respectively)
 in the name of `sharp front'. In particular the exclusion of sharp
 front in the 2D quasi-geostrophic equation is proved in
 \cite{cor6}. The following notion of saddle collapse singularity is relevant only
for 2D flows.
 \begin{definition}[saddle collapse singularity] We consider
 foliation of a neighborhood of the origin(with coordinates $x_1,
 x_2$) whose leaves are given by equations of the form
 $$
 \rho:= (y_1 \beta(t) +y_2 )\cdot (y_1 \delta (t) +y_2 )=cons
 $$
 and $(y_1 ,y_2 )=F_t (x_1 ,x_2 )$, where $\beta , \delta : [0, T)
 \to \Bbb R^+$ are $C^1$ foliations, $F$ is a $C^2$ function of
 $x,t,$ for a fixed $t$, and $F_t$ is an orientation preserving
 diffeomorphism.
 We say that the foliation experiences a saddle collapse when
 $$\lim\inf_{t\to T} \beta (t) +\delta (t) =0.
 $$
 If the leaves of the foliation are transported by a vector field
 $u$, we say that the vector field $u$ experiences a saddle
 collapse.
 \end{definition}
The exclusion of saddle point singularity in the 2D
quasi-geostrophic equation(see Section 4.3 below) was proved by
C\'{o}rdoba in \cite{cor2}. The following `unified' theorem is
proved in \cite{cor7}.
\begin{theorem}
If $u$ has a squirt singularity at $T$, then $\int_s ^T \sup_x
|u(x,t)|dt =\infty$ for all $ s\in (0, T)$. Moreover, if $u$ has a
potato chip singularity, then
$$
\int_s ^T \sup_x |\Pi u (x,t)|dt =\infty.
$$
 \end{theorem}

\subsection{Self-similar blow-up}

In this subsection we review the scenario of self-similar
singularity studied in \cite{cha13}. We first observe that the Euler
system (E) has scaling property that
  if $(v, p)$ is a
solution of the system (E), then for any $\lambda >0$ and $\alpha
\in \Bbb R $ the functions
 \bb
 \label{self}
  v^{\lambda, \alpha}(x,t)=\lambda ^\alpha v (\lambda x, \l^{\a +1}
  t),\quad p^{\l, \a}(x,t)=\l^{2\a}p(\l x, \l^{\a+1} t )
  \ee
  are also solutions of (E) with the initial data
  $ v^{\lambda, \alpha}_0(x)=\lambda ^\alpha v_0
   (\lambda x)$.
 In view of the scaling
  properties in (\ref{self}), the  self-similar blowing up
  solution $v(x,t)$ of (E), if it exists,  should be of the form,
  \bq
  \label{vel}
 v(x, t)&=&\frac{1}{(T_*-t)^{\frac{\a}{\a+1}}}
V\left(\frac{x}{(T_*-t)^{\frac{1}{\a+1}}}\right)
 \eq
 for  $\a \neq -1$ and $t$  sufficiently
 close to $T_*$.
 If we assume that initial vorticity $\omega_0$ has compact support,
then the nonexistence of self-similar blow-up of the form given by
(\ref{vel}) is rather immediate from the well-known formula, $
\omega (X(a,t),t)=\nabla_a X(a,t) \omega_0 (a)$. We want to
generalize this to a nontrivial case. Substituting (\ref{vel}) into
(E), we find that $V$ should be
  a solution of the system
\[
 (SE)
 \left\{
 \begin{aligned}
 \frac{\a}{\a +1} V+&\frac{1}{\a+1} (x\cdot\nabla ) V
 +(V\cdot \nabla )V =-\nabla  P,
 \\
  \textrm{div }\, V& =0
  \end{aligned}
  \right.
  \]
  for some scalar function $P$,
  which could be regarded as the Euler version of the Leray
  equations introduced in \cite{ler}. The question of existence of
  nontrivial solution to (SE) is equivalent to the  that of existence of
  nontrivial self-similar finite time blowing up
  solution to the Euler system of the form (\ref{vel}).
  Similar question for the
  3D Navier-Stokes equations was raised by  Leray in \cite{ler}, and
  answered  negatively by Necas, Ruzicka and Sverak(\cite{nec}), the result of which was
  refined later by Tsai in \cite{tsa}(see also \cite{mil} for a generalization).
 Combining the energy conservation with a simple scaling argument,
  the author of this article showed that if there exists a nontrivial self-similar finite time
blowing up solution, then its helicity should be zero(\cite{cha4}).
   Mainly due to lack of the laplacian term
  in the right hand side of the first equations of (SE), we cannot
 expect the maximum principle, which was crucial in the works in
  \cite{nec} and \cite{tsa} for the 3D Navier-Stokes equations.
  Using a completely different argument from those  previous ones,
   in \cite{cha13} it is proved that  there cannot be  self-similar blowing up solution to (E) of the form
  (\ref{vel}), if the vorticity decays sufficiently fast near infinity.
  Given a smooth velocity field $v(x,t)$,
  the particle trajectory mapping $a\mapsto X(a,t)$
The inverse $A(x,t):=X^{-1} (x ,t)$  is called the back to label
map, which satisfies $
 A (X(a,t),t)=a,$ and $ X(A
 (x,t),t)=x$.
The existence of the back-to-label
 map $ A(\cdot, t)$ for our smooth velocity $v(x,t)$  for
 $t\in (0, T_*)$, is guaranteed
 if we assume a uniform decay of  $v(x,t)$ near infinity, independent of the
 decay rate(see \cite{con5}). The following is proved in
 \cite{cha13}.
\begin{theorem}
There exists no finite time blowing up self-similar solution
$v(x,t)$ to the 3D Euler equations of the form
  (\ref{vel}) for $t \in (0 ,T_*)$ with $\a \neq -1$, if $v$ and $V$
  satisfy the following conditions:
  \begin{itemize}
  \item[(i)]For all $t\in (0, T_*)$ the particle trajectory mapping $X(\cdot ,t)$ generated
  by the classical solution $v\in C([0,T_*);C^1(\Bbb R^3;\Bbb R^3 ))$ is a
  $C^1$ diffeomorphism from $\Bbb R^3$ onto itself.
  \item[(ii)] The vorticity satisfies $\Omega=$curl $V\neq 0$, and there exists $p_1 >0$ such that
 $\Omega \in L^p (\Bbb R^3
  )$ for all $p\in (0, p_1)$.
  \end{itemize}
\end{theorem}
We note that  the condition (ii) is satisfied, for example, if $
\Omega\in L^{1}_{loc}(\Bbb R^3)$ and there exist constants $R, K$
and $\vare_1, \vare_2
>0$ such that $ |\Omega(x)|\leq K e^{-\vare_1 |x|^{\vare_2}}$ for $|x|>R$,
then we have
 $\Omega\in L^p (\Bbb R^3;\Bbb R^3)$ for all $p\in (0, 1)$.
 Indeed,
for all $p\in (0, 1)$, we have
 \bqn
\int_{\Bbb R^3} |\Omega(x)|^p dx
&=& \int_{|x|\leq R} |\Omega(x)|^p dx +\int_{|x|>R}  |\Omega (x)|^p \,dx\\
&\leq&|B_R |^{1-p}\left(\int_{|x|\leq R} |\Omega(x)| dx\right)^{p} +
K^p \int_{\Bbb R^3} e^{-p\vare_1|x|^{\vare_2}}dx <\infty ,
 \eqn
where $|B_R|$ is the volume of the ball $B_R$ of radius $R$.\\
In the zero vorticity case $\Omega =0$, from  div $V=0$ and  curl
$V=0$,
  we have $V=\nabla h$,  where $h(x)$ is a harmonic function in $\Bbb R^3$. Hence,
  we have an easy example of self-similar blow-up,
 $$v(x,t)=\frac{1}{(T_*-t)^{\frac{\a}{\a+1}}}\nabla
 h\left(\frac{x}{(T_*-t)^{\frac{1}{\a+1}}}\right),$$
 in $\Bbb R^3$, which is also the case for the 3D Navier-Stokes with $\alpha=1$.
 We do not consider this case in the theorem.\\

The above theorem is  actually a corollary of the following more
general theorem.

\begin{theorem}
 Let   $v\in C([0,T);C^1 (\Bbb R^3 ;\Bbb R^3))$ be a classical solution
  to the 3D Euler equations generating the particle trajectory mapping $X(\cdot ,t)$ which is a
   $C^1$ diffeomorphism from $\Bbb R^3$
  onto itself for all $t\in (0,T)$. Suppose we have
 representation of the vorticity of the solution,
 by
 \bb\label{thm13}
 \o (x,t) = \Psi(t)\Omega (\Phi(t)x) \qquad \forall t\in [0, T)
 \ee
 where $\Psi(\cdot)\in C([0, T );(0, \infty))$,
 $\Phi(\cdot)\in C([0, T );\Bbb R^{3\times 3})$
 with $\mathrm{det}(\Phi(t))\neq 0$ on $[0, T)$;  $\Omega = \mathrm{curl}\, V$ for some
 $V$, and   there exists $p_1 >0$ such that
 $|\Omega |$ belongs to $ L^p
 (\Bbb R^3)$ for all $p\in (0, p_1 )$.  Then,
necessarily
 either $\mathrm{det}(\Phi(t))\equiv \mathrm{det}(\Phi(0))$ on $[0, T)$, or
 $\Omega=0$.
\end{theorem}

For the detailed proof of Theorem 3.3 and 3.4 we refer \cite{cha13}.

\subsection{Asymptotic self-similar blow-up}

In this subsection we consider the possibility of more refined
scenario of self-similar singularity than in the previous
subsection, called the asymptotic self-similar singularity. This
means that the local in time smooth solution evolves into a
self-similar profile as the possible singularity time is approached.
The similar notion was considered previously by Giga and Kohn in
their study of semilinear heat equation(\cite{gig}). Their sense of
convergence of solution to the self-similar profile is  the
pointwise sense with a time difference weight to make it scaling
invariant, and cannot apply directly to the case of Euler system. It
is found in \cite{cha14} that if we make the sense of convergence
strong enough, then we can apply the notion of asymptotic
self-similar singularity to the Euler and the Navier-Stokes
equations. The following theorem is proved in \cite{cha14}.

\begin{theorem}
Let $v\in C([0, T);B^{\frac{3}{p}+1}_{p, 1}(\Bbb R^3))$  be a
classical solution to the 3D Euler equations.
 Suppose there exist $p_1>0$, $\a > -1$, $\bar{V}\in C^1 (\Bbb R^3)$
with $ \lim_{R\to \infty} \sup_{|x|=R} |\bar{V}(x)|=0$
 such
that $\bar{\Omega}=$curl $\bar{V}\in L^q(\Bbb R^3)$ for all $q\in
(0,p_1)$, and the following convergence holds true:
$$
 \lim_{t\nearrow T}
(T-t)^{\frac{\a-3}{\a+1}}\left\|v(\cdot, t)
-\frac{1}{(T-t)^{\frac{\a}{\a+1}}} \bar{V}
\left(\frac{\cdot}{(T-t)^{\frac{1}{\a+1}}} \right)\right\|_{L^1 }
=0,
$$
  and
 $$
   \sup_{t\in (0,T)}
(T-t)\left\|\o(\cdot, t) -\frac{1}{T-t} \bar{\O}
\left(\frac{\cdot}{(T-t)^{\frac{1}{\a+1}}}
\right)\right\|_{\dot{B}^0_{\infty, 1} }<\infty.
$$
 Then, $\bar{\O}=0$, and $v(x,t)$ can be extended to a solution of the 3D Euler system in $[0, T+\delta]\times \Bbb R^3$,
 and belongs to $C([0, T+\delta]; B^{\frac{3}{p}+1}_{p, 1}(\Bbb R^3))$ for some $\delta
 >0$.
\end{theorem}
We note that the above theorem still
 does not exclude the possibility that the sense of vorticity convergence to the
 asymptotically self-similar singularity is weaker than $L^\infty$
 sense. Namely, a self-similar vorticity profile could be approached
 from a local classical solution in the pointwise sense in space,
  or in the $L^p(\Bbb R^3)$ sense for
some $p$  with $1\leq p<\infty$. In \cite{cha14} we also proved
nonexistence of
 asymptotic self-similar solution to the 3D Navier-Stokes equations with appropriate
 change of functional setting(see also \cite{hou2}
 for related results).\\
\ \\
The proof of the above theorem follows without difficulty from the
following blow-up rate estimate(\cite{cha14}), which is interesting
in itself.

\begin{theorem}
Let $p\in [1, \infty)$ and $v\in C([0, T);B^{\frac{3}{p}+1}_{p,
1}(\Bbb R^3))$ be a classical solution to the 3D Euler equations.
There exists an absolute constant $\eta >0$ such that if
 \bb\label{th11}
 \inf_{0\leq t<T} (T-t) \|\o (t)\|_{\dot{B}^0_{\infty, 1} } <\eta
 ,
 \ee
 then $v(x,t)$ can be extended to a solution of the 3D Euler system in $[0, T+\delta]\times \Bbb R^3$,
 and belongs to $C([0, T+\delta]; B^{\frac{3}{p}+1}_{p,
1}(\Bbb R^3))$ for some $\delta
 >0$.
\end{theorem}
We note that the proof of the local existence for $v_0 \in
B^{\frac{3}{p}+1}_{p, 1}(\Bbb R^3)$  is
 done  in \cite{cha1, cha6}(see also \cite{vis1}).
  The above theorem implies that
 if $T_*$ is the first time of singularity, then  we have
 the lower estimate of the blow-up rate,
 \bb\label{blow}
 \|\o (t)\|_{\dot{B}^0_{\infty, 1} } \geq \frac{C}{T_*-t}\quad
 \forall t\in [0, T_*)
 \ee
 for an absolute constant $C$. The estimate (\ref{blow}) was actually
 derived previously by a different argument in \cite{cha4}.
 We observe that (\ref{blow}) is consistent  with  both of the
 BKM criterion(\cite{bea}) and
the Kerr's numerical calculation in \cite{ker1} respectively.\\

 The above continuation principle
 for a local solutions in $B^{\frac{3}{p}+1}_{p,
1}(\Bbb R^3)$ has  obvious applications to the  solutions belonging
to more conventional function spaces, due to the embeddings,
$$ H^m (\Bbb R^3 )\hookrightarrow C^{1, \gamma} (\Bbb R^3)\hookrightarrow
 B^{\frac{3}{p}+1}_{p,
1}(\Bbb R^3) $$
 for $m>5/2$ and  $\gamma = m-5/2$.
For example the local solution $v\in C([0, T); H^m (\Bbb R^3 ))$ can
be continued to be $v\in C([0, T+\delta];  H^m (\Bbb R^3 ))$ for
some $\delta$, if
 (\ref{th11}) is satisfied. Regarding other implication of the above theorem on the self-similar
blowing up solution to the 3D Euler equations, we have the following
corollary(see \cite{cha14} for the proof).

\begin{cor}
Let $v\in C([0, T_*);B^{\frac{3}{p}+1}_{p, 1}(\Bbb R^3))$  be a
classical solution to the 3D Euler equations. There exists $\eta
>0$ such that if we have representation for the
velocity by (\ref{vel}),  and $\bar{\O}$=curl $\bar{V}$ satisfies
  $\|\bar{\O}\|_{\dot{B}^0_{\infty ,1}} <\eta$, then $\bar{\Omega}=0$, and
   $v(x,t)$ can be extended to a solution of the 3D Euler system in $[0, T_*+\delta]\times \Bbb R^3$,
 and belongs to $C([0, T_*+\delta]; B^{\frac{3}{p}+1}_{p,
1}(\Bbb R^3))$ for some $\delta
 >0$.
\end{cor}

\section{Model problems}
 \setcounter{equation}{0}

Since the blow-up problem of the 3D Euler equations looks beyond the
 capability of current analysis, people proposed simplified model
 equations to get insight on the original problem.
 In this section we review some of them.
 Besides those results presented  in the following subsections there
 are also studies on the other model problems.
 In \cite{din} Dinaburg,  Posvyanskii and  Sinai analyzed a quasi-linear
 approximation of the infinite system of ODE arising when we write the
 Euler equation in a Fourier mode.  Fridlander and
 Pavlovi\'{c}(\cite{fri2}) considered  a vector model, and
  Katz and Pavlovi\'{c}(\cite{katz}) studied dyadic model, both of
  which are resulted from the representation of the Euler equations
  in the wave number space. Okamoto and Ohkitani proposed
  model equations in \cite{oka},  and a `dual' system to the Euler
  equations was considered in \cite{cha6a}.

\subsection{Distortions of the Euler equations}

Taking trace of the matrix equation (\ref{basfor1} )  for $V$, we
obtain
 $\Delta p=-tr V^2$, and hence the hessian of the pressure is given
 by
 $$  P_{ij}=- \partial_i \partial _j (\Delta )^{-1} tr V^2 =-R_i R_j tr
 V^2,$$
 where $R_j$ denotes the Riesz transform(see Section 1).
 Hence we can rewrite the Euler equations as
 \bb\label{matrixeuler}
 \frac{D V}{Dt}=-V^2 -R[tr V^2], \quad R[\cdot]=(R_iR_j[\cdot])
 \ee
In \cite{con1} Constantin studied a distorted version of the above
  system,
 \bb\label{coneq}
 \frac{\partial V}{\partial t}
 =-V^2 -R[tr V^2], \quad R[\cdot]=(R_iR_j[\cdot]),
\ee
  where the convection term of the original Euler equations is
  deleted, and showed that a solution indeed blows up in finite time.
  Note that the incompressibility condition, tr$V=0$, is
  respected in the system (\ref{coneq}).
  Thus we find that the convection term should have significant role
 in the study of the blow-up problem of the original Euler
 equations.\\
On the other hand, in \cite{liu} Liu and Tadmor studied another
 distorted version of (\ref{matrixeuler}), called the restricted Euler
equations,
  \bb\label{liueq}
   \frac{D V}{Dt}=-V^2 +\frac{1}{n} (tr V^2 ) I.
   \ee
   We observe that in (\ref{liueq}) the convection term is kept, but
   the non-local operator $R_iR_j(\cdot)$ is changed into a
   local one $-1/n \delta_{ij}$, where the numerical factor $-1/n$ is to keep the
   incompressibility condition.
  Analyzing the dynamics of eigenvalues of the matrix $V$, they
    showed that the system (\ref{liueq}) actually blows up in finite
    time(\cite{liu}).

  \subsection{The Constantin-Lax-Majda equation}

In 1985 Constantin, Lax and Majda constructed a one dimensional
model of the vorticity formulation of the 3D Euler equations, which
preserves the feature of nonlocality in vortex stretching term.
Remarkably enough this model equation has an explicit solution for
general initial data(\cite{con11}).  In this subsection we briefly
review their result.
  We first observe from section 1 that vorticity  formulation of the Euler equations is
 $ \frac{D \omega }{Dt}  = S \omega ,$ where
$S=\mathcal{P} (\omega )$ defines a singular integral operator of
the Calderon-Zygmund type on $\o$.
 Let us replace
 ${ \omega (x,t)\Rightarrow \theta (x,t)}$,
 ${ \frac{D}{Dt}\Rightarrow\frac{\partial}{\partial t}}$,
  ${\mathcal{P}(\cdot
 )\Rightarrow H (\cdot )}$, where
  $\theta (x,t)$ is a scalar
 function on $\Bbb R \times \Bbb R_+$, and
 $H(\cdot)$ is the \textit{Hilbert transform} defined by
 $$Hf (x)=\frac{1}{\pi} p.v. \int_{\Bbb R}
 \frac{f(y)}{x-y}dy.
 $$
 Then we obtain, the following 1D scalar equation from the 3D Euler
 equation,
 $$ (CLM): \frac{\partial \theta}{\partial t}=\theta H\theta .$$
This model preserve the feature of  { nonlocality} of the Euler
system (E), in contrast to the more traditional one dimensional
model, the  inviscid Burgers equation.
 We recall the identities for the Hilbert transform:
\bb\label{hilbert}
 H(Hf)=-f,\quad
 H(fHg+gHf)=(Hf)(Hg)-fg,
\ee
  which imply $ H(\theta H\theta )=\frac12 [ (H\theta)^2 -\theta^2
]$.
 Applying $H$ on both sides of the first equation of (CLM), and using the formula
 (\ref{hilbert}), we
obtain
$$
(CLM)^* : (H\theta )_{t}+\frac{1}{2}((H\theta )^{2}-(\theta
)^{2})=0.
$$
 We introduce the complex valued function,
$$ { z(x,t)=H\theta
(x,t)+i\theta (x,t)}.
$$
  Then, (CLM) and (CLM)$^*$ are the imaginary and the real parts
of the complex Riccati equation,
$$
z_{t} (x,t) =\frac12 z^2 (x,t)
$$
 The explicit solution to the complex equation is
 $$
 z(x,t)=\frac{z_0}{1-\frac12 t z_0 (x)} .
 $$
 Taking the imaginary part, we obtain
 $$
 \theta (x,t) = \frac{4\theta_0 (x)}{(2-tH\theta_0 (x))^2 +t^2
 \theta_0 ^2(x)}.
 $$
 The finite time blow-up occurs if and only if
 $$
 {Z=\{ x  \, | \, \theta_0 (x)=0 \, \mbox{and}\,
 H\theta_0(x) >0 \}\neq \varnothing }.
 $$
In \cite{scho} Schochet find that even if we add viscosity term to
(CLM) there is a finite time blow-up. See also \cite{sak1, sak2} for
the studies of other variations of (CLM).

 \subsection{The 2D quasi-geostrophic equation and its 1D model}

 The 2D quasi-geostrophic
equation(QG) models the dynamics of the mixture of cold and hot air
and the fronts between them.
$$
(QG) \left\{\aligned
& \theta _{t}+(u\cdot \nabla )\theta =0, \\
&v=-\nabla ^{\bot }(-\Delta )^{-\frac{1}{2}}\theta ,\\
&\theta (x,0)=\theta _{0}(x),
\endaligned
\right.
 $$
 where $\nabla ^{\bot }=(-\partial _{2},\partial _{1})$.
Here, $\theta (x,t)$ represents the temperature of the air at
$(x,t)\in \Bbb R^2 \times \Bbb R_+$. Besides its direct physical
significance (QG) has another important feature that it has very
similar structure to the 3D Euler equations. Indeed, taking $\nabla
^\bot$ to (QG), we obtain
 $$
\left(\frac{\partial }{\partial t} +v\cdot \nabla\right){ \nabla
^\bot \theta}
  =({\nabla ^\bot \theta }\cdot
  \nabla ) v,
  $$
  where
  $$
v(x,t)=\int_{\Bbb R^2} \frac{ \nabla ^\bot \theta (y,t)}{|x-y|} dy
  .
  $$
 This is exactly the vorticity formulation of 3D Euler
 equation if we identify
 $$ \nabla^\bot \theta \Longleftrightarrow \omega
 $$
After first observation and pioneering analysis of these feature by
Constantin, Majda and Tabak(\cite{con12}) there have been so many
research papers devoted to the study of this equation(also the
equation with the viscosity term, $-(-\Delta )^{\alpha}\theta$,
$\alpha
>0$, added)(\cite{con2,con6a, con12a, con12b, cor2, cor1, cor3, cor5,
cor6, cor7, cor8, wu1,wu2,wu3, cha2a, cha10, cha10a, cha24a, ohk,
din, kis, caff}). We briefly review some of them
here concentrating on the inviscid equation (QG). \\
The local existence can be proved easily by standard method. The BKM
type of blow-up criterion proved by Constantin, Majda and Tabak in
\cite{con12} is
 \bb\label{cmtcri}
 \lim\sup_{t\nearrow T_*} \|\theta(t)\|_{H^s }=\infty
\quad
 \mbox{if and only if}\quad
\int_0 ^{T_*} \|{ \nabla ^\bot \theta  (s)}\|_{L^\infty} ds =\infty.
\ee
 This criterion has been refined, using the Triebel-Lizorkin spaces
 \cite{cha2a}.
 The question of finite time singularity/global regularity is still
 open. Similarly to the Euler equations case we also have the following
geometric type of condition for the regularity. We define the
direction field $\xi =\nao \t /|\nao \t|$ whenever $|\nao \theta
(x,t)|\neq 0$.
 \begin{definition}
 We say that
a set $\Omega_0$ is smoothly directed if there exists $\rho >0$ such
that
$$
\sup_{q\in \Omega_0} \int_0 ^T |v(X(q,t),t)|^2 dt < infty,
$$
$$
\sup_{q\in \Omega^*_0} \int_0 ^T \|\nabla \xi (\cdot ,
t)\|_{L^\infty (B_\rho (X(q,t)))} ^2 dt <\infty,
$$
where $B_\rho (X)$ is the ball of radius $\rho$ centered at $X$ and
$$ \Omega ^*_0 =\{ q \in \Omega_0 \, ;\, |\nabla \theta_0 (q)| \neq
0\}.$$
 \end{definition}
We denote $ \mathfrak{O}_T(\Omega_0 )=\{ (x,t) \, |\, x\in
X(\Omega_0 ,t), 0\leq t\leq T\}$. Then, the following theorem is
proved(\cite{con12}).
\begin{theorem}
Assume that $\Omega_0$ is smoothly directed. Then
$$ \sup_{(x,t) \in \mathfrak{O}_T(\Omega_0 )} |\nabla \theta (x,t)|<\infty,$$
and no singularity occurs in $\mathfrak{O}_T(\Omega_0 )$.
\end{theorem}

Next we present an `interpolated' result between the criterion
(\ref{cmtcri}) and  Theorem 4.1, obtained in \cite{cha10}. Let us
denote bellow,
$$
 D_0 =\{ x\in \Bbb R^2 \, |\, |\nao\th_0 (x)| \neq 0\},\quad
 D_t = X(D_0 ,t).
 $$
The following theorem(\cite{cha10}) could be also considered as the
(QG) version of Theorem 2.9.
\begin{theorem}
Let $\th (x,t)$ be the local classical solution to  $(QG)$ with
initial data $\theta_0 \in H^m (\Bbb R^2)$,
 $m>3/2$, for which $D_0
\neq\emptyset$. Let $\xi (x,t)=\nao\th (x,t)/ |\nao\th (x,t)| $ be
the direction field defined for $x\in D_t $.
 Then, the  solution can be continued up to $T <\infty$ as the classical solution, namely
$\th (t)\in C([0,T];{H^m} (\Bbb R^2 ))$,
 if there exist parameters
  $p,p', q,q', s, r_1,r_2,r_3$ satisfying the following
conditions,
 \bb\label{31}
  \frac{1}{p} +\frac{1}{p'} =1, \qquad  \frac{1}{q} +\frac{1}{q'} =1,
  \ee
  and
  \bb\label{32}
\frac{1}{r_1}+\frac{p'}{r_2}\left(1-\frac{sq'}{2} \right)+
 \frac{1}{r_3} \left\{ 1-p' \left( 1-\frac{sq'}{2} \right)\right\}
 =1
 \ee
 with
 \bb\label{33}
 0<s<1, \qquad 1\leq \frac{2}{sq'} <p \leq \infty, \qquad  \quad 1\leq
 q\leq\infty,
 \ee
 and
 \bb\label{34}
 r_1\in [1, \infty],\,  r_2 \in \left[p'\left(1-\frac{sq'}{2}
 \right),\infty\right],\,  r_3 \in \left[1-p'\left(1-\frac{sq'}{2}
 \right), \infty\right]
\ee
 such that the followings hold:
\bb\label{35}
 \int_0 ^{T} \|\xi (t )\|^{r_1}_{ \mathcal{\dot{F}}^s_{\infty,q} (D_t) } dt <\infty,
 \ee
 and
 \bb\label{36}
  \int_0 ^{T}\|\nao \th(t) \|^{r_2}_{L^{pq'}(D_t)} dt +
   \int_0 ^{T}\|\nao\th (t) \|^{r_3}_{L^{q'}(D_t)} dt <\infty.
 \ee
\end{theorem}

In order to compare this theorem with the Constantin-Majda-Tabak
 criterion (\ref{cmtcri}), let
us consider the case of $p=\infty, q=1$. In this case the conditions
(\ref{35})-(\ref{36}) are satisfied if
 \bb\label{37}
 \xi (x,t )\in L^{r_1}(0,T; C^{ s} (\Bbb R^2 )),
 \ee
 \bb\label{38}
  |\nao\th
(x,t)| \in L^{r_2} (0, T; L^{\infty} (\Bbb R^2 ))\cap L^{r_3} (0, T;
L^{\infty} (\Bbb R^2 )).
 \ee
with
 $$
 \frac{1}{r_1}+\frac{1}{r_2}\left(1-\frac{s}{2} \right)+
 \frac{s}{2r_3}
 =1.
 $$
If we {\em formally} passing $s\to 0$, and choosing
 $r_1=\infty, r_2=r_3
 =1$, we find that the conditions (\ref{37})-(\ref{38})
  are satisfied if the
 Constantin-Majda-Tabak condition in (\ref{cmtcri}) holds, since the condition
 $$\xi (x,t )
 \in L^{\infty}(0,T; C^{0} (\Bbb
 R^2))\cong L^\infty ((0,T)\times \Bbb R^2 )$$
 is automatically satisfied.
 The other is the case $q'=2/s$, where
(\ref{35})-(\ref{36}) are satisfied if
 \bb
 \xi (x,t )\in L^{r_1}(0,T;
 \mathcal{\dot{F}}^s_{\infty,\frac{2}{2-s}}(\Bbb R^2 )
 ),\quad
  |\nao\th
(x,t)| \in L^{r_2} (0, T; L^{\frac{2}{s}} (\Bbb R^2 ))
 \ee
 with $1/r_1+1/r_2=1,$ which shows mutual compensation of the
 regularity  of the direction field $\xi (x,t)$ and the
 integrability of the
 magnitude of gradient $|\nao \theta (x,t)|$ to obtain smoothness
 of $\theta (x,t)$.\\
\ \\
There had been a conjectured scenario of singularity in (QG) in the
form of hyperbolic saddle collapse of level curves of $\theta
(x,t)$(see Definition 3.3). This was excluded by C\'{o}rdoba in
1998(\cite{cor2}, see also  Section 3.2 of this article). Another
 scenario of singularity, the sharp front singularity,
 which is a two dimensional version of
  potato chip singularity(see Definition 3.2 with $n=2$)
  was excluded by C\'{o}rdoba and Fefferman in \cite{cor6} under
  the assumption of suitable velocity control(see Section 3.2).\\
\ \\
We can also consider the possibility of self-similar singularity for
(QG). We first note that (QG) has the scaling property that
  if $\t$ is a
solution of the system, then for any $\lambda >0$ and $\alpha \in
\Bbb R $ the functions
 \bb
 \label{self3a}
  \t^{\lambda, \alpha}(x,t)=\lambda ^\alpha \t(\lambda x, \l^{\a +1}
  t)
  \ee
  are also solutions of (QG) with the initial data
  $ \t ^{\lambda, \alpha}_0(x)=\lambda ^\alpha \t_0
   (\lambda x)$.
Hence, the self-similar blowing up solution should be of the
  form,
  \bb
  \label{self3b}
 \t(x, t)=\frac{1}{(T_*-t)^{\frac{\a}{\a+1}}}
\Theta\left(\frac{x}{(T_*-t)^{\frac{1}{\a+1}}}\right)
 \ee
 for $t$ sufficiently close $T_*$ and $\a \neq -1$.
 The following theorem is proved in \cite{cha13}.
\begin{theorem}
Let $v $ generates a particle trajectory, which is a $C^1$
diffeomorphism from $\Bbb R^2$ onto itself for all $t\in (0, T_*)$.
There exists no nontrivial solution $ \t  $ to the system $(QG)$ of
the form (\ref{self3b}),  if there exists $p_1, p_2 \in (0, \infty
]$, $p_1 < p_2$, such that $\Theta\in L^{p_1} (\Bbb R^2 ) \cap
L^{p_2} (\Bbb R^2 )$.
\end{theorem}
We note that the integrability condition on the self-similar
representation function $\Theta$ in the above theorem is `milder'
than the case of the exclusion of self-similar Euler equations in
Theorem 3.3, in the sense that the decay condition is of $\Theta$
(not $\nao \Theta$) near infinity is weaker than that of
$\Omega=$curl $V$.\\
\ \\
 In the remained part of
this subsection we discuss a 1D model of the 2D quasi-geostrophic
equation studied in \cite{cha19}(see \cite{mor} for related
results). The construction of the one dimensional model can be done
similarly to the Constantin-Lax-Majda equation introduced in section
4.2. We
 first note that
$$
v=-R^{\bot }\theta =(-R_2 \theta , R_1 \theta ),
$$
where $R_{j}$, $j=1,2$, is the two dimensional Riesz transform(see
Section 1).
 We can rewrite the dynamical equation of (QG) as
 $$
 \theta _{t}+ \rm{div }\, [(R^{\bot }\theta )\theta ]=0,
 $$
  since $\rm{div} (R^{\bot }\theta )=0$.
 To construct the one dimensional model we  replace:
$${ R^{\bot
}(\cdot ) \Rightarrow  H(\cdot ), \qquad  \rm{div} (\cdot)
\Rightarrow
\partial_x}
$$
to obtain
$$\theta _{t}+(H(\theta )\theta )_{x}=0. $$
 Defining  the complex valued function
 ${ z(x,t)=H\theta (x,t)
 +i\theta (x,t)}$, and following Constantin-Lax-Majda(\cite{con11}), we
 find that our equation is the imaginary part of
$$z_t +zz_x =0, $$
which is complex Burgers'
 equation. The characteristics method does not work here.
Even in that case we can show that  the finite time blow-up occurs
for the generic initial data as follows.
 \begin{theorem}
Given a periodic non-constant initial data $\theta_0\in
C^1([-\pi,\pi])$ such that $\int^{\pi}_{-\pi}\theta_0(x) dx =0$,
there is no $C^1([-\pi,\pi]\times[0,\infty))$ periodic solution to
the model equation.
\end{theorem}
For the proof we refer \cite{cha19}. Here we give a brief outline of
 the construction of an explicit blowing up solution. We
begin with the complex Burgers equation:
$$
{z_{t}+zz_{x}=0, \quad z=u+i\theta}
$$
with  $u(x,t) \equiv H\theta(x,t)$.  Expanding it to  real and
imaginary parts, we obtain the system:
$$
\left\{ \aligned
u_{t}+uu_{x}-\theta \theta _{x} &=0, \\
\theta _{t}+u\theta _{x}+\theta u_{x} &=0
\endaligned
\right.
$$
In order to perform the \textit{hodograph transform} we consider ${
x(u,\theta )}$ and ${ t(u,\theta )}$
 We have,
 \bqn
u_{x} &=&Jt_{\theta }\;, \quad
\theta _{x} =-Jt_{u}\;, \\
u_{t} &=&-Jx_{\theta }\;, \quad \theta _{t} =Jx_{u}\;,
 \eqn where $J=(x_{u}t_{\theta }-x_{\theta
}t_{u})^{-1}$.
 By  direct substitution we obtain,
 $$
\left\{ \aligned -x_{\theta }+ut_{\theta }+\theta t_{u} &=0\;,  \\
x_{u}-ut_{u}+\theta t_{\theta } &=0\; \endaligned \right. $$
 as far
as $J^{-1}\neq 0$.  This system can be written more compactly in the
form:
 \bqn
-(x-tu)_{\theta }+(t\theta )_{u} &=&0\;, \\
(x-tu)_{u}+(t\theta )_{\theta } &=&0\;,
 \eqn
 which leads to the
following Cauchy-Riemann system,
$$
\xi _{u} =\eta _{\theta },\quad \xi _{\theta } =-\eta _{u},
$$
 where we set
  $
 \eta (u,\theta ):=
x(u,\theta )-t(u,\theta )u, \quad  \xi (u,\theta ):= t(u,\theta
)\theta . $
 Hence, $f(z)=\xi (u,\theta )+i\eta (u,\theta )$ with
$z=u+i\theta $ is an analytic function. Choosing $f(z)=\log\, z, $
we find,
 \bb\label{1drel}
 t\theta =\log \sqrt{u^{2}+\theta ^{2}} , \quad
 x-tu =\arctan \frac{\theta }{u},
\ee
 which corresponds  to the initial data,$ z(x,0)=\cos x+i\sin x .$
The relation (\ref{1drel}) defines implicitly the real and imaginary
parts $(u(x,t),\theta (x,t))$ of the solution.
 Removing  $\theta$ from the system, we obtain
$$
tu\tan (x-tu)=\log\left| \frac{u}{\cos (x-tu)}\right|,
$$
 which defines $u(x,t)$ implicitly.
By elementary computations we find  both $u_x$ and $\theta_x$ blow
up at  $t=e^{-1}$.

\subsection{The 2D Boussinesq system and Moffat's problem}

The 2D Boussineq system for the incompressible fluid flows in $\Bbb
R^2$ is
 $$
 (B)_{\nu, \kappa }\left\{
 \aligned
 &\frac{\partial v}{\partial t} +(v\cdot \nabla )v =-\nabla p  +\nu \Delta  v+
 \theta {e_2} ,
  \\
&\frac{\partial \theta}{\partial t} +(v\cdot \nabla )\theta =\kappa \Delta \theta\\
 &\textrm{div }\, v =0 ,\\
 &v(x,0)=v_0 (x), \qquad \theta(x,0)=\theta_0 (x),
 \endaligned
 \right.
 $$
where $v=(v_1, v_2)$, $v_j =v_j (x,t), j=1,2$, $(x,t)\in \Bbb R^2
\times (0, \infty)$, is the velocity vector field, $p=p(x,t)$ is the
scalar pressure, $\theta(x,t)$ is the scalar temperature, $\nu \geq
0$ is the viscosity, and $\kappa \geq 0$ is the thermal diffusivity,
and ${e_2} =(0,1)$. The Boussinesq system has important roles in the
atmospheric sciences(see e.g. \cite{maj4}). The global in time
regularity of $(B)_{\nu, \kappa }$ with $ \nu >0$ and $\kappa
>0$ is well-known(see e.g. \cite{can}). On the other hand, the
regularity/singularity questions of the fully inviscid case of
$(B)_{0,0}$ is an outstanding open problem in the mathematical fluid
mechanics. It is well-known that inviscid 2D Boussinesq system has
exactly same structure to the axisymmetric 3D Euler system {\em off}
the axis of symmetry(see e.g. \cite{maj2} for this observation).
This is why the inviscid 2D Boussinesq system can be considered as a
model equation of the 3D Euler system.  The problem of the finite
time blow-up for the fully inviscid Boussinesq system is an
outstanding open problem. The BKM type of blow-up criterion,
however, can be obtained without difficulty(see \cite{cha22, cha23,
e, tan} for various forms of blow-up criteria for the Boussinesq
system)). We first consider the partially viscous cases, i.e. either
the zero diffusivity case, $\kappa =0$ and $\nu>0$, or the zero
viscosity case, $\kappa  >0$ and $\nu=0$. Even the regularity
problem for partial viscosity cases has been open recently.
Actually, in an
 article  appeared in 2001, M. K. Moffatt
raised a question of finite time singularity in the case $\kappa =0,
\nu>0$ and its possible development in the limit $\kappa \to 0$ as
one of the 21th century problems(see the Problem no. 3 in
\cite{mof}). For this problem  Cordoba, Fefferman and De La
LLave(\cite{cor7}) proved that special type of singularities, called
`squirt singularities', is absent.  In \cite{cha11} the author
considered the both of two partial viscosity cases, and prove the
global in time regularity for both of the cases. Furthermore it is
proved that as diffusivity(viscosity) goes to zero the solutions of
$(B)_{\nu. \kappa}$ converge strongly to those of zero
diffusivity(viscosity) equations\cite{cha11}. In particular the
Problem no. 3 in \cite{mof} is solved. More precise statements of
these results are stated in Theorem 1.1 and Theorem 1.2 below.

\begin{theorem}
Let $\nu >0$ be fixed, and div $v_0=0$.
 Let $m >2$ be an integer, and
 $(v_0 , \theta_0 )\in H^m( \Bbb R^2
 )$. Then, there
 exists unique solution $(v, \th )$ with $\theta \in C([0, \infty );H^m (\Bbb R^2
))$ and  $ v\in  C([0, \infty );H^m (\Bbb R^2 ))\cap
L^2(0,T;H^{m+1}(\Bbb R^2 ))$  of the system $(B)_{\nu ,0}$.
  Moreover,  for
each $s<m$, the solutions $(v, \theta )$ of $(B)_{\nu ,\kappa}$
converge to the corresponding solutions of $(B)_{\nu ,0}$ in $C
([0,T];H^{s} (\Bbb R^2 ))$ as $\kappa \to 0$.
\end{theorem}
We note that Hou and Li also obtained the existence part of the
above theorem independently in \cite{hou1}. The following theorem is
concerned with zero viscosity problem with fixed positive
diffusivity.

 \begin{theorem}
Let $\kappa >0$ be fixed, and div $v_0=0$.
  Let $m >2$ be an integer. Let $m >2$ be an integer, and
 $(v_0 , \theta_0 )\in H^m( \Bbb R^2
 )$.  Then, there
exists unique solutions  $(v, \th)$ with $v\in C([0, \infty );H^m
(\Bbb R^2 ))$ and $ \theta\in  C([0, \infty );H^m (\Bbb R^2 ))\cap
L^2(0,T;H^{m+1}(\Bbb R^2 ))$  of the system $(B)_{0 ,\kappa}$.
 Moreover,  for
each $s<m$, the solutions $(v, \theta )$ of $(B)_{\nu,\kappa}$
converge to the corresponding solutions of $(B)_{0,\kappa}$ in $C
([0,T];H^{s} (\Bbb R^2 ))$ as $\nu \to 0$.
\end{theorem}
The proof of the above two theorems in \cite{cha11} crucially uses
the Brezis-Wainger inequality  in \cite{bre, eng}. Below we consider
the fully inviscid Boussinesq system, and show that there is no
self-similar singularities under milder decay condition near
infinity than the case of 3D Euler system. The inviscid Boussinesq
system (B)$=(B)_{0 ,0}$ has scaling property that
  if $(v, \t , p)$ is a
solution of the system (B), then for any $\lambda >0$ and $\alpha
\in \Bbb R$ the functions
 \bb
 \label{self2a}
  v^{\l, \a }(x,t)=\lambda ^\alpha v(\lambda x, \l^{\a +1}
  t),\quad
  \t ^{\l, \a } (x,t)= \l ^{2\a +1} \t (\lambda x, \l^{\a +1}
  t),
  \ee
  \bb
  p^{\l, \a }(x,t)= \lambda ^{2\alpha}p (\lambda x,
\l^{\a +1}
  t)
  \ee
  are also solutions of (B) with the initial data
  $$ v_0^{\l,\a}(x)=\lambda ^\alpha v_0
   (\lambda x),\quad  \t_0^{\l,\a}(x)=
  \l ^{2\a +1} \t_0(\lambda x ).
  $$
  In view of the scaling
  properties in (\ref{self2a}), the
 self-similar blowing-up
  solution $(v(x, t), \t (x,t))$ of (B)  should of the form,
  \bq
  \label{self2b1}
 v(x, t)&=&\frac{1}{(T_*-t)^{\frac{\a}{\a+1}}}
V\left(\frac{x}{(T_*-t)^{\frac{1}{\a+1}}}\right),\\
\label{self2b2}
 \t (x,t)&=&\frac{1}{(T_*-t)^{2\a+1}}
 \Theta \left(\frac{x}{(T_*-t)^{\frac{1}{\a+1}}}\right),
 \eq
 where $\alpha
\neq -1 $. We have the following nonexistence result of such type of
solution(see \cite{cha13}).

\begin{theorem}
Let $v $ generates a particle trajectory, which is a $C^1$
diffeomorphism from $\Bbb R^2$ onto itself for all $t\in (0, T_*)$.
There exists no nontrivial solution $(v, \t )$ of the system $(B)$
of the form (\ref{self2b1})-(\ref{self2b2}), if there exists $p_1,
p_2 \in (0, \infty]$, $p_1 < p_2$, such that $\Theta\in L^{p_1}
(\Bbb R^2 ) \cap L^{p_2} (\Bbb R^2 )$, and $V\in H^m (\Bbb R^2 )$,
$m>2$.
\end{theorem}
Recalling  the fact that the system $(B)$ has the similar form as
the axisymmetric 3D Euler system, we can also deduce the
nonexistence of self-similar blowing up solution to the axisymmetric
3D Euler equations  of the fome (\ref{vel}), if $\Theta = rV^\theta$
satisfies the condition of Theorem 4.7, and curl $V=\Omega \in H^m
(\Bbb R^3 )$, $m>5/2$, where $r=\sqrt{x_1^2 +x_2 ^2}$, and
$V^\theta$ is the angular component of $V$. Note that in this case
we do not need to assume strong decay of $\Omega$ as in Theorem 3.3.
See \cite{cha13} for more details.

\subsection{Deformations of the Euler equations}

Let us consider the following system considered in \cite{cha16}.
\[
\mathrm{ (P_1)}
 \left\{ \aligned
 \frac{\partial u}{\partial t} +(u\cdot \nabla )u &=-\nabla q
 +(1+\vare )\|\nabla
u (t)\|_{L^\infty}u,
  \\
 \textrm{div }\, u =0 , &\\
 u(x,0)=u_0 (x), &
  \endaligned
  \right.
  \]
  where $u=(u_1, \cdots, u_n )$, $u_j =u_j (x, t)$, $j=1, \cdots , n,
$, is the unknown vector field $q=q(x,t)$ is the scalar, and $u_0 $
is the given initial vector field satisfying div $u_0 =0$. The
constant $\vare >0$ is fixed. Below denote curl
  $u=\o$ for `vorticity' associated the `velocity' $u$.
  We first note that the system of $(P_1)$ has the similar nonlocal structure
  to the  Euler system (E), which is implicit in the pressure term  combined
  with the divergence free condition.
Moreover it has the same scaling properties
  as the original Euler system in $(E)$. Namely, if $u(x,t)$,
  $q(x,t)$ is a pair of solutions to $(P_1)$ with initial data $u_0 (x)$,
  then for any $\alpha \in \Bbb R$
  $$ u^\lambda (x,t) =\lambda ^\alpha u (\lambda x, \lambda
  ^{\alpha +1} t ), \quad q^\lambda (x,t) =
  \lambda ^{2\alpha} q (\lambda x, \lambda
  ^{\alpha +1} t )
  $$
 is also a pair of solutions to $(P_1)$ with initial data $u_0 ^\lambda
 (x)=\lambda ^\alpha u_0 (x)$. As will be seen below, we can have the local well-posedness in the
 Sobolev space, $H^m (\Bbb R^n ), m>n/2+2$, as well as the BKM type
 of blow-up criterion for $(P_1)$, similarly to the Euler system
 (E). Furthermore, we can prove actual finite time blow-up for
 smooth initial data if $\o_0 \neq 0$. This is rather surprising in
 the viewpoint that people working on the Euler system  often have
 speculation that the divergence free condition might have the role
 of `desingularization', and might make the singularity disappear.
 Obviously this is not the case for the system $(P_1)$.
 Furthermore, there is a canonical
 functional relation between the solution of $(P_1)$ and that of
 the Euler system (E); hence the word `deformation'. Using this relation we can translate the
 blow-up condition of the Euler system in terms of the solution of
 $(P_1)$. The precise contents of the above results on $(P_1)$ are stated
 in the
 following theorem.
\begin{theorem}
Given $u_0 \in H^m (\Bbb R^n )$ with div $u_0 =0$, where $m>
\frac{n}{2}+2$, the following statements hold true for $(P_1)$.
\begin{itemize}
\item[(i)] There exists a local in time unique solution $u(t)\in
C([0,T]:H^m (\Bbb R^n ))$ with $T=T (\|u_0 \|_{H^m })$.
\item[(ii)]
The solution $u(x,t)$ blows-up at $t=t_*$, namely
 $$
 \lim\sup_{t\to t_*}\|u(t)\|_{H^m} =\infty\quad \mbox{ if and
only if}\quad \int_0 ^{t_*} \|\o (t)\|_{L^\infty} dt =\infty,
 $$
 where
$\o=$ curl $u$. Moreover, if the solution $u(x,t)$ blows up at
$t_*$, then necessarily,
 $$
 \int_0 ^{t_*} \exp\left[ (2+\vare)\int_0 ^\tau \|\nabla
 u (s)\|_{L^\infty}ds\right] d\tau =\infty
 $$
 for $n=3$, while
 $$
\int_0 ^{t_*} \exp\left[ (1+\vare)\int_0 ^\tau \|\nabla
 u (s)\|_{L^\infty}ds\right] d\tau=\infty
 $$
 for $n=2$.
\item[(iii)] If $\|\o_0
\|_{L^\infty }\neq 0$, then there exists time $ t_* \leq
\frac{1}{\vare \|\o_0 \|_{L^\infty }}$ such that solution $u(x,t)$
of $(P_1)$ actually blows  up at $t_*$. Moreover, at such $t_*$ we
have
 $$
  \int_0
^{t_*} \exp\left[ (1+\vare)\int_0 ^\tau \|\nabla
 u (s)\|_{L^\infty}ds\right] d\tau=\infty .
 $$
 \item[(iv)] The functional
relation between the  solution $u(x,t)$ of $(P_1)$ and the solution
$v(x,t)$ of the Euler system (E) is given by
$$
 u(x,t)=\varphi'(t) v(x,  \varphi (t)),
$$
 where
 $$
 \varphi (t)=\lambda\int_0 ^t\exp\left[ (1+\vare )\int_0 ^\tau \|\nabla
 u(s)\|_{L^\infty}ds\right]d\tau.
 $$
 (The relation between the two initial datum is
 $u_0 (x)=\lambda v_0 (x)$.)
 \item[(v)]
 The solution $v(x,t)$ of the Euler system (E) blows up at
 $T_* <\infty$ if and only if  for $t_*:=\varphi^{-1} (T_*) < \frac{1}{\vare \|\o_0 \|_{L^\infty}}$
 both of the followings hold true.
$$
 \int_0 ^{t_*}\exp\left[ (1+\vare )\int_0 ^\tau \|\nabla
 u(s)\|_{L^\infty}ds\right]d\tau <\infty,
 $$
 and
 $$
 \int_0 ^{t_*}\exp\left[ (2+\vare )\int_0 ^\tau \|\nabla
 u(s)\|_{L^\infty}ds\right]d\tau =\infty.
 $$
\end{itemize}
\end{theorem}
For the proof we refer \cite{cha16}. In the above theorem the result
(ii) combined with (v) shows indirectly that there is no finite time
blow-up in 2D Euler equations, consistent with the well-known
result.  Following the argument on p. 542 of \cite{cha4}, the
following fact can be
verified without difficulty:\\
 We set
\bb\label{aoft}
  a(t)=\exp \left( \int_0 ^t (1+\vare) \|\nabla
u(s)\|_{L^\infty} ds \right).
 \ee
  Then, the solution $(u, q)$ of
$(P_1)$ is given by
$$ u(x,t)=a(t) U(x,t), \quad q(x,t)=a(t)P(x,t), $$
where $(U, P)$ is a solution of the following system,
$$
(aE)\left\{ \aligned
 &\frac{\partial U}{\partial t} +a(t) (U\cdot \nabla )U =-\nabla P,\\
 &\qquad \mathrm{div}\, U=0,\\
 &U(x,0)=U_0 (x)
 \endaligned
 \right.
 $$
The system $(aE)$ was studied in \cite{cha4}, when $a(t)$ is a
prescribed function of $t$, in which case the proof of local
existence of $(aE)$ in \cite{cha4} is exactly same as the case of
$(E)$.  In the current case, however, we need an extra proof of
local existence, as is done in the next section, since the function
$a(t)$ defined by (\ref{aoft}) depends on the solution $u(x,t)$
itself. As an application of Theorem 4.8 we can prove the following
lower estimate of the possible blow-up time(see \cite{cha16} for the
detailed proof).
 \begin{theorem}
 Let $p\in (1, \infty)$ be fixed.
 Let $v(t)$ be the local classical solution of the 3D Euler
 equations with initial data $v_0\in H^m (\Bbb R^3 )$, $m>7/2$. If
 $T_*$ is the first blow-up time, then
\bb\label{loweres}
 T_* -t\geq \frac{1}{ C_0\|\o (t) \|_{\dot{B}^{\frac{3}{p}}_{p,1}}},
 \quad \forall t\in (0, T_*)
\ee
 where $C_0$ is the absolute constant in $(Q_2)$.
  \end{theorem}
 In \cite{cha4} the following form of lower estimate for the blow-up
 rate is derived.
\bb\label{loweres1}
 T_*-t \geq \frac{1}{ \tilde{C}_0\|\o (t)\|_{\dot{B}^{0}_{\infty,1}}},
\ee
 where $\tilde{C}_0$ is another absolute constant(see also the remarks after Theorem 3.6).
  Although there
is (continuous) embedding relation,
$\dot{B}^{\frac{3}{p}}_{p,1}(\Bbb R^3) \hookrightarrow
\dot{B}^{0}_{\infty,1}(\Bbb R^3)$ for $p\in [1, \infty]$(see Section
1), it is difficult to compare the two estimates (\ref{loweres}) and
(\ref{loweres1}) and decide which one is sharper, since the precise
evaluation of the optimal constants
$C_0, \tilde{C}_0$ in those inequalities could be very difficult problem.\\
\ \\
Next, given $\vare \geq 0$, we consider the following problem.
\[
\mathrm{ (P_2)}
 \left\{ \aligned
 \frac{\partial u}{\partial t} +(u\cdot \nabla )u &=-\nabla q
 -(1+\vare )\|\nabla
u (t)\|_{L^\infty}u ,
  \\
 \textrm{div }\, u =0 , &\\
 u(x,0)=u_0 (x), &
  \endaligned
  \right.
  \]
  Although the system of $(P_2)$ has also the same nonlocal structure and the scaling
  properties as the Euler system and $(P_1)$, we have the result
  of the global regularity stated in the following theorem(see \cite{cha16} for the proof).
  \begin{theorem}
Given $u_0 \in H^m (\Bbb R^n)$ with div $u_0 =0$, where $m>
\frac{n}{2}+2$, then the solution $u(x,t)$ of $(P_2)$  belongs to
$C([0, \infty ):H^m (\Bbb R^n ))$. Moreover, we have the following
decay estimate for the vorticity,
$$
 \| \o (t)\|_{L^\infty }
 \leq \frac{\|\o _0 \|_{L^\infty}}{1+\vare\|\o_0 \|_{L^\infty} t} \qquad \forall
 t\in [0, \infty ).
$$
\end{theorem}
We also note that solution of the system $(P_2)$ has also similar
functional relation with that of the Euler system as given in (iv)
of Theorem 5.8 as will be clear in the proof of Theorem 1.1 in the
next section.\\
 \ \\
Next, given $\vare >0 $, we consider the following perturbed systems
of (E).
\[
 (E)^\vare _\pm\left\{ \aligned
 &\frac{\partial u}{\partial t} +(u\cdot \nabla )u =-\nabla q
 \pm \vare \|\nabla u\|_{L^\infty}^{1+\vare}u,
  \\
 &\textrm{div}\, u =0 , \\
 & u(x,0)=u_0 (x).
  \endaligned
  \right.
  \]
If we set $\vare=0$ in the above, then the system $(E)^0 _\pm$
becomes $(E)$. For $\vare >0$  we have finite time blow-up for the
system $(E)^\vare_+$ with certain initial data, while  we have the
global regularity for $(E)^\vare_-$ with all solenoidal initial data
in $H^m (\Bbb R^3 )$, $m>5/2$. More precisely we have the following
theorem(see \cite{cha17} for the proof).
 \begin{theorem}
  \begin{itemize}
  \item[(i)] Given $\vare>0$, suppose
  $u_0 =u_0^\vare\in H^m (\Bbb R^3 )$  with div $u_0=0$ satisfies
  $\|\o_0 \|_{L^\infty} > (2/\vare )^{1/\vare}$, then
 there exists  $T_* $ such that the
  solution $u(x,t)$ to $(E)_\vare ^+$ blows up at $T_*$, namely
  $$
   \lim\sup _{t\nearrow T_*} \|u(t)\|_{H^m} =\infty.
   $$
  \item[(ii)] Given $\vare >0$ and  $u_0 \in H^m (\Bbb R^3)$ with div $u_0=0$, there exists
   unique global in time classical solution $u(t)\in
   C([0, \infty ); H^m (\Bbb R^3 ))$ to
   $(E)^\vare _-$. Moreover, we have the global in time vorticity
   estimate for the solution of $(E)^\vare _-$,
   $$
    \|\o(t)\|_{L^\infty} \leq
 \max\left\{ \|\o_0\|_{L^\infty}, \left(\frac{1}{\vare}\right)^{\frac{1}{\vare}}\right\}
  \qquad \forall t\geq 0.
   $$
   \end{itemize}
\end{theorem}
The following theorem relates the finite time blow-up/global
regularity of the Euler system with those of  the system
$(E)^\vare_\pm$.
\begin{theorem}
Given $\vare >0$, let $u^\vare_\pm$ denote the solutions of
$(E)^\vare _\pm$ respectively with the same initial data $u_0 \in
H^m (\Bbb R^3 )$, $m>5/2$. We define
$$
\varphi ^\vare_\pm (t, u_0):=\int_0 ^t \exp \left[ \pm \vare \int_0
^\tau \|\nabla u^\vare_\pm  (s)\|_{L^\infty}^{1+\vare} ds\right]
d\tau.
$$
\begin{itemize}
\item[(i)] If  $\varphi ^\vare_- (\infty , u_0 ) =\infty$, then the
solution of the Euler system with initial data $u_0$ is regular
globally in time.
\item[(ii)] Let $t_*$
be the first blow-up time  for a solution  $u^\vare_+$  of
$(E)^\vare_+$ with initial data $u_0$ such that
$$\int_0 ^{t_*} \|\o ^\vare _+(t)\|_{L^\infty} dt =\infty, \qquad \mbox{where}\quad
\o ^\vare _+ =\textrm{curl}\, \,\, u^\vare_+ .$$
 If
$\varphi ^\vare_+ (t_* , u_0 )<\infty$, then the solution of the
Euler system blows up at the finite time $T_*=\varphi ^\vare_+ (t_*
, u_0 )$.
\end{itemize}
\end{theorem}
We refer \cite{cha17} for the proof of the above theorem.

 \section{Dichotomy: singularity or global regular dynamics?}
 \setcounter{equation}{0}

In this section we review  results in \cite{cha12}. Below  $S$, $P$
and $\xi (x,t)$ are the deformation tensor , the Hessian of the
pressure and the vorticity direction field, associated with the
flow, $v$, respectively as introduced in section 1.
 Let  $\{ (\lambda _k , \eta^k)\}_{k=1}^3$ be
the eigenvalue and the normalized eigenvectors of  $S$.
 We set $\lambda =(\lambda_1, \lambda_2 ,\lambda_3 )$, and
 $$|\lambda|= \left(\sum_{k=1}^3 \lambda_k ^2 \right)^{\frac12},\quad
 \rho_k =\eta^k \cdot P \eta^k \quad\mbox{ for $k=1,2,3$}.
$$
 We also denote
$$\eta^k (x,0)=\eta^k_0 (x),\quad
\lambda_k (x,0)=\lambda_{k,0}(x),\quad \lambda (x,0)=\lambda_0 (x),
\quad \rho_k (x,0)=\rho_{k,0}(x)
$$
for the quantities at $t=0$.
  Let  $\o (x,t)\neq 0$. At such point $(x,t)$ we define the scalar fields
$$
 \a  = \xi\cdot  S\xi,\quad
 \rho = \xi \cdot P\xi.
$$
 At
the points where $\o (x,t)=0$ we define $\a(x,t)=\rho(x,t)=0$. We
denote $\a_0 (x)=\a (x,0)$, $\rho_0 (x)=\rho (x,0 )$.
 Below we denote $ f(X(a,t),t)'=\frac{Df}{Dt} (X(a,t),t)$
 for simplicity.

Now, suppose that there is no blow-up of the solution on $[0, T_*]$,
and
 the inequality
 \bb\label{contra1}
\a(X(a,t),t) |\o (X(a,t),t)|  \geq \vare |\o (X(a,t),t)|^2
 \ee
 persists on $[0,T_*]$. We will see that this leads to a
 contradiction. Combining (\ref{contra1}) with (\ref{basfor3}), we have
 $$
 |\o |' \geq \vare |\o|^2.
 $$
 Hence, by Gronwall's lemma, we obtain
 $$
 |\o (X(a,t),t)|\geq \frac{|\o_0 (a)|}{1-\vare |\o_0 (a)| t},
 $$
which implies that
 $$ \lim\sup_{t\nearrow T_*}|\o (X(a,t),t)|=\infty.$$
Thus we are lead to the following lemma.
\begin{lemma}
 Suppose $\a_0 (a) > 0$, and there exists $\vare>0$ such that
 \bb\label{210}
\a _0 (a) |\o_0 (a)|\geq \vare |\o_0 (a)|^2.
 \ee
 Let us set
 \bb
T_* =\frac{1}{\vare \a_0 (a)}.
 \ee
 Then, either  the vorticity blows up no later than $T_*$, or there
 exists $t \in (0, T_*)$ such that
 \bb\label{th11} \a(X(a,t),t) |\o (X(a,t),t)|  < \vare |\o (X(a,t),t)|^2.
  \ee
  \end{lemma}

From this lemma we can derive the following:

\begin{theorem}[vortex dynamics]
Let $v_0 \in H^m (\Omega )$, $m>5/2$, be given.  We define
$$\Phi_1(a,t)= \frac{\a (X(a,t),t)}{|\o (X(a,t),t)|}$$
and
$$\Sigma_1 (t)=\{ a\in \Omega \,|\, \a (X(a,t),t)> 0\}$$
 associated with the classical solution $v(x,t)$.
 Suppose $a\in
\Sigma_1 (0)$ and $\o_0 (a) \neq 0$. Then one of the following holds
true.
\begin{itemize}
\item[(i)] {\rm{(finite time singularity)}}
The solution of the Euler equations blows-up in finite time along
the trajectory $\{ X(a,t)\}$.
\item[(ii)] {\rm{(regular dynamics)}} On of the following holds
true:
\begin{description}
\item[(a)] {\rm{(finite time extinction of $\a$)}}
 There exists $t_1 \in (0, \infty)$ such that $\a (X(a, t_1),t_1)=0$.
\item[(b)] {\rm{(long time behavior of $\Phi_1$)}} There exists an infinite sequence $\{
t_j\}_{j=1}^\infty$ with $t_1<t_2<\cdots <t_j <t_{j+1} \to \infty$
as $j\to \infty$ such that for all $j=1,2,\cdots$ we have $\Phi_1
(a,0)> \Phi_1 (a,t_1)>\cdots>\Phi_1 (a, t_j)>\Phi_1 (a,t_{j+1})>0$
and $\Phi_1 (a,t)\geq \Phi _1 (a,t_j ) >0$ for all $t\in [0,t_j]$.
\end{description}
\end{itemize}
 \end{theorem}
 As an illustration of proofs for the Theorem 5.2 and 5.3 below, we
 give outline of the proof of the above theorem.
Let us first observe that the formula
 $$ |\o (X(a,t),t)|=\exp\left[ \int_0 ^t \a (X(a,s),s) ds\right]
 |\o_0 (a)|,
 $$
 which is obtained from (\ref{basfor3}) immediately  shows that $\o (X(a,t),t)\neq
 0$ if and only if $\o_0 (a)\neq 0$ for the particle trajectory $\{ X(a,t)\}$  of the classical
 solution $v(x,t)$ of the Euler equations. Choosing
$\vare = \a _0 (a)/|\o_0 (a)|$ in Lemma 4.1, we see that
 either the vorticity  blows up no later than $T_* =1/\a_0 (a)$,
 or there exists $t_1\in (0, T_* )$ such that
$$
 \Phi_1 (a,t_1)=\frac{\a (X(a,t_1),t_1)}{|\o (X(a,t_1),t_1 )|} < \frac{\a_0 (a)}{|\o_0
(a)|}=\Phi_1 (a,0).
$$
Under the hypothesis that (i) and (ii)-(a) do not hold true, we may
assume $a\in \Sigma_1 (t_1)$ and  repeat the above argument
 to find $t_2 >t_1$ such that $\Phi_1(a, t_2 )<\Phi_1 (a, t_1 )$,
 and also $a \in \Sigma_1 (t_2)$. Iterating the argument, we find
 a monotone increasing sequence $\{ t_j\}_{j=1}^\infty$ such that $\Phi _1 (a,t_{j})>\Phi
_1 (a,t_{j+1})$ for all $j=1,2,3,\cdots$. In particular we can
choose each $t_j$ so that $\Phi_1 (a, t)\geq \Phi_1 (a, t_j )$
 for all $t\in (t_{j-1}, t_j]$. If $t_j \to t_\infty <\infty$ as $j\to \infty$,
  then we can proceed further to have $t_* > t_\infty$ such that
$\Phi_1 (a, t_\infty )>\Phi_1 (a, t_* )$. Hence,  we may set
$t_\infty =\infty$, which finishes the proof.\\

The above argument can be extended to prove the following theorems.
\begin{theorem}[dynamics of $\a$]
Let $v_0 \in H^m (\Omega )$, $m>5/2$, be given. In case
$\a(X(a,t),t)\neq 0$ we define
$$\Phi_2(a,t)= \frac{|\xi\times S\xi
|^2(X(a,t),t)-\rho (X(a,t),t)}{\a
  ^2(X(a,t),t)},
  $$
  and
$$\Sigma_2 ^+(t)=\{ a\in \Omega \,|\, \a (X(a,t),t)> 0,  \Phi _2
(X(a,t),t)>1\},$$
$$\Sigma_2 ^-(t)=\{ a\in \Omega \,|\, \a (X(a,t),t)<0,  \Phi _2
(X(a,t),t)<1\},$$

 associated with $v(x,t)$.
 Suppose $a\in \Sigma_2 ^+(0)\cup \Sigma_2 ^-(0)$.
 Then one of the following
holds true.
\begin{itemize}
\item[(i)] {\rm{(finite time singularity)}}
 The solution of the Euler equations blows-up in finite time along
the trajectory $\{ X(a,t)\}$.
\item[(ii)] {\rm{(regular dynamics)}} One of the following holds
true:
\begin{description}
 \item[(a)] {\rm{(finite time extinction of $\a$)}} There exists $t_1 \in (0, \infty)$ such that $\a (X(a, t_1),t_1)=0$.
\item[(b)] {\rm{(long time behaviors of $\Phi_2$)}} Either there exists $T_1\in (0, \infty)$ such that
$ \Phi_2 (a,T_1) =1$, or there exists an infinite sequence $\{
t_j\}_{j=1}^\infty$ with $t_1<t_2<\cdots <t_j <t_{j+1} \to \infty$
as $j\to \infty$ such that  one of the followings hold:
 \begin{description}
 \item[(b.1)] In the case  $a\in \Sigma_2 ^+(0)$,  for
all $j=1,2,\cdots$ we have $\Phi_2 (a,0)> \Phi_2
(a,t_1)>\cdots>\Phi_2 (a, t_j)>\Phi_2 (a,t_{j+1})
>1$  and $\Phi_2 (a,t)\geq \Phi_2
(a,t_j )>1 $ for all $t\in [0,t_j]$.
\item[(b.2)]  In the case $a\in \Sigma_2 ^-(0)$,  for all $j=1,2,\cdots$ we
have $\Phi_2 (a,0)<\Phi_2 (a,t_1)<\cdots <\Phi_2 (a, t_j)< \Phi_2
(a,t_{j+1}) <1$  and $\Phi_2 (a,t)\leq \Phi_2 (a,t_j )<1 $ for all
$t\in [0,t_j]$.
\end{description}
\end{description}
\end{itemize}
 \end{theorem}

\begin{theorem}[spectral dynamics]
Let $v_0 \in H^m (\Omega )$, $m>5/2$, be given. In case
$\lambda(X(a,t),t)\neq 0$ we define
$$
\Phi_3(a,t)= \frac{\sum_{k=1}^3\left[ -\lambda_k ^3 +\frac14 |\eta_k
\times \o |^2\lambda_k - \rho_k \lambda_k \right](X(a,t),t)}{
|\lambda (X(a,t),t)|^3},
$$
and
$$\Sigma_3 (t)=\{ a\in \Omega \,|\, \lambda (X(a,t),t)\neq 0,  \Phi _3
(X(a,t),t)>0\}$$
 associated with  $v(x,t)$.
 Suppose $a\in \Sigma_3 (0)$.
 Then one of the following
holds true:
\begin{itemize}
\item[(i)] {\rm{(finite time singularity)}}
 The solution of the Euler equations blows-up in finite time along
the trajectory $\{ X(a,t)\}$.

\item[(ii)]{\rm{(regular dynamics)}} One of the followings hold
true:
 \begin{description}
 \item[(a)] {\rm{(finite time extinction of $\lambda$)}}
 There exists $t_1 \in (0, \infty)$ such that $\lambda (X(a, t_1),t_1)=0$.

\item[(b)]{\rm{(long time behavior of $\Phi_3 $)}} Either there exists $T_1\in (0, \infty)$ such that
$ \Phi_3 (a,T_1) =0$, or there exists an infinite sequence $\{
t_j\}_{j=1}^\infty$ with $t_1<t_2<\cdots <t_j <t_{j+1} \to \infty$
as $j\to \infty$ such that for all $j=1,2,\cdots$ we have $\Phi_2
(a,0)> \Phi_3 (a,t_1)>\cdots>\Phi_3 (a, t_j)>\Phi_3(a,t_{j+1})
>0$  and $\Phi_3 (a,t)\geq \Phi_3
(a,t_j )>0 $ for all $t\in [0,t_j]$.
\end{description}
\end{itemize}
 \end{theorem}
 For the details of the proof of Theorem 5.2 and Theorem 5.3 we refer \cite{cha12}.

 Here we present a refinement of Theorem 2.1 of \cite{cha15}, which is proved in
 \cite{cha12}.

\begin{theorem} Let $v_0 \in H^m (\Omega )$, $m>5/2$, be given.
For such $v_0$ let us define a set $\Sigma \subset \Omega$ by
 \bqn
 \lefteqn{ \Sigma=\{ a\in
\Omega \,|\, \a_0 (a) > 0,  \o_0 (a) \neq 0, \exists \,\,\vare \in
(0, 1) \,\, \mbox{such
that}}\hspace{.3in}\\
&&  \rho_0(a)  + 2\a ^2_0 (a)-|\xi_0 \times S_0\xi_0 |^2(a)\leq
(1-\vare)^2 \a^2 _0 (a)\}.
 \eqn
  Let us set
 \bb
 T_* =\frac{1}{\vare \a_0 (a)}.
 \ee
Then, either the solution blows up no later than $T_*$, or there
exists $t\in (0, T_*)$ such that
 \bb
  \rho(X(a,t),t) + 2\a ^2(X(a,t),t) -|\xi \times S\xi |^2(X(a,t),t)
  > (1-\vare)^2 \a ^2(X(a,t),t).
\ee
  \end{theorem}
We note that if we ignore the term $|\xi_0 \times S_0\xi_0
  |^2(a)$, then we have the condition,
$$ \rho_0(a)  + \a ^2_0 (a)\leq  (-2\vare +\vare^2) \a^2 _0 (a) < 0,
 $$
 since $\vare \in (0, 1)$.
 Thus $\Sigma \subset \mathcal{S}$, where $\mathcal{S}$ is the
 set defined in Theorem 2.1 of \cite{cha15}. One can verify without
 difficulty that $\Sigma =\emptyset$ for the  2D Euler flows.
 Regarding the question if $\Sigma \neq\emptyset$ or not for 3D Euler flows, we
have the following proposition(see \cite{cha15} for more details).
\begin{pro}
Let us consider the system  the domain $\Omega =[0,2\pi ]^3$ with
the periodic boundary condition. In $\Omega$ we consider the
Taylor-Green vortex flow defined by
 \bb\label{ta}
 u(x_1,x_2,x_3 )=(\sin x_1 \cos x_2 \cos x_3 ,
-\cos x_1 \sin x_2 \cos x_3 , 0).
 \ee
Then, the set
$$\mathcal{S}_0 =\left\{(0,
\frac{\pi}{4} , \frac{7\pi}{4} ),(0, \frac{7\pi}{4} , \frac{\pi}{4}
)\right\} $$
 is  included in $\Sigma$ of Theorem 4.4.  Moreover, for $x\in
 \mathcal{S}_0$ we have the explicit values of $\alpha$ and $\rho$,
 $$
 \alpha(x)=\frac12,\quad
 \rho(x)=-\frac12.
$$
 \end{pro}
We recall that the Taylor-Green vortex has been the first candidate
proposed for a finite time singularity for the 3D Euler equations,
and there have been many numerical calculations of solution of (E)
with the initial data given by it(see e.g. \cite{bra}).

\section{Spectral dynamics approach}
 \setcounter{equation}{0}
Spectral dynamics  approach in the fluid mechanics was initiated by
Liu and Tadmor(\cite{liu}). They analyzed the restricted Euler
system (\ref{liueq}) in terms of (poitwise) dynamics of the
eigenvalues of the velocity gradient matrix $V$. More specifically,
multiplying
 left and right eigenvectors of $V$ to (\ref{liueq}), they derived
   $$
  \frac{D \lambda_j}{Dt}=-\lambda_j^2
   +\frac{1}{n} \sum_{k=1}^n \lambda_k ^2 ,\quad j=1,2,\cdots, n,
   $$
   where  $\lambda_j$, $j=1,2, ..,n$ are eigenvalues $V$, which are not necessarily
   real values.
   In this model system they proved finite time blow-up for suitable
   initial data. In this section we review the results in
   \cite{cha8}, where the full Euler system is concerned.
Moreover, the we are working on the dynamics of eigenvalues of the
deformation tensor $S$(hence real valued), not the velocity gradient
matrix. We note that there were also application of the spectral
dynamics of the deformation tensor in the study of regularity
problem of the Navier-Stokes equations by Neustupa and
Penel(\cite{neu}).
 In this section for
simplicity we consider the 3D Euler system (E) in the periodic
domain, $\Omega =\Bbb T^3(=\Bbb R^3/\Bbb Z^3)$.
 Below we denote ${ \lo ,\lt  , \ltt }$ for the eigenvalues of
 the deformation tensor
 $S=( S_{ij})$ for the velocity fields of the 3D Euler system.
 We will first establish the following formula,
 \bb\label{spec1}
 \frac{d}{dt} \intrr (\lo^2 +\lt^2 +\ltt^2 )dx=
 -4\intrr \lo\lt \ltt \, dx,
\ee
 which  has important implications(Theorem 6.1- Theorem 6.3 below).
  Indeed, using (\ref{basfor2}), we can compute
 \bqn
\lefteqn{\frac12\frac{d}{dt}\intrr S_{ij}S_{ij}dx=\intrr S_{ij}
 \frac{D S_{ij}}{Dt}dx}\hspace{.0in} \\
 &&=-\intrr S_{ik}S_{kj}S_{ij} dx -\frac14 \intrr \o_i S_{ij}\o_j
+\frac14\intrr |\o |^2 S_{ii} dx +\intrr P_{ij} S_{ij} dx\\
 &&=-\intrr S_{ik}S_{kj}S_{ij} dx-\frac18\frac{d}{dt}\intrr |\o|^2dx,
 \eqn
where we used the summation
  convention for the repeated indices, and used the
  $L^2$-version of the vorticity equation,
\bb\label{spec2}
  \frac12\frac{d}{dt}\intrr |\o
 |^2dx=\intrr \o_i S_{ij}\o_j ,
 dx
 \ee
 which is immediate from (\ref{basfor2a}).
 We note
 \bqn
 \lefteqn{{ \intrr |\o |^2 dx} = \intrr |\nabla v|^2 dx=\intrr V_{ij}
 V_{ij}dx= \intrr (S_{ij}+A_{ij})(S_{ij}+A_{ij})dx}\hspace{.3in}\\
 &&=\intrr S_{ij}S_{ij} dx +\intrr A_{ij}A_{ij}dx
 =\intrr S_{ij}S_{ij} dx+\frac12 { \intrr |\o|^2dx}.
 \eqn
  Hence,
 $$\intr S_{ij}S_{ij} dx =\frac12 \intrr |\o |^2 dx$$
 Substituting this into (\ref{spec2}), we obtain that
$$
\frac{d}{dt}\intrr S_{ij}S_{ij}dx=-\frac43\intrr S_{ik}S_{kj}S_{ij}
dx,
$$
which, in terms of the spectrum of $S$, can be written as
 \bb\label{spec3}
\frac{d}{dt} \intrr (\lo^2 +\lt^2 +\ltt^2 )dx=-\frac43\intrr (\lo^3
+ \lt^3 + \ltt^3 ) dx.
 \ee
 We observe from the divergence free condition, $0=$div $v=Tr
 S=\lo+\lt+\ltt$,
 \bqn
 0&=&(\lo+\lt +\ltt )^3\\
 &=&\lo^3 +\lt^3 +\ltt^3 +3\lo^2 (\lt+\ltt
 )+3\lt^2 (\lo +\ltt ) +3\ltt (\lo+\lt )+ 6\lo\lt\ltt\\
 &=&\lo^3 +\lt^3 +\ltt^3 -3(\lo^3 +\lt^3 +\ltt^3)+ 6\lo\lt\ltt.
 \eqn
 Hence, $\lo^3 + \lt^3 + \ltt^3=3\lo\lt\ltt$. Substituting this into
 (\ref{spec3}), we completes the proof of (\ref{spec1}).\\
 \ \\
Using the formula (\ref{spec1}), we can first prove  the following
new a priori estimate for the $L^2$ norm of vorticity for the 3D
incompressible Euler equations(see \cite{cha8} for the proof). We
denote
 $$ \mathbb{H}^m _\s =\{ v\in [H^m (\Bbb T^3 )]^3\, |\, \mathrm{div
 }\, v=0\}.$$
\begin{theorem}
Let $v(t)\in C([0,T); \mathbb{H}^m _\s )$, $m>5/2$ be the local
classical solution of the 3D Euler equations with initial data $v_0
\in \mathbb{H}^m_\s $ with $\o_0 \neq 0$. Let $\lo (x,t)\geq \lt
(x,t)\geq \ltt (x,t)$ are the eigenvalues of the deformation tensor
$S_{ij}(v)=\frac12 ( \frac{\partial v_j}{\partial x_i}
+\frac{\partial v_i}{\partial x_j})$. We denote $\lt^+ (x,t)=\max\{
\lt (x,t), 0\}$, and $\lt^- (x,t)=\min\{ \lt (x,t), 0\}$. Then, the
following  estimates hold.
 \bqn
  \lefteqn{ \exp\left[
  \int_0 ^t \left(\frac12 \inf_{x\in \Bbb T^3}\lt^+ (x,t)-\sup_{x\in \Bbb T^3}
  |\lt^-(x,t)|\right)dt\right]\leq
  \frac{\|\o (t)\|_{L^2}}{\|\o_0\|_{L^2}}}\hspace{1.in}\n \\
  &&\leq   \exp\left[
  \int_0 ^t \left( \sup_{x\in \Bbb T^3}\lt^+ (x,t)-\frac12 \inf_{x\in \Bbb T^3}
  |\lt^- (x,t)|\right)dt\right]\n \\
\eqn
 for all $t\in (0,T)$.
\end{theorem}
The above estimate says, for example, that if we have the following
compatibility conditions,
$$\sup_{x\in \Bbb T^3}\lt^+ (x,t)\simeq  \inf_{x\in \O}|\lt^- (x,t)|\simeq  g(t)$$
for some time interval $[0,T]$, then
$$ \|\o (t)\|_{L^2}
\lesssim O\left(\exp\left[C \int_0 ^t g(s)ds\right]\right) \qquad
\forall t\in [0,T]
$$
for some constant $C$. On the other hand, we note the following
connection of the above result to the previous one. From the
equation
$$ \frac{D |\o |}{Dt}=\alpha |\o |, \qquad \alpha (x,t)=\frac{\o \cdot S\o
}{|\o |^2}$$ we immediately have
 \bqn
 \|\o (t)\|_{L^2}&\leq& \|\o_0 \|_{L^2}
 \exp \left(\int_0 ^t \sup_{x\in \Bbb T^3} \alpha (x,s)
  ds \right)\\
 &\leq&\|\o_0 \|_{L^2} \exp \left(\int_0 ^t  \sup_{x\in \Bbb T^3}
 \lambda_1(x,s )
 d\tau \right),
 \eqn
 where we used the fact $\lambda_3 \leq \alpha \leq \lambda_1$, the well-known
 estimate for the Rayleigh quotient.  We note that
$\lt^+(x,t)>0$ implies we have stretching of infinitesimal fluid
volume in two directions and compression in the other one
direction(planar stretching) at $(x,t)$, while $|\lt^-(x,t)|>0$
implies stretching in one direction and compressions in two
directions(linear stretching). The above estimate says that the
dominance competition between planar stretching and linear
stretching is an important mechanism controlling  the  growth/decay
in time of the $L^2$ norm of  vorticity.

 In order to state our next
theorem  we introduce some definitions. Given a differentiable
vector field $f=(f_1 ,f_2 ,f_3 )$ on $\Bbb T^3$, we denote by the
scalar field $\l_i (f)$, i=1,2,3, the eigenvalues of the deformation
tensor associated with $f$. Below we always assume the ordering, $
\lo (f)\geq \lt (f)\geq \ltt (f). $ We also fix $m>5/2$ below. We
recall that if $f\in \mathbb{H}^m _\s $, then $\lo (f)+\lt (f)+\ltt
(f)=0$, which is another representation of div $f=0$.

Let us begin with introduction of  admissible classes
$\mathcal{A}_\pm$ defined by
$$\mathcal{A}_+=\{ f\in \mathbb{H}^m _\s (\Bbb T^3)\, | \, \inf_{x\in \Bbb T^3}\lt (f)(x)
>0 \,\},$$
and
$$\mathcal{A}_-=\{ f\in \mathbb{H}^m_\s (\Bbb T^3)\, |
\sup_{x\in \Bbb T^3}\lt (f)(x)<0 \, \}.$$
 Physically $\mathcal{A}_+$ consists of solenoidal vector fields
 with planar stretching everywhere, while $\mathcal{A}_-$ consists
 of everywhere linear stretching vector fields. Although they do not represent
 real physical flows, they might be useful in the study of searching
 initial data leading to finite time singularity for the 3D Euler
 equations. Given $v_0 \in \mathbb{H}_\s ^m$, let $T_*(v_0)$ be the maximal time
of unique existence of solution in $\mathbb{H}_\s ^m$ for the system
(E).
 Let $S_t : \mathbb{H}^m_\s \to
\mathbb{H}^m_\s$ be the solution operator, mapping from initial data
to the
 solution $v(t)$.
Given $f\in \mathcal{A}_+$, we define the first zero touching time
of $\lt (f)$ as
$$ T(f)=
\inf\{ t\in (0, T_* (v_0)) \, | \, \mbox{$\exists x\in \Bbb T^3 $
such that $\lt (S_t f ) (x)<0 $}\}.
$$
Similarly for $f\in \mathcal{A}_-$, we define
$$ T(f)=
\inf\{ t\in (0, T_* (v_0)) \, | \, \mbox{$\exists x\in \Bbb T^3 $
such that $\lt (S_t f ) (x)>0 $}\}.
$$
The following theorem is actually an immediate corollary of Theorem
6.1, combined with the above definition of $\mathcal{A}_\pm$ and
$T(f)$. We just observe that for $v_0 \in \mathcal{A}_+ $(resp.
$\mathcal{A}_- $) we have $\lt^-=0, \lt^+=\lt$(resp. $\lt^+ =0,
\lt^-=\lt$) on $\O \times (0,T( v_0 ))$.
\begin{theorem}
Let $v_0\in \mathcal{A}_\pm$ be given. We set $\lo(x,t)\geq
\lt(x,t)\geq \ltt(x,t)$ as the eigenvalues of the deformation tensor
associated with $v(x,t)=(S_t v_0)(x)$ defined $t\in (0, T(v_0 ))$.
Then, for all $t\in (0, T(v_0 ))$ we have the following
estimates:\\
 (i) If $v_0 \in \mathcal{A}_+$, then
 $$
 \exp\left( \frac12\int_0 ^t
 \inf_{x\in\Bbb T^3} |\lt (x,s)| ds\right)\leq \frac{\|\o (t)\|_{L^2}}{ \|\o_0\|_{L^2}}
  \leq  \exp\left( \int_0 ^t
 \sup_{x\in\Bbb T^3} |\lt (x,s)| ds\right) .
 $$
(ii) If $v_0 \in \mathcal{A}_-$, then
 $$
 \exp\left( - \int_0 ^t
 \sup_{x\in\Bbb T^3} |\lt (x,s)| ds\right)\leq \frac{\|\o (t)\|_{L^2}}{ \|\o_0\|_{L^2}}
  \leq  \exp\left( -\frac12 \int_0 ^t
 \inf_{x\in\Bbb T^3} |\lt (x,s)| ds\right).
$$
\end{theorem}
(see \cite{cha8} for the proof) If we have the compatibility
conditions,
 \bqn
\inf_{x\in\Bbb T^3} |\lt (x,t)|&\simeq& \sup_{x\in\Bbb T^3} |\lt
(x,t)|\simeq g(t) \quad \forall t\in (0, T(v_0 )), \eqn
 which is the case for sufficiently small box $\Bbb T^3$, then we have
 \bqn
 \frac{\|\o (t)\|_{L^2}}{ \|\o_0\|_{L^2}} \simeq\left\{
 \aligned &\exp\left(\int_0 ^t g(s)ds \right) \quad \mbox{if}\quad v_0 \in
 \mathcal{A}_+\\
 &\exp\left(-\int_0 ^t g(s)ds \right) \quad \mbox{if}\quad v_0
\in
 \mathcal{A}_- \endaligned \right.
 \eqn
 for $t\in (0, T(v_0 ))$.
In particular, if we could find $v_0 \in \mathcal{A}_+$ such that
 $$
\inf_{x\in\Bbb T^3} |\lt (x,t)| \gtrsim
O\left(\frac{1}{t_*-t}\right)
 $$
 for time
interval near $t_*$, then such data would lead to singularity at
$t_*$.

As another application of the formula (\ref{spec1}) we have some
decay in time estimates for some ratio of eigenvalues(see
\cite{cha8} for the proof).

\begin{theorem}
Let $v_0 \in \mathcal{A}_\pm$ be given, and we set $\lo(x,t)\geq
\lt(x,t)\geq \ltt(x,t)$ as  in Theorem 3.1.
  We define
 $$
 \e (x,t)=\frac{|\lt (x,t) |}{\l (x,t)}\quad \forall (x,t)\in
\Bbb T^3 \times (0, T(v_0)),
  $$ where we set
$$\l (x,t)=\left\{\aligned \l _1 (x,t)
\quad \mbox{ if} \quad v_0 \in  \mathcal{A}_+\\
 -\l _3 (x,t) \quad \mbox{ if} \quad v_0 \in  \mathcal{A}_- .
\endaligned
\right.
$$
  Then, there exists a constant $C=C(v_0)$ such that
  $$
 \inf_{(x,s)\in \Bbb T^3\times (0,t)} \e (x,s)< \frac{C}{\sqrt{t}}\quad
 \forall t\in (0, T(v_0)).
 $$
 \end{theorem}
 Regarding the problem of searching  finite time
 blowing up solution,  the proof of
 the above theorem suggests the following:\\
 Given $\delta >0$, let us suppose we could find $v_0 \in
 \mathcal{A}_+$ such that for the associated solution $v(x,t)=(S_t v_0)(x)$
 the estimate
 \bb\label{spec3.3}
 \inf_{(x,s)\in \Bbb T^3 \times (0,t)}\e (x,s)\gtrsim
 O\left(\frac{1}{t^{\frac12+\delta}}\right),
 \ee
 holds true, for sufficiently
 large time $t$.
Then such $v_0$ will lead to the finite time singularity. In order
to check the  behavior (\ref{spec3.3}) for a given solution we need
a sharper and/or localized
 version of the equation (\ref{spec1}) for the dynamics of eigenvalues
 of the deformation tensor.

 \section{Conservation laws for  singular solutions }
 \setcounter{equation}{0}

For the smooth solutions of the Euler equations there are many
conserved quantities as described in Section 1 of this article. One
of the most important conserved quantities is the total kinetic
energy.  For nonsmooth(weak) solutions it is not at all obvious that
we still have energy conservation. Thus, there comes very
interesting question of how much smoothness we need to assume for
the solution to have energy conservation property. Regarding this
question L. Onsager conjectured that a H\"{o}lder continuous weak
solution with the H\"{o}lder exponent $1/3$ preserve the energy, and
this is sharp(\cite{ons}). Considering Kolmogorov's scaling argument
on the energy correlation in the homogeneous turbulence the exponent
$1/3$ is natural.  A sufficiency part of this conjecture is proved
in a positive direction by an ingenious argument due to
Constantin-E-Titi\cite{con7}, using a special Besov type of space
norm, $\mathcal{\dot{B}}^s_{3, \infty}$ with $s>1/3$(more precisely,
the Nikolskii space norm) for the velocity. See also \cite{caf2} for
related results in the magnetohydrodynamics. Remarkably enough
Shnirelman\cite{shn2} later constructed an example of weak solution
of 3D Euler equations, which does not preserve energy. The problem
of finding optimal regularity condition for a weak solution to have
conservation property can also be considered for the helicity. Since
the helicity is closely related to the topological invariants, e.g.
the knottedness of vortex tubes, the non-conservation of helicity is
directly related to the spontaneous apparition of singularity from
local smooth solutions, which is the main theme of this article. In
\cite{cha5} the author of this article obtained a sufficient
regularity condition for the helicity conservation, using the
function space $\mathcal{\dot{B}}^s_{\frac95, \infty}$, $s>1/3$, for
the vorticity. These results on the energy and the helicity are
recently refined in \cite{cha9}, using the Triebel-Lizorkin type of
spaces, $\mathcal{\dot{F}}^s_{p,q}$, and the Besov spaces
$\mathcal{\dot{B}}^s_{p,q}$(see Section 1 for the definitions) with
similar values for $s, p$, but allowing full range of values for
$q\in [1, \infty]$.

By a weak solution of $(E)$ in $\Bbb R^n \times (0, T)$ with initial
data $v_0$ we mean a vector field $v\in C([0, T); L^2_{loc} (\Bbb
R^n))$ satisfying the integral identity:
 \bq
 \lefteqn{-\int_0 ^T \intr v(x,t)\cdot \frac{\partial\phi (x,t)}{\partial
 t} dx dt -\intr v_0(x)\cdot \phi (x, 0 )dx}\hspace{.0in}\n \\
 &&\quad -\int_0 ^T \intr v(x,t)\otimes v(x,t) : \nabla \phi (x,t)
 dxdt\n \\
&&\qquad -\int_0 ^T \intr \mathrm{div }\, \phi (x,t ) p(x,t )dxdt
=0, \eq \bb
 \int_0 ^T \intr v(x,t)\cdot \nabla \psi (x,t)dxdt=0
\ee
 for every vector test function $\phi =(\phi_1 , \cdots
 ,\phi_n )\in C_0 ^\infty (\Bbb R^n\times [0, T))$, and for every scalar
 test function $\psi \in C_0 ^\infty (\Bbb R^n \times [0, T))$. Here we
 used the notation $(u\otimes v)_{ij}= u_i v_j$, and
 $A:B=\sum_{i,j=1} ^n A_{ij} B_{ij}$ for $n\times n$ matrices $A$ and $ B$.
In the case when we discuss the helicity conservation of the weak
solution  we impose further regularity for the
 vorticity, $\o (\cdot ,t)\in L^{\frac32}
 (\Bbb R^3)$ for almost every $t\in [0,T]$ in order to define the
 helicity for such weak solution. {\em Hereafter, we use the notation
${\dot{X}}^s_{p,q}$(resp. ${X}^s_{p,q}$) to represent
$\mathcal{\dot{F}}^s_{p,q}$(resp. $\mathcal{{F}}^s_{p,q}$) or
$\mathcal{\dot{B}}^s_{p,q}$(resp. $\mathcal{{B}}^s_{p,q}$)}. The
following is proved in \cite{cha9}.
\begin{theorem}
 Let $s>1/3$ and $q\in [2, \infty]$  be given. Suppose $ v$ is a weak
 solution of the  $n-$dimensional Euler equations
 with  $v\in C([0, T]; L^{2} (\Bbb R^n))\cap L^3 (0, T;
 \dot{X}^s_{3,q}(\Bbb R^n ))$. Then,
 the energy is preserved in time, namely
  \bb
  \intr |v(x,t)|^2 dx =\intr |v_0 (x)|^2 dx
  \ee
  for all $t\in [0, T)$.
  \end{theorem}
 When we restrict $q=\infty$,  the above theorem
 reduce to the one in \cite{con7}.  On the other hand, the
results for Triebel-Lizorkin type of space are completely new.

  \begin{theorem}
 Let $s>1/3$,  $q\in [2, \infty]$, and $r_1\in [2, \infty], r_2\in [1, \infty]$ be given, satisfying
 $2/r_1 +1/r_2 =1$. Suppose $ v$ is a weak
 solution of the 3-D Euler equations
 with $v\in C([0, T]; L^{2} (\Bbb R^3 ))\cap L^{r_1} (0, T;
\dot{ X}^s_{\frac92,q} (\Bbb R^3 ))$ and $\o \in L^{r_2} (0, T;
{\dot{X}}^s_{\frac95,q}(\Bbb R^3 ))$,
 where the curl operation is  in the sense of distribution. Then,
 the helicity is preserved in time, namely
  \bb\label{helicityy}
  \int_{\Bbb R^3} v(x,t)\cdot \o (x,t )dx =\int_{\Bbb R^3} v_0 (x)\cdot \o _0 (x)dx
  \ee
  for all $t\in [0, T)$.
  \end{theorem}
 Similarly to the case of Theorem 7.1, when  we restrict $q=\infty$,  the above theorem
 reduce to the one in \cite{cha5}.  The
results for the case of the Triebel-Lizorkin type of space, however,
is  new in \cite{cha9}.

As an application of  the above theorem we have the following
  estimate from below of the vorticity by a constant depending on
   the initial data
  for the weak solutions of the 3D Euler equations.
  We estimate the helicity,
\bqn
 \lefteqn{\int_{\Bbb R^3} v(x,t )\cdot \o (x,t )dx
 \leq \|v (\cdot , t)\|_{L^3} \|\o (\cdot , t
 )\|_{L^{\frac32}}}\hspace{.2in} \\
 &\leq &C \|\nabla v (\cdot , t)\|_{L^{\frac32}} \|\o (\cdot , t )\|_{L^{\frac32}}
 \leq C\|\o (\cdot , t )\|_{L^{\frac32}}^2,
 \eqn
 where we used the Sobolev inequality and the Calderon-Zygmund
 inequality. Combining this estimate with  (\ref{helicityy}), we obtain
 the following:
\begin{cor}
  Suppose $ v$ is a weak
 solution of the 3D Euler equations satisfying the conditions of
 Theorem 7.2. Then,
 we have the following estimate:
 $$
 \|\o (\cdot ,t)\|_{L^{\frac32}}^2 \geq C H_0, \quad \forall
 t\in [0, T)
 $$
 where $H_0 =\intt v_0(x)\cdot \o _0(x) dx$ is the initial
 helicity, and $C$ is an absolute constant.
\end{cor}

Next we are concerned on the $L^p$-norm conservation for the weak
solutions of (QG).  Let $p\in [2, \infty)$. By a weak solution of
$(QG)$ in $D\times (0, T)$ with initial data $v_0$ we mean a scalar
field $\theta \in C([0, T); L^p (\Bbb R^2)\cap L^{\frac{p}{p-1}}
(\Bbb R^2))$ satisfying the integral identity:
 \bb\label{qgg1}
-\int_0 ^T \int_{\Bbb R^2} \theta (x,t)\left[\frac
{\partial}{\partial
 t} +v \cdot \nabla\right] \phi (x,t)dx dt -\int_{\Bbb R^2}  \th_0(x) \phi (x, 0
 )dx=0
 \ee
 \bb
 \label{qgg2}
 v(x,t) =-\nao\int_{\Bbb R^2}
 \frac{ \th (y,t)}{|x-y|} dy
 \ee
 for every test function $\phi \in C_0 ^\infty (\Bbb R^2\times [0,
 T))$, where $\nao$ in (\ref{qgg2}) is in the sense of
 distribution. We note that contrary to the case of 3D Euler
 equations there is a global existence result for the weak
 solutions of (QG) for $p=2$ due to Resnick(\cite{res}). The following is proved in
 \cite{cha9}.

 \begin{theorem}
 Let $s>1/3$,  $p\in [2, \infty)$, $q\in [1, \infty]$, and
 $r_1\in [p, \infty], r_2\in [1, \infty]$
 be given, satisfying
 $p/r_1 +1/r_2 =1$. Suppose $ \th $ is a weak
 solution of (QG)
 with $\theta \in C([0, T]; L^{p} (\Bbb R^2)\cap L^{\frac{p}{p-1}}(\Bbb R^2))\cap L^{r_1} (0, T;
 X^s_{p+1,q}(\Bbb R^2))$ and $v \in L^{r_2} (0, T;
{\dot{X}}^s_{p+1,q} (\Bbb R^2))$. Then,
 the $L^p$ norm of $\theta(\cdot,t)$ is preserved, namely
  \bb
  \|\theta (t)\|_{L^p}=\|\theta_0\|_{L^p}
  \ee
  for all $t\in [0, T]$.
\end{theorem}


\begin{thebibliography}{1}

\bibitem{arn}V. I. Arnold and B. A. Khesin, {\it Topological Methods
in Hydrodynamics,} Springer-Verlag, (1998).
\bibitem{bab}A. Babin, A. Mahalov and B. Nicolaenko, {\it 3D
Navier-Stokes and Euler equations with initial data characterized by
uniformly large vorticity,} Indiana Univ. Math. J., {\bf 50}, no. 1,
(2001),  pp. 1-35.
\bibitem{bah}H. Bahouri and B. Dehman, {\it Remarques sur
l'apparition de singularit\'{e}s dans les \'{e}coulements
Eul\'{e}riens incompressibles \`{a} donn\'{e}es initiales
H\"{o}lderiennes,} J. Math. Pure Appl., {\bf 67}, (1994), pp.
335-346.
\bibitem{bea}J. T. Beale, T. Kato and A. Majda,  {\it Remarks on the
breakdown of smooth solutions for the 3-D Euler equations}, Comm.
Math. Phys., {\bf 94},  (1984), pp. 61-66.
\bibitem{bei1}H. Beir$\tilde{a}$o da Veiga, {\it Vorticity and
Smoothness in Incompressible Viscous Flows,} in ``Wave Phenomena and
Asymptotic Analysis", RIMS, Kokyuroku 1315 ,(2003), pp. 37-45.
\bibitem{bei2}H. Beir$\tilde{a}$o da Veiga and L. C. Berselli, {\it On
the Regularizing Effect of the Vorticity Direction in Incompressible
Viscous Flows,} Diff. Int. Eqns, {\bf 15}, No. 3, (2002), pp.
345-356.
\bibitem{bra}M. E. Brachet, D. Meiron, S. Orszag, B. Nickel, R. Morf
and U. Frisch, {\it Small-scale structure of the Taylor-Green
vortex,} J. Fluid. Mech., {\bf 130}, (1983), pp. 411-452.
\bibitem{bren}Y. Brenier, {\it Topics on hydrodynamics and area
preserving maps,} Handbook of mathematical fluid dynamics(S.
Friedlaner and D. Serre eds.), Vol. II, North-Holland, Amsterdam,
(2003), pp. 55-86.
 \bibitem{bre} H. Brezis and S.  Wainger, {\it A note on limiting
  cases of Sobolev embeddings and convolution
inequalities,} Comm. P.D.E., {\bf 5 }, no. 7, (1980), pp.773-789.
\bibitem{caff}L. Caffarelli and A. Vasseur, {\it Drift diffusion
equations with fractional diffusion and the quasi-geostrophic
equation,} arXiv:Math.AP/0608447, (2006).
\bibitem{caf1}R. E. Caflisch, {\it Singularity formation for complex
solutions of the 3D incompressible Euler equations,} Physica D, {\bf
67}, (1993), pp. 1-18.
\bibitem{caf2}R. E. Caflisch, I. Klapper and G. Steele, {\it Remarks
on singularities, dimension and energy dissipation for ideal
hydrodynamics and MHD,} Comm. Math. Phys., {\bf 184}, (1997), pp.
443-455.
\bibitem{can} J.R. Cannon and E. DiBenedetto, {\it The initial
 problem for the Boussinesq equations with data in $L^p$,}
 Lect. Note Math., no. 771, Springer, Berlin, (1980), pp. 129-144.
\bibitem{cha1}D. Chae, {\it On the Well-Posedness of the Euler
Equations in the Besov and Triebel-Lizorkin Spaces}, ``Tosio Kato's
Method and Principle for Evolution Equations in Mathematical
Physics", pp. 42-57, Yurinsha, Tokyo (2002)--A Proceddings of the
workshop  held at Hokkaido University, Japan on June 27-29, 2001.
\bibitem{cha2}D. Chae, {\it On the well-posedness of the Euler
equations in the Triebel-Lizorkin spaces,} Comm. Pure and Appl.
Math., {\bf 55}, no. 5, (2002), pp. 654-678.
\bibitem{cha2a}D. Chae, {\it The quasi-geostrophic
equation in the Triebel-Lizorkin spaces,} Nonlinearity, {\bf 16},
(2003), pp. 479-495.
\bibitem{cha3} D. Chae, {\it On the Euler Equations in the Critical
 Triebel-Lizorkin Spaces},
  Arch. Rational Mech. Anal., {\bf 170}, no. 3, (2003), pp.185-210.
\bibitem{cha4} D. Chae, {\it Remarks on the blow-up of the Euler
 equations and the related equations,} Comm. Math. Phys., {\bf 245}, no. 3,
 (2003), 539-550.
\bibitem{cha5} D. Chae,
{\it Remarks on the helicity of the 3-D incompressible Euler
equations,} Comm. Math. Phys., {\bf 240}, (2003), pp. 501-507.
\bibitem{cha6}D. Chae, {\it  Local Existence and Blow-up Criterion for the
Euler Equations in the Besov Spaces}, Asymp. Anal., {\bf 38}, no.
3-4, (2004), pp. 339-358: a printed version of RIM-GARC(Seoul
National University, Korea) preprint no., {\bf 8} (June, 2001).
\bibitem{cha6a}D. Chae, {\it On the Dual Systems to the Euler and
the Navier-Stokes Equations in $\Bbb R^3$,} Proc. Roy. Soc. London,
Ser. A, {\bf 460}, no. 2044, (2004), pp. 1153-1168.
\bibitem{cha7}D. Chae, {\it Remarks on the blow-up criterion of the 3D Euler
 equations,} Nonlinearity, {\bf 18}, (2005), pp. 1021-1029.
 \bibitem{cha8}D. Chae, {\it On the spectral dynamics of the
 deformation tensor and new a priori estimates for the 3D Euler
 equations,} Comm. Math. Phys., {\bf 263}, (2006), pp. 789-801.
 \bibitem{cha9}D. Chae, {\it On the conserved quantities for the weak
 solutions of the Euler equations and the quasi-geostrophic
 equations,} Comm. Math. Phys., {\bf 266}, (2006), pp. 197-210.
 \bibitem{cha10}D. Chae, {\it On the continuation principles for  the Euler equations and the
quasi-geostrophic equation,} J. Diff. Eqns, {\bf 227}, (2006), pp.
640-651.
\bibitem{cha10a}D. Chae, {\it On the Regularity Conditions for the Dissipative
Quasi-geostrophic Equations,} SIAM J. Math. Anal., {\bf 37}, no. 5,
(2006), pp. 1649-1656.
\bibitem{cha11}D. Chae, {\it Global regularity for the 2D Boussinesq equations with partial
viscosity terms,} Advances  in Math., {\bf 203}, (2006), pp.
497-513.
  \bibitem{cha12}D. Chae, {\it On the Lagragian dynamics for the 3D
  incompressible Euler equations,} Comm. Math. Phys., {\bf 269}.
  (2006), pp. 557-569.
  \bibitem{cha12a}D. Chae, {\it Incompressiblle Euler Equations: Mathematical Theory,}
  Encyclopedia of Mathematical Physics,  Elsevier Science Ltd., {\bf 3},
  (2006), pp. 10-17.
   \bibitem{cha15}D. Chae, {\it On the finite
time singularities of the 3D incompressible Euler equations}, Comm.
Pure Appl. Math., {\bf 60}, no. 4, (2007), pp. 597-617.
 \bibitem{cha17}D. Chae, {\it Notes on
perturbations of the 3D Euler equations,} Nonlinearity, {\bf 20},
(2007), pp. 517-522
\bibitem{cha13} D. Chae, {\it Nonexistence of
self-similar singularities for the 3D incompressible Euler
equations,} Comm. Math. Phys., in press.
   \bibitem{cha14} D. Chae, {\it Nonexistence of
asymptotically self-similar singularities in the  Euler and the
 Navier-Stokes equations,} Math. Ann., in press.
\bibitem{cha16}D. Chae, {\it On the deformations of the
incompressible  Euler equations,} Math. Z., in press.
\bibitem{cha18}D. Chae, {\it On the Regularity Conditions for the Navier-Stokes
and the Related Equations,} Revista Mat. Iberoamericana, {\bf 23},
no. 1, (2007), pp. 373-386.
\bibitem{cha19}D. Chae, A. C\'{o}rdoba, D. C\'{o}rdoba and M. Fontelos, {\it Finite time
singularities in a 1D model of the quasi-geostrophic equation,}
Advances in Math., {\bf 194}, (2005), pp. 203-223.
\bibitem{cha19a}D. Chae and P. Dubovskii, {\it Functional and
measure-valued solutions of the Euler equations for flows of
incompressible fluids,}  Arch. Rational Mech. Anal., {\bf 129}, no.
4, (1995),  pp. 385-396.
\bibitem{cha22}D. Chae and O. Yu. Imanivilov, {\it Generic solvability
of the axisymmetric 3D Euler equations and 2D Boussunesq equations,}
J. Diff.  Eqns.,  {\bf 156}, no. 1, (1999),  pp. 1-17.
\bibitem{cha20}D. Chae,  K. Kang and J. Lee, {\it On the interior regularity
 of suitable weak solutions of the
Navier-Stokes equations,}  Comm. P.D.E., in press.
\bibitem{cha21}D. Chae and N. Kim, {\it On the breakdown of
axisymmetric smooth solutions for the 3-D Euler equations, } Comm.
Math. Phys., {\bf 178}, (1996), pp. 391-398.
\bibitem{cha24}D. Chae, S.-K. Kim, and H. -S. Nam, {\it Local
 existence and blow-up criterion of H\"{o}lder continuous
 solutions of the Boussinesq equations,} Nagoya Math. J., {\bf
 155}, (1999), pp. 55-80.
  \bibitem{cha24a}D. Chae and J. Lee, {\it Global well-posedness in the super-critical
dissipative quasi-geostrophic equations,} Comm. Math. Phys. , {\bf
233}, Issue 2, (2003), pp. 297-311.
\bibitem{cha23} D. Chae, H. -S. Nam, {\it Local existence and blow-up
 criterion for the Boussinesq equations,}
 Proc. Roy. Soc. Edinburgh, Sect. A, {\bf 127},  no. 5, (1997),
 pp. 935-946.
\bibitem{che1}J. Y. Chemin, {\it
R\'{e}gularit\'{e} des trajectoires des particules d'un fluide
incompressible remplissant l'espace,} J. Math. Pures Appl., {\bf
71}, (5) , (1992), pp. 407-417.
 \bibitem{che2}J. Y. Chemin, {\it Perfect incompressible fluids,}
 Clarendon Press, Oxford, (1998).
 \bibitem{cho}A. J. Chorin, {\it The evolution of a turbulent
 vortex,} Comm. Math. Phys., {\bf 83}, (1982), pp. 517-535.
 %\bibitem{cle}A. Clebsch, {\it \"{U}ber eine allgemeine
 %transformation hydrodynamischen gleichungen,} Crelle, liv., (1857)
 %and (1858).
 \bibitem{con1}P. Constantin, {\it Note on loss of regularity for
 solutions of the 3D incompressible and related equations,} Comm.
 Math. Phys., {\bf 106}, (1986), pp. 311-325.
\bibitem{con2}P. Constantin, {\it Geometric Statistics in
Turbulence}, SIAM Rev.,{\bf 36}, (1994), pp. 73-98.
\bibitem{con3}P. Constantin,  {\it A few results and open problems
regarding incompressible fluids,} Notices Amer. Math. Soc., {\bf
42}, no. 6, (1995), pp. 658-663.
\bibitem{con4}P. Constantin, {\it Absence of proper nondegenerate
generalized self-similar singularities,} J. Stat. Phys., {\bf 93},
no. 3/4, (1998), pp. 777-786.
\bibitem{con5}P. Constantin,  {\it An Eulerian-Lagrangian approach
for incompressible fluids: local theory,} Journal of AMS, {\bf 14},
(2001), pp. 263-278.
\bibitem{con6}P. Constantin, {\it Euler equations, Navier-Stokes
 equations and Turbulence,} Mathematical Foundation of Turbulent
 Viscous Flows, Lecture Notes in Mathematics, no. 1871, (2006), pp.
 1-43.
 \bibitem{con6a}P. Constantin, D. Cordoba and J. Wu, {\it On the
 critical dissipative quasi-geostrophic equation,} Indiana Univ.
 Math. J., {\bf 50}, (2001), pp. 97-107.
 \bibitem{con7}P. Constantin, W. E., and E. S. Titi, {\it Onsager's
 conjecture on the energy conservation for solutions of Euler's
 equation,} Comm. Math. Phys., {\bf 165}, (1994), pp. 207-207.
\bibitem{con8}P. Constantin and C. Fefferman, {\it Direction of
Vorticity and the Problem of Global Regularity for the Navier-Stokes
Equations,} Indiana Univ. Math. J., {\bf 42}, (1993), pp. 775-789.
\bibitem{con9}P. Constantin, C. Fefferman and A. Majda,
 {\it Geometric constraints on potential singularity formulation in the
 3-D Euler equations}, Comm. P.D.E., {\bf 21}, (3-4), (1996), pp.
 559-571.
 \bibitem{con10}P. Constantin and C. Foias, {\it Navier-Stokes Equations,}
 Chicago Lectures in Mathematics Series, Univ. Chicago Press (1988).
 \bibitem{con11}P. Constantin, P. Lax and A. Majda, {\it A simple
 one-dimensional model for the three dimensional vorticity
 equation,} Comm. Pure Appl. Math., {\bf 38}, (1985), pp. 715-724.
 \bibitem{con12}P. Constantin, A. Majda and E. Tabak, {\it Formation
 of strong fronts in the 2-d quasi-geostrophic thermal active
 scalar,} Nonlinearity, {\bf 7}, (1994), pp. 1495-1533.
 \bibitem{con12a}P. Constantin and J. Wu, {\it Regularity of
 H\"{o}lder continuous solutions of the supercritical
 quasi-geostrophic equation,} arXiv:math.AP/0701592, (2007).
\bibitem{con12b}P. Constantin and J. Wu, {\it H\"{o}lder continuity
of solutions of supercritical dissipative hydrodynamic transport
equations,} arXiv:math.AP/0701594, (2007).
\bibitem{cor3}A. C\'{o}rdoba and D. C\'{o}rdoba, {\it A Maximum
Principle Applied to Quasi-geostrophic Equations,} Comm. Math. Phys.
{\bf 249}, no. 3, (2004), pp. 511-528.
\bibitem{cor8}A. C\'{o}rdoba, D. C\'{o}rdoba and M. A. Fontelos,
{\it Formation of singularities for a transport equation with
nonlocal velocity,} Ann. of Math. (2), {\bf 162}, no. 3, (2005), pp.
1377-1389.
 \bibitem{cor1}D. C\'{o}rdoba, {\it On the geometry of solutions of the
 quasi-geostrophic and Euler equations,} Proc. Natl. Acad. Sci.,
 {\bf 94}, (1997), pp. 12769-12770.
 \bibitem{cor2} D. C\'{o}rdoba, {\it Nonexistence of simple hyperbolic blow-up for the
quasi-geostrophic equation,} Ann. of Math., {\bf 148}, (1998), pp.
1135-1152.
\bibitem{cor4}D. C\'{o}rdoba and C. Fefferman, {\it On the collapse of
tubes carried by 3D incompressible flows,} Comm. Math. Phys., {\bf
222}, (2), (2001), pp. 293-298.
\bibitem{cor4a}D. C\'{o}rdoba and C. Fefferman, {\it Potato chip
singularities of 3D flows,} SIAM J. Math. Anal., {\bf 33}, (2001),
pp. 786-789.
\bibitem{cor5} D. C\'{o}rdoba and C. Fefferman, {\it Growth of solutions for QG and 2D
Euler equations,} Journal Amer. Math. Soc., {\bf 15}, no. 3 (2002),
pp.665-670.
\bibitem{cor6}D. C\'{o}rdoba and C. Fefferman, {\it Scalars convected
by a two-dimensional incompressible flow,} Comm. Pure Appl. Math.,
{\bf 55}, no. 2, (2002), pp. 255-260.
\bibitem{cor7}D. C\'{o}rdoba, C. Fefferman and R. de la Llave, {\it On
Squirt Singularities in Hydrodynamics,} SIAM J. Math. Anal., {\bf
36}, no.1, (2004), pp. 204-213.
%\bibitem{dan} R. Danchin, {\it Local theory in critical
%spaces for compressible viscous and heat-conductive gases,} Comm.
%PDE, {\bf 26}, (2001), pp. 1183-1233.
\bibitem{den1}J. Deng, T. Y. Hou and X. Yu, {\it Geometric and Nonblowup
of 3D Incompressible Euler Flow,} Comm. P.D.E, {\bf 30}, (2005), pp.
225-243.
\bibitem{den2}J. Deng, T. Y. Hou and X. Yu, {\it Improved geometric
conditions for non-blow upof the 3D incompressible Euler equations,}
Comm. P.D. E., {\bf 31}, no. 1-3, (2006), pp. 293-306.
\bibitem{din}E. I. Dinaburg, V. S. Posvyanskii and Ya. G. Sinai,{\it  On
some approximations of the Quasi-geostrophic equation,} Geometric
methods in dynamics. II. Asterisque, {\bf xvii}, no. 287, (2003),
pp. 19-32.
\bibitem{e}W. E. and C. Shu, {\it Small scale structures un
Boussinesq convection,}  Phys. Fluids, {\bf 6}, (1994), pp. 48-54.
\bibitem{ebi}D. G. Ebin and J. E. Marsden, {\it Groups of
diffeomorphisms and the motion of an incompressible fluid,} Ann.
Math., {\bf 92}, (1970), pp. 102-163.
\bibitem{eng} H. Engler, {\it An alternative proof of the
Brezis-Wainger inequality,}  Comm. P.D.E., {\bf 14 }, no. 4, (1989),
pp.541-544.
\bibitem{eul}L. Euler, {\it Principes g\'{e}n\'{e}raux du mouvement des fluides,}
M\'{e}moires de l'acad\'{e}mie des sciences de Berlin, {\bf 11},
(1755), pp. 274-315.
\bibitem{eyi}G. Eyink, {\it Energy dissipation without viscosity in
ideal hydrodynamics, I. Fourier anaysis and local energy trandfer,}
Physica D, {\bf 78}, no. 3-4, (1994), pp. 222-240.
\bibitem{fri1}S. Friedlander, {\it Lectures on Stability and
Instability of an Ideal Fluid,} IAS/Park City Mathematics Series
{\bf 5}, `Hyperbolic Equations and Frequency Interactions' Edited by
L. Caffarelli and W. E, AMS/IAS, (1998).
\bibitem{fri2}S. Friedlander, N. Pavlovi\'{c}, {\it Blow-up in a three
dimensional vector model for the Euler equations,} Comm. Pure Appl.
Math., {\bf 42}, (2004), pp. 705-725.
\bibitem{fris1}U. Frisch, {\it Turblence}, Cambridge University
Press, (1995).
\bibitem{fris2}U. Frisch, T. Matsumoto and J. Bec, {\it
Singularities of Euler Flow? Not Out of Blue!}, J. Stat. Phys., {\bf
113}, no. 5-6, (2003), pp. 761-781.
\bibitem{gal}B. Galanti, J.D. Gibbon and M. Heritage, {\it Vorticity
alignment results for the three-dimensional Euler and Navier-Stokes
equations,} Nonlinearity, {\bf 10}, (1997), pp. 1675-1694.
\bibitem{gald}G. Galdi, {\it An Introduction to the mathematical
theory of Navier-Stokes equations, I,II}, Springer-Verlag, (1994).
\bibitem{gib1}J. D. Gibbon, {\it A quaternionic structure in the
three-dimensional Euler and ideal magneto-hydrodynamics equations,}
Physica D., {\bf 166}, (2002), pp. 17-28.
\bibitem{gib2}J. D. Gibbon, {\it Ortho-normal quaternion frames,
Lagrangian evolution equations and the three-dimensional Euler
equations,} arXiv:math-ph/0610004, (2006).
\bibitem{gib3}J. D. Gibbon, D. D. Holm, R. M. Kerr and I. Roulstone,
{\it Quaternions and particle dynamics in Euler fluid flow,}
Nonlinearity, {\bf 19}, (2006), pp. 1969-1983.
 \bibitem{gig}Y. Giga and R. V. Kohn, {\it Asymptotically
 Self-Similar Blow-up of Semilinear Heat Equations,} Comm. Pure
 Appl. Math., {\bf 38}, (1985), pp. 297-319.
 \bibitem{gra1}R. Grauer and T. Sideris, {\it Numerical computation
of three dimensional incompressible ideal fluids with swirl,} Phys.
Rev. Lett., {\bf 67}, (1991), pp. 3511-3514.
\bibitem{gra2}R. Grauer and T. Sideris, {\it  Finite time singularities
in ideal fluids with swirl,} Physica D, {\bf 88}, no. 2, (1995),
pp.116-132.
 \bibitem{gre}J. M. Greene and R. B. Pelz, {\it Stability of postulated,
self-similar, hydrodynamic blowup solutions}, Phys. Rev. E, {\bf
62}, no. 6, pp. 7982-7986.
\bibitem{hou1} T. Y. Hou and C. Li, {\it Global well-posedness of
the viscous Boussinesq equations,} Discrete and Continuous Dynamical
System, {\bf 12}, No. 1, (2005), pp. 1-12.
\bibitem{hou2}T. Y. Hou and R. Li, {\it Nonexistence of local
self-similar blow-up for the 3D incompressible Navier-Stokes
equations,} arXiv-preprint, math.AP/0603126.
\bibitem{hou3} T. Y. Hou and R. Li, {\it Dynamic depletion of
vortex stretching and non-blowup of the 3-D incompressible Euler
equations}, to appear in J. Nonlinear Sciences.
 \bibitem{jaw}B. Jawerth, {\it Some observations on Besov and
 Lizorkin-Triebel Spaces}, Math. Sacand., {\bf 40}, (1977), pp.
 94-104.
\bibitem{kat1}T. Kato, {\it On classical solutions of the two
dimensional nonstationary Euler equations,} Arch. Rat. Mech. Anal.,
{\bf 25}, (1967), pp. 188-200.
\bibitem{kat2}T. Kato, {\it Nonstationary flows of viscous and ideal
 fluids in $\Bbb R^3$}, J. Funct. Anal.,{\bf 9}, (1972),  pp.
 296-305.
%\bibitem{kat3} T. Kato, {\it Strong $L^p$ solutions of the
%Navier-Stokes equations in ${\mathbb R}^m$ with applications to
%weak solutions,} Math. Z., {\bf 187}, (1984), pp. 471-480.
\bibitem{kat3}T. Kato and G. Ponce, {\it Well posedness of the
 Euler and Navier-Stokes equations in Lebesgue spaces $L^s_p (\Bbb
 R^2 )$}, Revista Ibero-Americana, {\bf 2},  (1986), pp. 73-88.
\bibitem{kat4}T. Kato and G. Ponce, {\it Commutator estimates and
the Euler and Navier-Stokes equations,} Comm. Pure Appl. Math., {\bf
41}, (1988), pp. 891-907.
\bibitem{katz}N. Katz and N. Pavlovic, {\it Finite time blowup
for a dyadic model of the Euler equations,} Trans. AMS, {\bf 357},
no. 2, (2005), pp. 695-708.
\bibitem{ker1}R. M. Kerr, {\it Evidence for a singularity of the
3-dimensional, incompressible Euler equations,} Phys. Fluids A, {\bf
5}, (1993), pp. 1725-1746.
\bibitem{ker2} R. M. Kerr, {\it Computational Euler history},
arXiv:physics/0607148, (2006).
\bibitem{kis}A. Kiselev, F. Nazarov and A. Volberg, {\it Global
well-posedness for the critical 2D dissipative quasi-geostrophic
equation}, arXiv:Math.AP/0604185, (2006).
\bibitem{kla}S. Klainerman and A. Majda, {\it Singular limits of
quasilinear hyperbolic systems with large parameters and the
incompressible limit of compressible fluids,} Comm. Pure Appl. Math.
{\bf 34}, (1981), pp. 481-524.
 \bibitem{koz1}H. Kozono and Y. Taniuchi, {\it Limiting case of the
Sobolev inequality in BMO, with applications to the Euler
equations,} Comm. Math. Phys., {\bf 214}, (2000), pp. 191-200.
\bibitem{koz2}H. Kozono, T. Ogawa, and T. Taniuchi, {\it
The critical Sobolev inequalities in Besov spaces and regularity
criterion to some semilinear evolution equations,} Math Z., {\bf
242}(2), (2002), pp. 251-278.
\bibitem{lad}O. A. Ladyzenskaya, {\it The mathematical theory of
viscous incompressible flow,} Gordon and Breach, (1969).
\bibitem{lam}H. Lamb, {\it Hydrodynamics}, Cambridge Univ. Press,
(1932).
\bibitem{lem}P. G. Lemari\'{e}-Rieusset, {\it Recent developments in
the Navier-Stokes problem,}   research notes in mathematics series,
{\bf 431}, Chapman \& Hall/CRC , (2002).
\bibitem{ler}J. Leray, {\it Essai sur le mouvement d'un fluide
visqueux emplissant l'espace,} Acta Math., {\bf 63}, (1934), pp.
193-248.
\bibitem{lic}L. Lichtenstein, {\it \"{U}ber einige
Existenzprobleme der Hydrodynamik homogener
unzusammendr\"{u}ckbarer, reibunglosser Fl\"{u}ssikeiten und die
Helmholtzschen Wirbelsalitze,} Math. Zeit., {\bf 23}, (1925), pp.
89-154; {\bf 26}, (1927), pp. 193-323, pp. 387-415 ; {\bf 32},
(1930), pp. 608-725.
\bibitem{lio}P. L. Lions, {\it Mathematical Topics in Fluid
Mechanics}, {\bf Vol 1.} Incompressible Models, Oxford University
Press, (1996).
\bibitem{liu} H. Liu and E. Tadmor,
{\it Spectral Dyanamics of the Velocity Gradient Field in Restricted
Flows,} Comm. Math. Phys., {\bf 228}, (2002), pp. 435-466.
\bibitem{maj1}A. Majda, {\it Compressible Fluid Flow and Systems of
Conservation Laws in Several Space Variables, } Appl. Math. Sci.,
{\bf 53}, Springer, (1984).
\bibitem{maj2}A. Majda, {\it Vorticity and the mathematical theory
of incompressible fluid flow,} Comm. Pure Appl. Math., {\bf 39},
(1986), pp. 187-220.
\bibitem{maj3}A. Majda and A. Bertozzi, {\it Vorticity and
Incompressible Flow,} Cambridge Univ. Press. (2002).
\bibitem{maj4} A. Majda, {\it Introduction to PDEs and Waves for
the Atmosphere and Ocean,} Courant Lecture Notes in Mathematics, no.
9, AMS/CIMS, (2003).
\bibitem{mar}C. Marchioro and M. Pulvirenti, {\it Mathematical
Theory of Incompressible Nonviscous Fluids,} Springer-Verlag,
(1994).
\bibitem{mil}J. R. Miller, M. O'Leary and M. Schonbek, {\it
Nonexistence of singular pseudo-self-similar solutions of the
Navier-Stokes system,} Math. Ann., {\bf 319}, (2001), no. 4, pp.
809-815.
\bibitem{mof}H. K. Moffatt, {\it Some remarks on topological fluid
mechanics,} in An Introduction to the Geometry and Topology of Fluid
Flows, R. L. Ricca, ed., Kluwer Academic Publishers, Dordrecht, The
Netherlands, (2001), pp. 3-10.
\bibitem{mor}A. Morlet, {\it Further properties of a continuum of model
equations with globally defined flux,} J. Math. Anal.  Appl.,
\textbf{221}, (1998), pp. 132-160.
\bibitem{nec}J. Necas, M.  Ruzicka and V. Sverak, {\it On Leray's
self-similar solutions of the Navier-Stokes equations,} Acta Math.
{\bf 176}, no. 2, (1996), pp. 283-294.
\bibitem{neu}J.  Neustupa and
P. Penel, {\it Regularity of a weak solution to the Navier-Stokes
equation in dependence on eigenvalues and eigenvectors of the rate
of deformation tensor,} Progr. Non.  Diff. Eqns. Appl., {\bf 61},
Birkhauser, Basel, (2005), pp. 197-212.
\bibitem{ohk}K. Ohkitani and M. Yamada, {\it Inviscid and inviscid-limit behavior
of a surface quasi-geostrophic flow,} Phys. Fluids, {\bf 9}, (1997),
pp. 876-882.
\bibitem{ons}L. Onsager, {\it statistical hydrodynamics,} Nuovo
Cimento Suppl., {\bf 6}, (1949), pp. 279-287.
\bibitem{oka}H. Okamoto and K. Ohkitani, {\it On the Role of the Convection Term in
the Equations of Motion of Incompressible Fluid,} J. Phys. Soc. of
Japan, {\bf 74}, no. 10, (2005), pp.2737-2742.
\bibitem{pel}R. Pelz, {\it Symmetry and hydrodynamic blow-up
problem}, J. Fluid Mech., {\bf 444}, pp. 299-320.
\bibitem{pon}G. Ponce, {\it Remarks on a paper by J. T. Beale, T.
Kato and A. Majda,} Comm. Math. Phys., {\bf 98}, (1985), pp.
349-353.
\bibitem{res}S. Resnick, {\it Dynamical problems in nonlinear advective partial
differential equations,}  Ph.D. Thesis,  University of Chicago,
Chicago, (1995).
\bibitem{run}T. Runst and W. Sickel, {\it Sobolev Spaces of
Fractional Order, Nemytskij Operators, and Nonlinear Partial
Differential Equations,} Walter de Gruyter, Berlin, (1996).
%\bibitem{ruz}A. Ruzmaikina and Z. Gruji\'{c}, {\it On Depletion of
%the Vortex-strectching Term in the 3D Navier-Stokes Equations,}
%Comm. Math. Phys., {\bf 247}, (2004), pp. 601-611.
 \bibitem{sak1}T. Sakajo, {\it Blow-up solutions of the
 Constantin-Lax-Majda equation with a
generalized viscosity term,} J. Math. Sci. Univ. Tokyo, {\bf 10},
no. 1, (2003),  pp. 187-207.
 \bibitem{sak2}T. Sakajo,{\it On global solutions for the
 Constantin-Lax-Majda equation
with a generalized viscosity term,} Nonlinearity, {\bf 16}, no. 4,
(2003), pp. 1319-1328.
\bibitem{sch}V. Scheffer, {\it An inviscid flow with compact support
in space-time,} J. Geom. Anal, {\bf 3}, no. 4, (1993), pp. 343-401.
\bibitem{scho}S. Schochet, {\it Explicit solutions of the viscous model
vorticity equation,} Comm. Pure Appl. Math., {\bf  39}, no. 4,
(1986), pp. 531-537.
\bibitem{shn1}A. Shnirelman, {\it On the nonuniqueness of weak
solution of the Euler equations,} Comm. Pure Appl. Math., {\bf L},
(1997), pp. 1261-1286.
\bibitem{shn2}A. Shnirelman, {\it Weak solutions with decreasing
energy of incompressible Euler equations,} Comm. Math. Phys., {\bf
210}, (2000), pp. 541-603.
\bibitem{ste1}{ E. M. Stein}, {\it Singular Integrals and
Differentiability Properties of Functions,} Princeton Univ. Press,
Princeton, NJ, (1970).
\bibitem{ste2}{ E. M. Stein}, {\it Harmonic Analysis, Real Variable Methods,
Orthogonality, and Oscillatory Integrals,} Princeton Univ. Press,
Princeton, NJ, (1993).
\bibitem{tad}E. Tadmor, {\it On a new scale of regularity spaces
with applications to Euler's equations,} Nonlinearity, {\bf 14},
(2001), pp. 513-532.
\bibitem{tan} Y. Taniuchi, {\it A note on the blow-up criterion
for the inviscid 2-D Boussinesq equations,} Lecture Notes in Pure
and Applied Math., {\bf 223}, ``the Navier-Stokes equations: theory
and numerical methods", edited by R. Salvi.,(2002), pp. 131-140.
\bibitem{tay}M. E. Taylor, {\it Tools for PDE}, Mathematical Surveys
and Monographs, {\bf  81}, Amer. Math. Soc. (2000).
\bibitem{tem1}R. Temam, {\it On the Euler equations of
incompressible flows,} J. Funct.  Anal., {\bf 20}, (1975), pp.
32-43.
 \bibitem{tem2}R. Temam, {\it  Local existence of  solutions of the Euler
equations of incompressible perfect fluids,} Lecture Notes in
Mathematics {\bf 565}, Berlin, Heidelberg, New York, Springer,
(1976), pp. 184-195.
\bibitem{tem3}R. Temam, {\it Navier-Stokes equations,} 2nd ed.,
 North-Holland, Amsterdam,
(1986).
\bibitem{tri}H. Triebel,  {\it Theory of Function Spaces},
Birk\"{a}user Verlag, Boston, (1983).
\bibitem{tsa}T-P. Tsai, {\it  On Leray's self-similar solutions of the Navier-Stokes
equations satisfying local energy estimates,} Arch. Rational Mech.
Anal., {\bf 143}, no. 1, (1998), pp. 29-51.
\bibitem{vis1}M. Vishik, {\it Hydrodynamics in Besov spaces,}
Arch. Rational Mech. Anal, {\bf 145}, (1998), pp. 197-214.
\bibitem{vis2}M. Vishik, {\it Incompressible flows of an ideal
fluid with vorticity in borderline spaces of Besov type,} Ann.
Scient. \'{E}c. Norm. Sup., $4^e$ s\'{e}rie, t. 32, (1999), pp.
769-812.
\bibitem{wu1}J. Wu, {\it Inviscid limits and regularity estimates for the
solutions of the 2-D dissipative Quasi-geostrophic equations,}
Indiana Univ. Math. J., {\bf 46}, no. 4 (1997), pp. 1113-1124.
\bibitem{wu2} J. Wu, {\it Dissipative quasi-geostrophic equations with
$L^p$ data, } Electro. J. Dff.  Eqns,  {\bf 56} (2001), pp. 1-13.
\bibitem{wu3} J. Wu, {\it The  quasi-geostrophic equations and its two
regularizations,}  Comm. P.D.E. {\bf 27} $no.$ 5-6 (2002), pp.
1161-1181.
\bibitem{yu}X. Yu, {\it Localized Non-blow-up Conditions for 3D
Incompressible Euler Flows and Related Equations,} Ph.D Thesis,
California Institute of Technology, (2005).
\bibitem{yud1}V. I. Yudovich, {\it Non-stationary flow of an ideal
incompressible fluid,} Akademiya Nauk SSSR. Zhurnal Vychislitelnol
Matematiki I Matematicheskoi Fiziki, {\bf 3}, (1963), pp. 1032-1066.
\bibitem{yud2}V. I. Yudovich, {\it Uniqueness theorem for the basic
nonstationary problem in the dynamics of an ideal incompressible
fluid,} Math. Res. Lett., {\bf 2}, (1995), pp. 27-38.
  \end{thebibliography}
\end{document}